\documentclass[11pt]{article}
\usepackage{geometry}
\usepackage{float}
\usepackage{latexsym,amsmath,amsfonts,amscd,amssymb,subfigure}
\usepackage{graphicx}
\usepackage{color}
\usepackage{changebar}
\usepackage{bm}
\numberwithin{equation}{section}

\begin{document}

\vspace{.5in}

\begin{center}

{\Large \bf A fully convergent fixed-point fast sweeping method with the WENO-JS local solver for steady state of hyperbolic conservation laws}

\end{center}

\vspace{.15in}

\centerline{
Liang Li\footnote{
School of Mathematics and Statistics, Huang Huai University, Zhumadian, Henan 463000, P.R. China. E-mail: liliangnuaa@163.com. Research was supported by Natural Science Foundation of Henan Province (252300423517) and NSFC 12501548.},
Jun Zhu\footnote{
College of Science and Key Laboratory of Mathematical Modelling and High Performance Computing of Air Vehicles (NUAA), MIIT, Nanjing University of Aeronautics and Astronautics, Nanjing, Jiangsu 210016, P.R.
China. E-mail: zhujun@nuaa.edu.cn. Research was supported by NSFC 12472292 and Science Challenge Project grant TZ2025007.}
and
Yong-Tao Zhang\footnote{Department of Applied and Computational
Mathematics and Statistics, University of Notre Dame, Notre Dame,
IN 46556, USA. E-mail: yzhang10@nd.edu. Research was partially supported by Simons Foundation MPS-TSM-00007854.}$\textsuperscript{,}
\renewcommand*{\thefootnote}{\fnsymbol{footnote}}
\setcounter{footnote}{0}\footnote{Corresponding author.}$
}
\baselineskip=1.7pc

\begin{abstract}
The fixed-point fast sweeping methods with weighted essentially non-oscillatory (WENO) local solvers are a class of efficient and high-order accuracy numerical methods for solving steady-state solutions of hyperbolic conservation laws. However, with the classical WENO-JS local solver, the iteration residue of high-order fixed-point fast sweeping scheme often hangs at a truncation error level, or even higher error levels, instead of settling down
to round-off errors / machine zero. This is a typical issue when the high-order WENO-JS schemes are applied for solving steady-state problems of hyperbolic conservation laws. To achieve the full convergence in a fast sweeping method, the fixed-point fast sweeping methods with non-traditional WENO local solvers based on unequal-sized substencils (e.g., the multi-resolution WENO schemes) were designed. However, the WENO schemes based on unequal-sized stencils are more complex and in general more expensive in computational costs than the classical WENO-JS schemes. 
In this paper, we go back to the classical WENO-JS local solver and develop a new fully convergent fifth-order fixed-point fast sweeping method for solving steady-state problems of hyperbolic conservation laws. Based on recent studies on the nonlinear weighting process around discontinuities of solution, we apply the technique of frozen weights for avoiding unnecessary adjustment of nonlinear weights in the WENO-JS local solver, which freezes the nonlinear weights once the residual sequence in the fast sweeping iterations has stabilized. Different from the existing work on frozen weights, we design a simple and robust approach to judge stabilization of iteration residues and determine the iteration step when the nonlinear weights are frozen in the fast sweeping method. Extensive numerical experiments on a wide range of challenging two-dimensional steady-state problems demonstrate that, unlike the previous fast sweeping method with the fifth-order WENO-JS local solver, the proposed new scheme consistently drives the iteration residues to round-off errors and achieves the full convergence. Furthermore, the proposed fast sweeping method is compared with both the fully convergent fast sweeping method using the multi-resolution WENO local solver and the popular third-order total variation diminishing Rung-Kutta (TVD-RK3) time-marching method for steady-state computations. Numerical results show that this new method is more efficient than them to fully converge to steady-state solutions.    

\end{abstract}

\normalsize

\section{Introduction}

Hyperbolic conservation laws play a fundamental role in modeling various steady-state phenomena across applications such as fluid dynamics, wave problems, optimal control, image analysis, etc. For spatial discretization of hyperbolic conservation laws, the finite difference weighted essentially non-oscillatory (WENO) schemes represent a class of widely adopted and computationally efficient approaches, which can achieve uniform high-order accuracy in smooth regions of the solution while maintaining sharp and essentially non-oscillatory (ENO) \cite{HartenOsher, RK2} transitions across discontinuities. The development of high-order WENO schemes was initiated in \cite{LiuOsherChen} for a third-order accurate finite volume method. A general framework for high-order finite difference WENO schemes was established by Jiang and Shu in \cite{JS}, which are often referred to as the ``WENO-JS'' schemes in the literature. Since then, high-order WENO schemes were further developed and improved extensively in many aspects. For example, to handle complex domain geometries, WENO schemes on unstructured meshes were designed in e.g. \cite{HS, DK2, ZS, Y.Liu, ZQ3, Bals1, TsDum}; Efforts have been made to improve the robustness, the accuracy or efficiency of high order WENO schemes, with strategies such as modifying the linear or nonlinear weights, modifying the smoothness indicators, using different approximation functions, etc. (see e.g. \cite{HAPo, BorgesCarmona, JunZhu2, BaZor, JYoon, CSHuang}).


High-order nonlinear spatial discretizations such as WENO schemes 
for steady-state problems of hyperbolic conservation laws yield large nonlinear systems of equations. A key factor which determines computational efficiency is to design fast iterative schemes for solving these highly coupled nonlinear systems. A class of efficient iterative schemes for solving hyperbolic steady-state partial differential equations (PDEs) are the fast sweeping methods. They use alternating sweeping strategy to cover a family of characteristics of hyperbolic PDEs in a certain direction simultaneously in each sweeping order. Coupled with the Gauss-Seidel iterations, these methods can achieve a fast convergence speed \cite{ZE}. High-order accuracy fast sweeping methods were originally developed to solve static Hamilton-Jacobi equations (e.g. \cite{ZZQ,XZZS,ZCLS,WZ}). For complex problems such as the nonlinear systems of hyperbolic conservation laws, high-order fixed-point fast sweeping WENO methods were developed \cite{WuLiang}, extending the fixed-point fast sweeping WENO methods for static Hamilton-Jacobi equations in \cite{ZZC}. Different from some other fast sweeping methods for hyperbolic conservation laws (e.g. \cite{Engqui, Lozano}), fixed-point
fast sweeping methods have explicit forms that lead to simple algorithms for complex nonlinear hyperbolic systems. They are free of inverse operation of nonlinear local systems, and can be applied in solving general nonlinear hyperbolic PDEs with any monotone numerical fluxes and high-order approximations easily. 

A typical issue when the high-order WENO-JS schemes are applied for solving steady-state problems of hyperbolic conservation laws is that the iteration residues do not fully converge to values of round-off error level for difficult steady-state problems \cite{SCW, SSCW}. Investigations in \cite{SCW} show that slight post-shock oscillations which exist in the high-order WENO-JS schemes and propagate from the regions near the shocks downstream to the smooth regions cause the issue. This kind of convergence-difficulties also appear in the high order fixed-point fast sweeping WENO methods \cite{WuLiang}, which are based on the WENO-JS local solver. The issue makes it difficult to determine the convergence criterion for the fast sweeping methods when different problems are solved, and challenging to apply the methods to complex problems in real applications. Furthermore, residue errors may dominate the simulation errors and affect the convergence of the numerical solution to the desired solution of the PDEs. On the other hand, non-traditional WENO local solvers based on unequal-sized substencils, for example the multi-resolution WENO (MRWENO) schemes \cite{JunZhu2}, are found to have much better convergence properties than the WENO-JS schemes \cite{taijie2,JUNZ3}. In \cite{LZZ}, with the adoption of MRWENO scheme as the local solver, fully convergent fixed-point fast sweeping methods (called ``absolutely convergent'' fixed-point fast sweeping WENO methods) were developed.  The new absolutely convergent fixed-point fast sweeping WENO methods improve the convergence property of the original fixed-point fast sweeping methods in \cite{WuLiang}, and their iteration residues fully converge to values of machine zero / round-off error level for difficult steady-state problems. 
However, these non-traditional WENO schemes based on unequal-sized substencils such as the MRWENO schemes are more complex and in general more expensive in computational costs than the classical WENO-JS schemes, as shown in \cite{TsybulZhang}.
How to preserve the inherent advantage of simple construction and computational efficiency offered by the WENO-JS local solver while achieving full convergence in the fixed-point fast sweeping methods for hyperbolic conservation laws is still a challenging and interesting topic. 

In this paper, we return to the classical WENO-JS local solver and develop a new fully convergent fifth-order fixed-point fast sweeping method for solving steady-state problems of hyperbolic conservation laws. Based on recent studies on the nonlinear weighting process around discontinuities of solution \cite{liangfu,FW}, we apply the technique of frozen weights for avoiding unnecessary adjustment of nonlinear weights in the WENO-JS local solver, which freezes the nonlinear weights once the residual sequence in the fast sweeping iterations has stabilized. Different from the existing work on frozen weights in \cite{FW}, we design a simple and robust approach to judge stabilization of iteration residues and determine the step when the nonlinear weights are frozen in iterations of the fast sweeping method. Careful studies on parameter-sensitivity are carried out to justify the chosen parameter in the proposed procedure.   
Extensive numerical experiments on various challenging two-dimensional steady-state problems verify that, unlike the original fixed-point fast sweeping method with the fifth-order WENO-JS local solver in \cite{WuLiang}, the new 
WENO-JS fast sweeping scheme consistently drives the iteration residues to round-off errors and achieves the full convergence. To demonstrate the improved performance on computational costs, the proposed WENO-JS fast sweeping method is compared with both the fully convergent fast sweeping method using the multi-resolution WENO local solver \cite{LZZ} and the popular third-order total variation diminishing Rung-Kutta (TVD-RK3) time-marching method \cite{RK1,RK2} for steady-state computations. Numerical results show that the new fast sweeping method is more efficient than them to fully converge to steady-state solutions. The rest of the paper is organized as follows. We describe the detailed numerical algorithm in Section 2. In Section 3, numerical experiments are performed to test and evaluate the proposed method, and to compare it with the other schemes. Concluding remarks and discussions are provided in Section 4.   


\section {Description of the numerical method}
\label{secfs}
\setcounter{equation}{0}
\setcounter{figure}{0}
\setcounter{table}{0}

We consider the steady state problems of two dimensional (2D) hyperbolic conservation laws with the form
\begin{equation}
F(U)_x+G(U)_y=R(x,y,U),  \qquad   (x,y)\in\Omega,
\label{(e2.1)}
\end{equation}
on a bounded domain $\Omega$ with suitable boundary conditions on $\partial\Omega$. Here $U$ denotes the vector of the conservative variables, $F(U)$ and $G(U)$ represent the flux vectors, and $R$ is the source term. For example, the steady Euler system of equations for compressible flows has
$U = (\rho, \rho u, \rho v, E)^T$, $F(U) = (\rho u, \rho u^2 + p, \rho u v, u (E + p))^T$, and $G(U) = (\rho v, \rho u v, \rho v^2 +p, v (E + p))^T$.
Here $\rho$ is the fluid density, $(u,v)^T $ is the velocity vector, and $p$ is the pressure. The total energy $E=\frac{p}{\gamma'-1}+\frac{1}{2}\rho(u^2+v^2)$ with the constant $\gamma'=1.4$ for the case of air. In the following, we first briefly review the fifth-order WENO-JS finite difference scheme for discretizing the equation (\ref{(e2.1)}) which serves as the local solver, and the fixed-point fast sweeping method for solving the resulted nonlinear algebraic system. Then the new fully convergent WENO-JS fast sweeping method is described in details.  

A uniform grid $\{(x_i, y_j)\}$ is used to partition the computational domain with the grid sizes $\Delta{x}=x_i-x_{i-1}$ and $\Delta{y}=y_j-y_{j-1}$ in the $x$-direction and the $y$-direction respectively. Denote $x_{i+1/2}=\frac{1}{2}(x_{i}+x_{i+1})$, and $y_{j+1/2}=\frac{1}{2}(y_{j}+y_{j+1})$; a computational cell $I_{i}=[x_{i-1/2},x_{i+1/2}]$ of the $x$-direction, and $J_{j}=[y_{j-1/2},y_{j+1/2}]$ of the $y$-direction. The flux terms in (\ref{(e2.1)}) are approximated using the conservative fifth-order WENO-JS finite difference scheme \cite{JS}  
\begin{equation}
(F(U)_{x}+G(U)_{y})|_{x=x_{i},y=y_{j}}\approx\frac{1}{\Delta{x}}(\hat{F}_{i+\frac{1}{2},j}-\hat{F}_{i-\frac{1}{2},j})+\frac{1}{\Delta{y}}(\hat{G}_{i,j+\frac{1}{2}}-\hat{G}_{i,j-\frac{1}{2}}),
\end{equation}
where $\hat{F}_{i+\frac{1}{2},j}$ and $\hat{G}_{i,j+\frac{1}{2}}$ are the numerical fluxes. For the computation of these numerical flux terms in a system of equations using high-order WENO schemes (e.g., the fifth-order WENO-JS), a robust approach is to first use the technique of characteristic decomposition, and then apply WENO approximations to the resulted scalar components individually. To save space, we omit the review of details on the technique of characteristic decomposition, and refer it to \cite{JS}. In the following, we focus on the scalar case of (\ref{(e2.1)}): $f(u)_x + g(u)_y = r(x, y, u)$. Due to the simple and efficient dimension-by-dimension framework of finite difference schemes, the discretizations of $f(u)_x$ and $g(u)_y$ follow a similar procedure, along the grid lines in the $x-$ and $y-$ directions respectively. Hence here we only describe the procedure for $f(u)_x$. To enforce upwinding or linear stability, the flux $f(u)$ is split to the positive and the negative parts: $f(u)=f^{+}(u)+f^{-}(u)$ with $\frac{df^{+}(u)}{du}\geq0$ and $\frac{df^{-}(u)}{du}\leq0$. The popular Lax-Friedrichs splitting is used here: $f^{\pm}(u)=\frac{1}{2}(f(u)\pm\alpha{u})$ and $\alpha=\max_{u}|f'(u)|$ over the relevant range of $u$. Along a grid line, e.g. \( y = y_j \), the fifth-order WENO-JS reconstruction algorithm is applied to compute the positive numerical flux $\hat{f}^{+}_{i+1/2,j}$ and the negative numerical flux $\hat{f}^{-}_{i+1/2,j}$ respectively at computational cells' boundaries, and the final numerical flux $\hat{f}_{i+1/2,j}=\hat{f}^{+}_{i+1/2,j}+\hat{f}^{-}_{i+1/2,j}$. The reconstruction algorithm, originated from the finite volume method, is applied to the finite difference approximations \cite{RK2,LiuOsherChen,JS}. For reconstructing 
$\{\hat{f}^{+}_{i+1/2,j}\}_{i=\ldots,1,2,\ldots}$, we first identify the numerical values $\{f^{+}(u_{i,j})\}_{i=\ldots,1,2,\ldots}$ to be cell averages of a function \( W (x) \) on the one dimensional (1D) computational cells \( I_i = [x_{i-1/2}, x_{i+1/2}], i=\ldots,1,2,\ldots, \) in the \( x \)-direction grid line \( y=y_j \). Here \( u_{i,j} \) denotes the numerical solution of \( u \) at the grid point \( (x_i, y_j) \). Then, using these information of cell averages $\overline{W}_{i}=\frac{1}{\Delta{x}}\int_{x_{i-1/2}}^{x_{i+1/2}}W(x)dx=f^{+}(u_{i,j}), i=\ldots,1,2,\ldots$, the numerical-flux values $\hat{f}^{+}_{i+1/2,j}$ are computed by the following WENO-JS reconstruction algorithm. Note that since the formulas of reconstruction procedure for $\hat{f}^{-}_{i+1/2,j}$ are mirror-symmetric with respect to the point $x_{i+1/2}$ of those for $\hat{f}^{+}_{i+1/2,j}$ in the following sub-section, they are omitted to be presented here to save space. 

\subsection{The fifth-order WENO-JS reconstruction for $\hat{f}^{+}_{i+1/2,j}$}

The fifth-order WENO-JS reconstruction for $\hat{f}^{+}_{i+1/2,j}$ is based on three substencils: 
$T_{1}=\{I_{i-2},I_{i-1},I_{i}\}$, $T_{2}=\{I_{i-1},I_{i},I_{i+1}\}$ and $T_{3}=\{I_{i},I_{i+1},I_{i+2}\}$. Three quadratic polynomials $p_{k}(x), k=1,2,3$ are  reconstructed to satisfy
\begin{equation}
\frac{1}{\Delta{x}}\int_{I_{m}}p_{k}(x)dx=\overline{W}_{m}, \quad\forall I_{m}\in T_{k},\quad k=1,2,3.
\end{equation}
Note that these $\overline{W}_{m}$'s are just the available values of $f^{+}(u_{m,j})$'s. The reconstructed polynomials provide three third-order reconstructed values at $x_{i+1/2}$:
\begin{equation}
\begin{cases}
p_1(x_{i+1/2})=\frac{1}{3}\overline{W}_{i-2}-\frac{7}{6}\overline{W}_{i-1}+\frac{11}{6}\overline{W}_{i},\\
p_2(x_{i+1/2})=-\frac{1}{6}\overline{W}_{i-1}+\frac{5}{6}\overline{W}_{i}+\frac{1}{3}\overline{W}_{i+1},\\
p_3(x_{i+1/2})=\frac{1}{3}\overline{W}_{i}+\frac{5}{6}\overline{W}_{i+1}-\frac{1}{6}\overline{W}_{i+2}.\\
\end{cases}
\label{3rdorder}
\end{equation}
To compute the nonlinear weights in the WENO-JS scheme, the smoothness indicators $\beta_{k}$ measuring the smoothness of the reconstructed polynomials $p_{k}(x), k=1, 2, 3$ on the cell $I_i$ are defined as \cite{JS}:
\begin{equation}
\beta_{k}=\sum^{2}_{\alpha=1}\int_{x_{i-1/2}}^{x_{i+1/2}}{\Delta{x}}^{2\alpha-1}(\frac{d^{\alpha}p_{k}(x)}{dx^{\alpha}})^{2}dx, \quad k=1,2,3.
\label{GHYZ}
\end{equation}
The implementation expressions of these $\beta_{k}$'s using the available cell averages of $W(x)$ are 
\begin{equation}
\begin{cases}
\beta_1=\frac{13}{12}(\overline{W}_{i-2}-2\overline{W}_{i-1}+\overline{W}_{i})^2+\frac{1}{4}(\overline{W}_{i-2}-4\overline{W}_{i-1}+3\overline{W}_{i})^2,\\
\beta_2=\frac{13}{12}(\overline{W}_{i-1}-2\overline{W}_{i}+\overline{W}_{i+1})^2+\frac{1}{4}(\overline{W}_{i-1}-\overline{W}_{i+1})^2,\\
\beta_3=\frac{13}{12}(\overline{W}_{i}-2\overline{W}_{i+1}+\overline{W}_{i+2})^2+\frac{1}{4}(3\overline{W}_{i}-4\overline{W}_{i+1}+\overline{W}_{i+2})^2.\\
\end{cases}
\end{equation}
The nonlinear weights are calculated based on the smoothness indicators and the linear weights, as the following:
\begin{equation}
\omega_{k}=\frac{\overline{\omega}_{k}}{\sum_{s=1}^{3}\overline{\omega}_{s}},\quad \overline{\omega}_{k}=\frac{d_{k}}{(\varepsilon+\beta_{k})^{2}},\quad k=1,2,3.
\end{equation}
The linear weights are $d_{1}=\frac{1}{10},d_{2}=\frac{3}{5},d_{3}=\frac{3}{10}$, which provide the coefficients to combine these three third-order reconstructions (\ref{3rdorder}) for a fifth-order reconstruction.
$\varepsilon$ is a small positive number taken to avoid the denominator to become $0$. 
The final fifth-order WENO-JS reconstruction of the numerical flux $\hat{f}^{+}_{i+1/2,j}$ is 
\begin{equation}
\hat{f}^{+}_{i+1/2,j}=\omega_{1}p_{1}(x_{i+1/2})+\omega_{2}p_{2}(x_{i+1/2})+\omega_{3}p_{3}(x_{i+1/2}).
\end{equation}

\subsection{The fixed-point fast sweeping method}
Discretization of the PDE by the fifth-order WENO-JS finite difference scheme leads to a nonlinear algebraic system
$$0=-(\hat{f}_{i+1/2,j}-\hat{f}_{i-1/2,j})/  \Delta x-(\hat{g}_{i,j+1/2}-\hat{g}_{i,j-1/2})/\Delta y+r(x_i,y_j,u_{ij}),$$
\begin{equation}
\hspace{2in}i=1,\cdots,N; \quad j=1,\cdots, M,
\label{eq3}
\end{equation}
where $N$ and $M$ are the number of grid points in the $x-$ and $y-$ directions respectively.
A popular approach to solve the system (\ref{eq3}) for steady state solution of hyperbolic conservation laws is to add a pseudo-time-derivative term and perform time-marching iterations using the third-order TVD Runge-Kutta method \cite{RK2}, which preserves both the linear and the nonlinear stability of high-order WENO spatial discretizations. Another class of efficient iterative methods for steady state hyperbolic conservation laws is the 
fixed-point fast sweeping scheme \cite{WuLiang,LZZ}, which has improved efficiency over the TVD Runge-Kutta method. Here we apply the fixed-point fast sweeping scheme to the nonlinear system (\ref{eq3}) obtained by the fifth-order WENO-JS scheme, which has the following form:
$$u_{i,j}^{n+1}=u_{i,j}^{n}+\frac{\gamma}{\alpha_{x}/\Delta x +\alpha_{y}/\Delta  y}L(u_{i-3,j}^{\ast},\cdots,u_{i+3,j}^{\ast};u_{i,j}^{n};u_{i,j-3}^{\ast},\cdots, u_{i,j+3}^{\ast}),$$
\begin{equation}
\hspace{2in} i=i_{1},\cdots,i_{N}; \quad j=j_{1},\cdots,j_{M}.
\label{eq11}
\end{equation}
$u_{i,j}^{n}$ is the $n$th step iteration value at the grid point $(x_i, y_j)$.  $\gamma$ is the Courant-Friedrichs-Lewy (CFL) number and $\frac{\gamma}{\alpha_{x}/\Delta{x}+\alpha_{y}/\Delta{y}}$ actually corresponds to the time-step size $\Delta{t_{n}}$ in time-marching iterations.
$\alpha_{x}=\max_{u}{|f'(u)|}$ and $\alpha_{y}=\max_{u}{|g'(u)|}$ represent the maximum characteristic speeds in each spatial direction. The nonlinear function $L$, which denotes the fifth-order WENO-JS discretization, depends on $13$ numerical values of the fifth-order WENO scheme's stencil. In the fast sweeping method, the Gauss-Seidel iterations and alternating direction sweepings are used.
The Gauss-Seidel philosophy requires that the newest
numerical values of $u$ are used in the fifth-order WENO-JS scheme's stencil whenever they are available.
The iterations do not just proceed in only one direction $i = 1:N, j = 1:M$ as in a usual time-marching iteration, but in
the following four alternating directions repeatedly:
$$\mbox{(1) } i=1:N, j=1:M;$$
$$\mbox{(2) } i=N:1, j=1:M;$$
$$\mbox{(3) } i=N:1, j=M:1;$$
$$\mbox{(4) } i=1:N, j=M:1.$$
This is expressed as the notation ``$i=i_{1},\cdots,i_{N};j=j_{1},\cdots,j_{M}$'' in the scheme (\ref{eq11})
to show the alternating sweeping directions in iterations.
Via using alternating direction sweepings, the characteristics property of hyperbolic PDEs is utilized, which is one key factor in fast sweeping methods \cite{ZE}. By combining it with the Gauss-Seidel philosophy, a fast convergence to steady state numerical solution of hyperbolic conservation laws is obtained. This will be verified for the proposed fully convergent WENO-JS fast sweeping method in the following numerical experiments. In the  Gauss-Seidel iterations, the newest numerical values on
the stencil of the WENO scheme are used if they are available. This important component of the method is demonstrated by
the notation $u^{\ast}$ in the scheme (\ref{eq11}) to represent the numerical values in the WENO stencil, with the understanding that $u_{k,l}^{\ast}$ could be $u_{k,l}^{n}$ or $u_{k,l}^{n+1}$, depending on the current sweeping direction.

\subsection{The fully convergent fifth-order fast sweeping WENO-JS scheme}
This scheme represents an improvement over the original fast sweeping WENO-JS method. When addressing steady-state problems, the conventional fast sweeping WENO-JS scheme encounters difficulties in reducing the residue to machine zero. Our investigation revealed that the iteration residual cannot drop below machine epsilon because the nonlinear weights fail to stabilize in shock regions. If the unnecessary adjustments of nonlinear weights can be avoided, the residue can continue to decrease. When the nonlinear weights near shock waves undergo meaningless adjustments, an effective approach is to freeze the nonlinear weights. 
Experimental results confirm that this approach is highly effective. It maintains the accuracy of the numerical solution while ensuring the residue converges to machine zero and preserving the essentially non-oscillatory property. The key to this method lies in determining the optimal moment to freeze the weights. Our strategy is to trigger this freeze once the residue sequence shows no further significant decrease. The following sections detail the two key components of this strategy: the criterion for assessing the residue sequence and the specific mechanism for freezing the weights.

\subsubsection{The criterion for assessing the residue sequence}
\begin{figure}
		\centering
\subfigure[Example2]
{\begin{minipage}[t]{0.3\linewidth}
\includegraphics[width=2.1in]{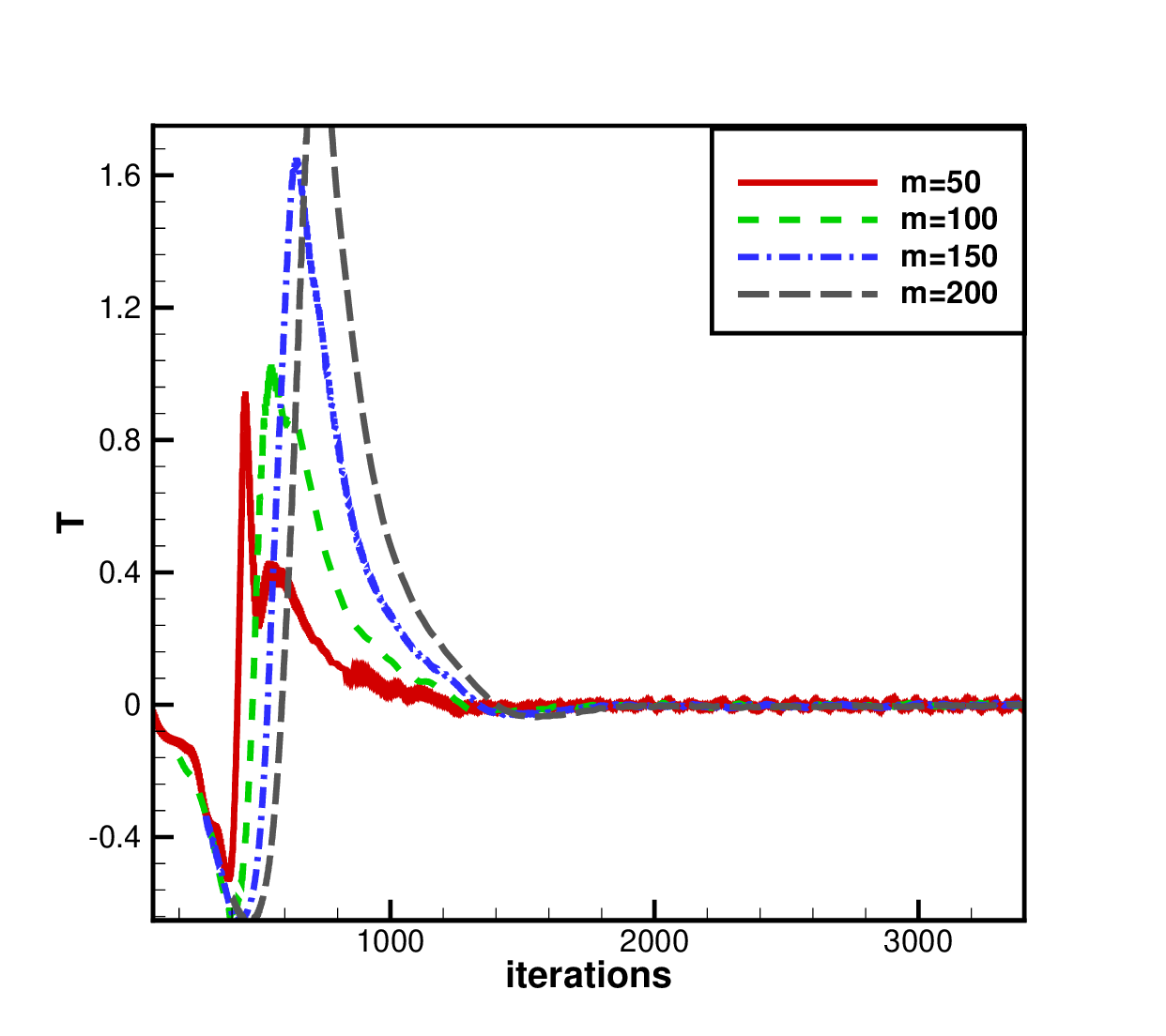}
\end{minipage}}
\subfigure[Example3]
{\begin{minipage}[t]{0.3\linewidth}
\includegraphics[width=2.1in]{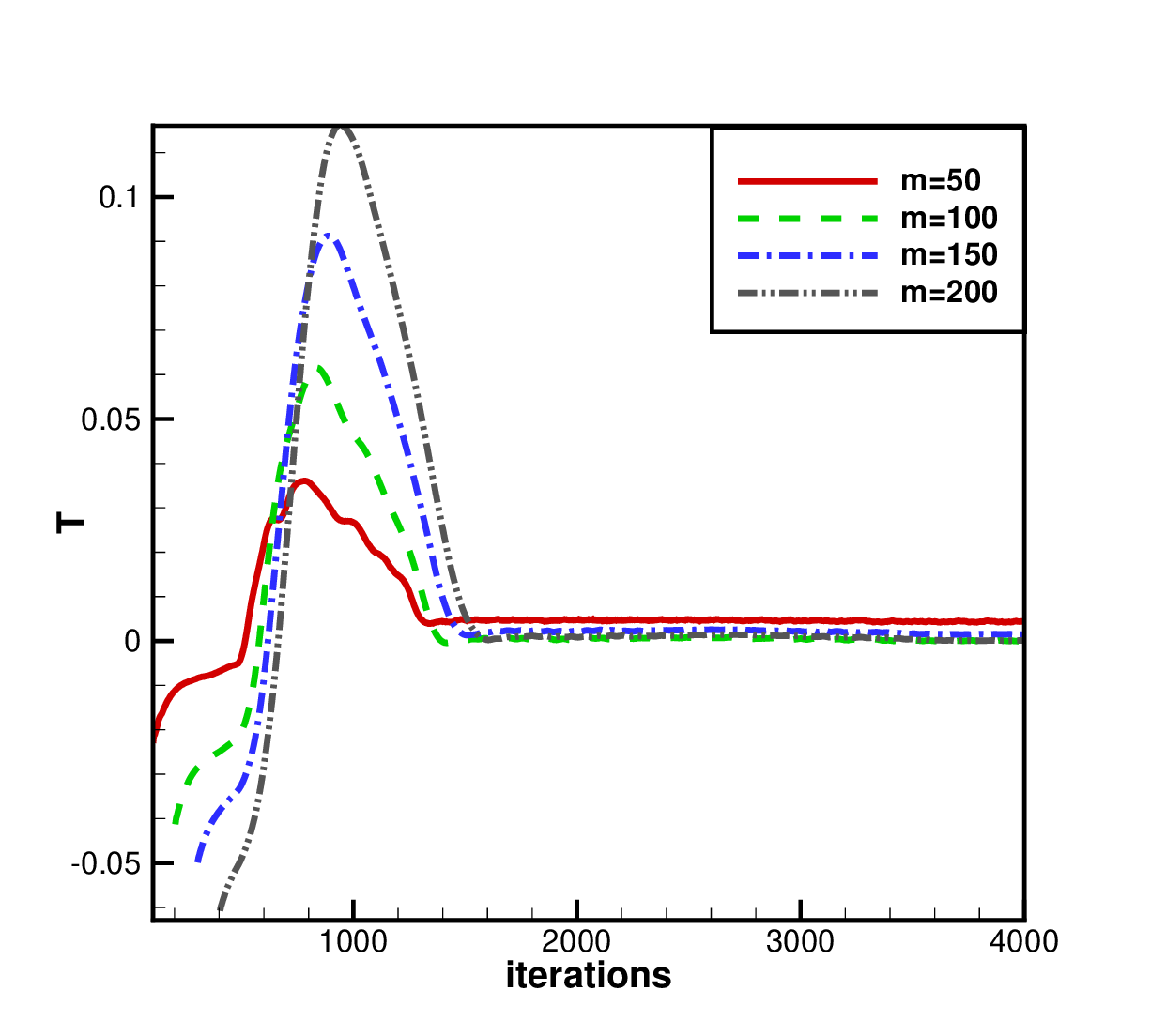}
\end{minipage}}
\subfigure[Example4]
{\begin{minipage}[t]{0.3\linewidth}
\includegraphics[width=2.1in]{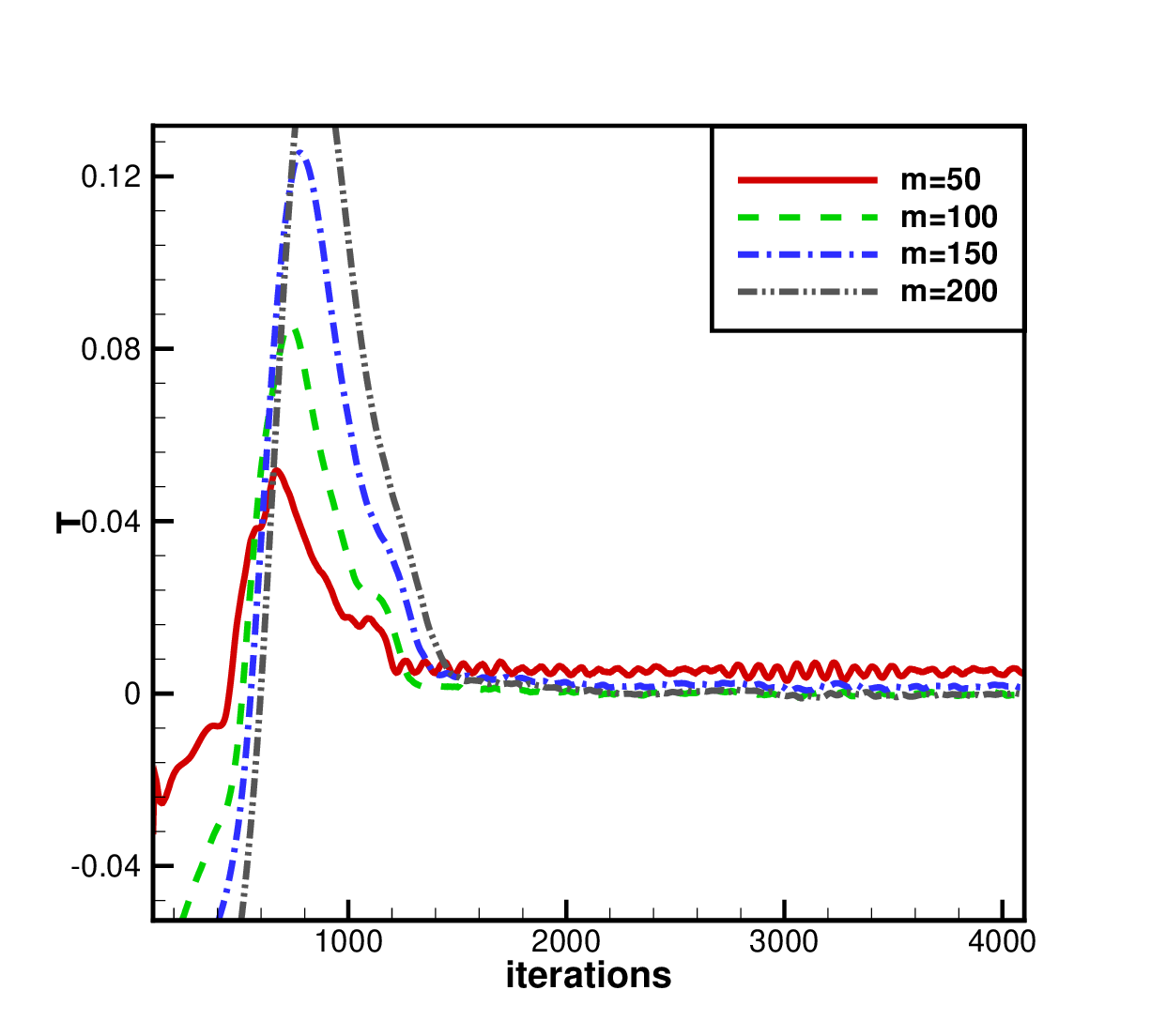}
\end{minipage}}

\subfigure[Example5]
{\begin{minipage}[t]{0.3\linewidth}
\includegraphics[width=2.1in]{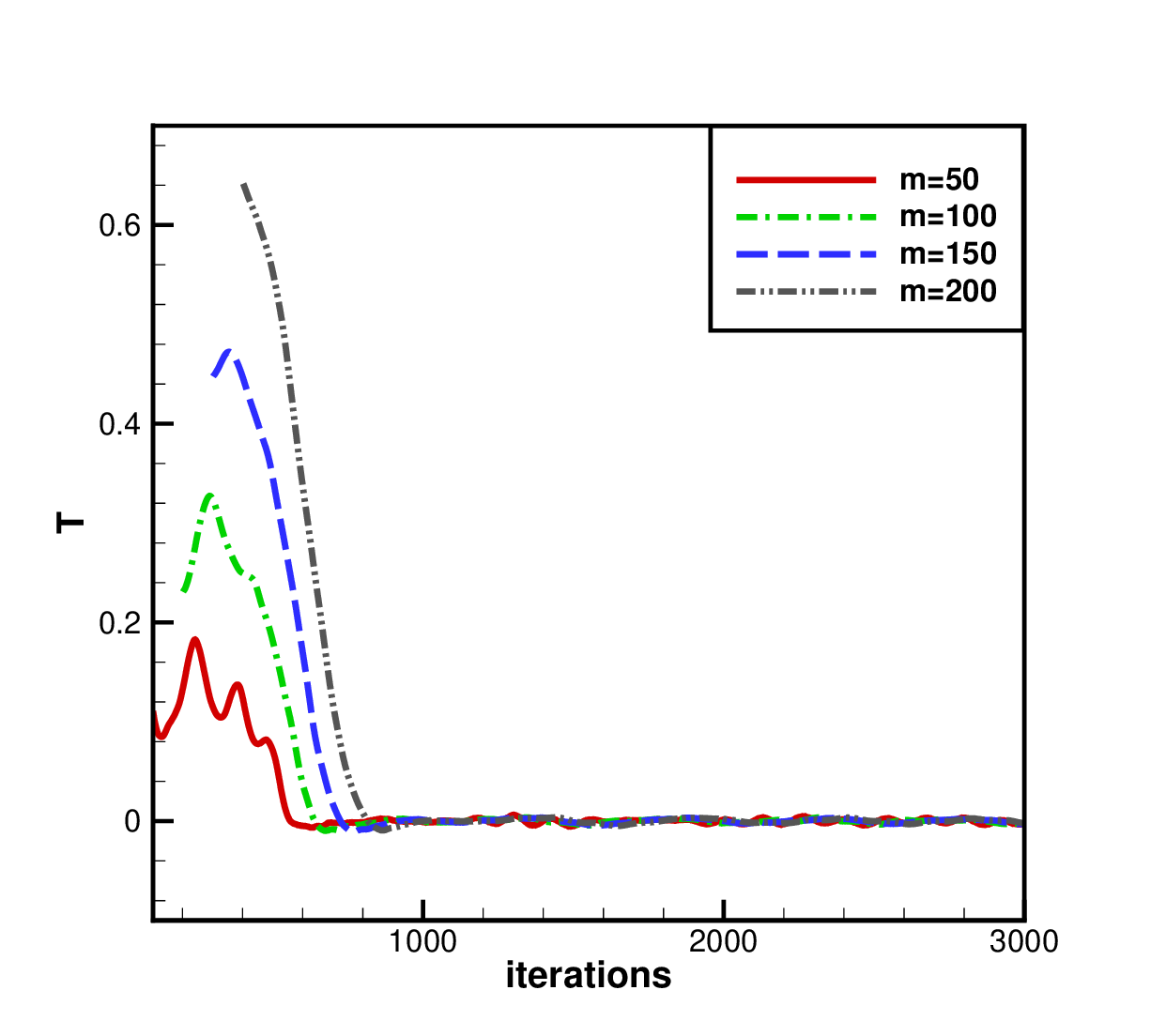}
\end{minipage}}
\subfigure[Example6]
{\begin{minipage}[t]{0.3\linewidth}
\includegraphics[width=2.1in]{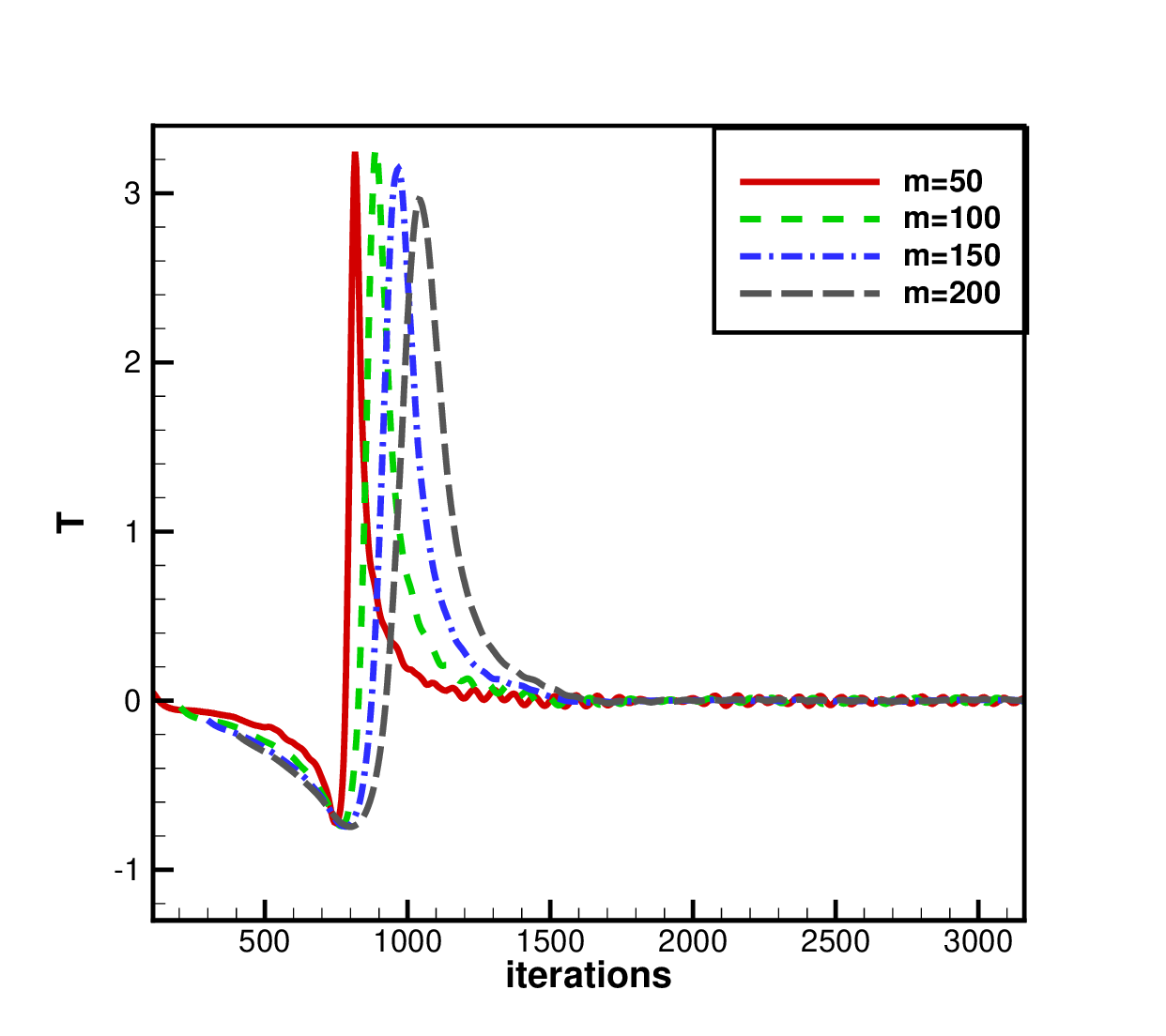}
\end{minipage}}
\subfigure[Example7]
{\begin{minipage}[t]{0.3\linewidth}
\includegraphics[width=2.1in]{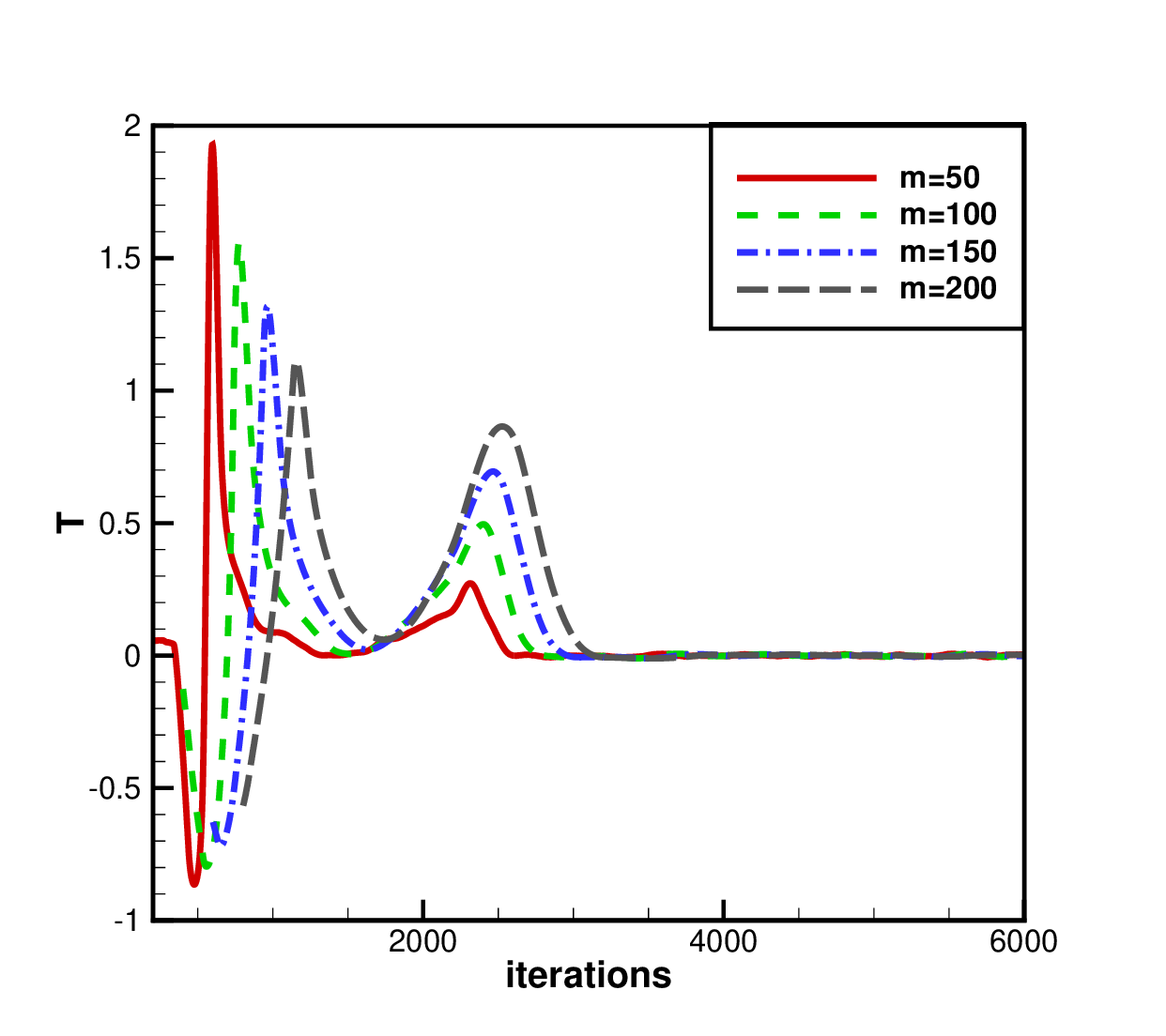}
\end{minipage}}
\caption{\label{aa}  T(n) as a function of iterations under different m in the FS-WENOJS scheme (Example 2-7, CFL=0.6).}
\end{figure}

Each iteration yields a average residue value, denoted as $ResA(n)$ for the n-th iteration. The value is defined as
\begin{equation}
ResA(n)=\sum_{i=1}^{N}\frac{|R1_{i}|+|R2_{i}|+|R3_{i}|+|R4_{i}|}{4\times N}
\end{equation}
where $R\ast_{i}$ are the local residues of different conservative variables, that is, $R1_{i}=\frac{\partial \rho}{\partial t}\mid_{i}=\frac{\rho_{i}^{n+1}-\rho_{i}^{n}}{\Delta t}, R2_{i}=\frac{\partial (\rho u)}{\partial t}\mid_{i}=\frac{(\rho u)_{i}^{n+1}-(\rho u)_{i}^{n}}{\Delta t}, R3_{i}=\frac{\partial (\rho v)}{\partial t}\mid_{i}=\frac{(\rho v)_{i}^{n+1}-(\rho v)_{i}^{n}}{\Delta t}, R4_{i}=\frac{\partial E}{\partial t}\mid_{i}=\frac{E_{i}^{n+1}-E_{i}^{n}}{\Delta t}$. Determining whether the iteration residue has stabilized essentially involves assessing when the sequence $ResA$ becomes stable. When computing each test case using the traditional WENOJS scheme, the iteration residue stabilize at different levels. For example, in Example 2 they stabilize at around $10^{-1}$, in Example 3 at around $10^{-3.3}$ and in Example 4 at around $10^{-2.8}$. Therefore, it is not straightforward to set a single iteration residue criterion to determine when the residue has stabilized. Therefore, we need to process the residues in order to more easily and accurately determine when the residue has stabilized. A straightforward and effective approach is to use a relative error criterion. We set
\begin{equation}
Q(n)=\sum_{k=1}^{m}\log(ResA(n+1-k)),P(n)=\sum_{k=1}^{m}\log(ResA(n+1-m-k))
\end{equation}
When the relative error satisfy $|T(n)|=|\frac{Q(n)-P(n)}{P(n)}|<\varepsilon$, the residue is considered stabilized. The parameter n in the formula represents the number of iterations, m is a positive integer, and $\varepsilon$ is a small positive value. Through our testing of different test cases, we have selected relatively appropriate parameters. In this paper, m is set at 100, and $\varepsilon$ is set to $10^{-3}$. The parameter selection accurately identifies when the residue has stabilized in all numerical experiments. Compared with the method in reference \cite{FW}, the technique proposed in this paper is much simpler. To illustrate the rationality of the parameter selection, we tested the evolution plot of $T(n)$ for the numerical example in figure \ref{aa} using different values of m. As can be seen from the figure \ref{aa}, the curves of $T(n)$ are similar for different values of m. When m is too large, the convergence becomes slower, while when m is too small, $T(n)$ becomes oscillatory. Therefore, we chose a moderate value of $m=100$. Furthermore, we present in the table \ref{cc} the number of iterations when the nonlinear weights are frozen, for different values of m. It can be observed that varying m has little effect on the iteration count. This again confirms that the choice of m is not overly sensitive for weight freezing.

\begin{figure}
		\centering
\subfigure[Example2]
{\begin{minipage}[t]{0.3\linewidth}
\includegraphics[width=2.1in]{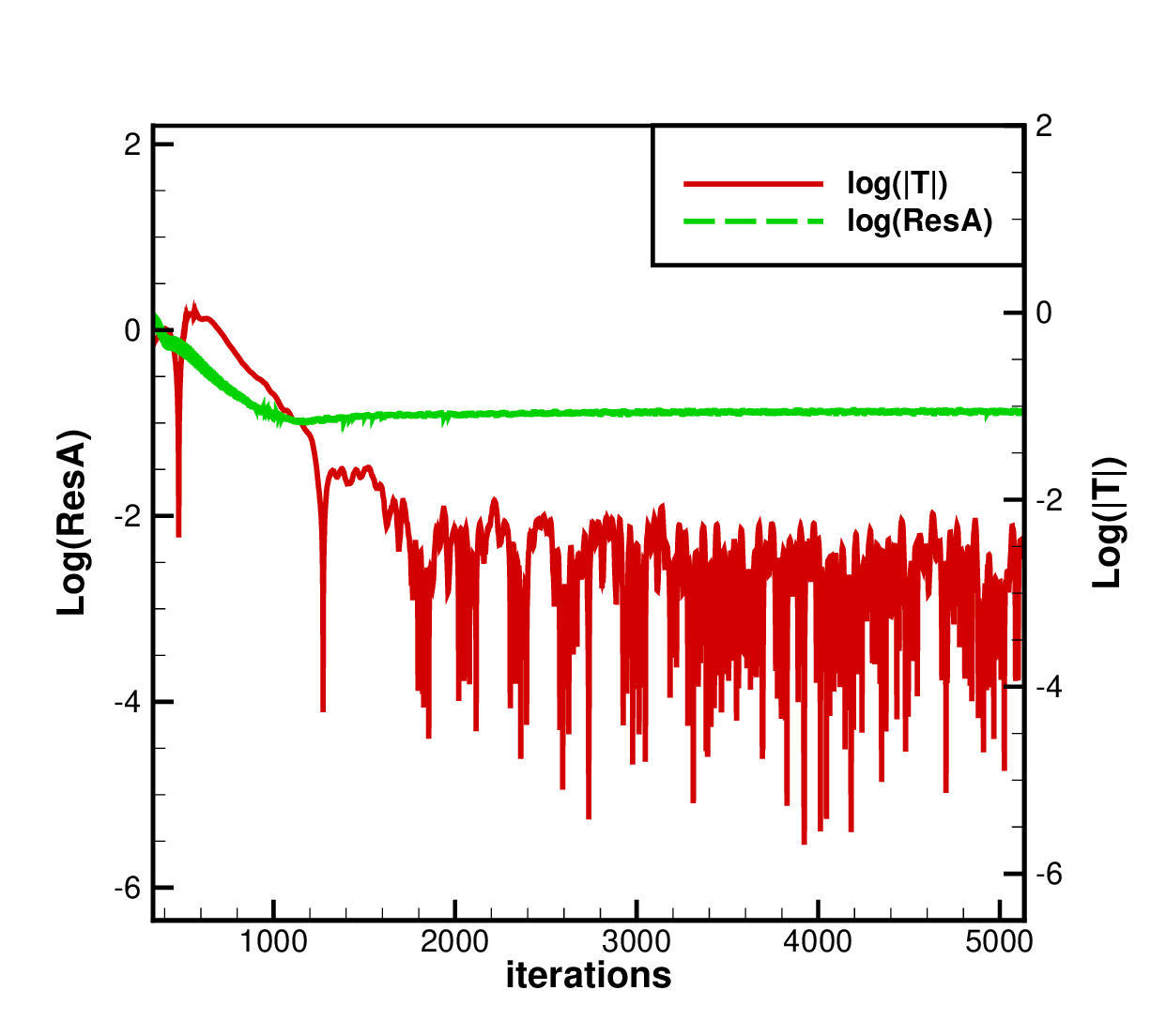}
\end{minipage}}
\subfigure[Example3]
{\begin{minipage}[t]{0.3\linewidth}
\includegraphics[width=2.1in]{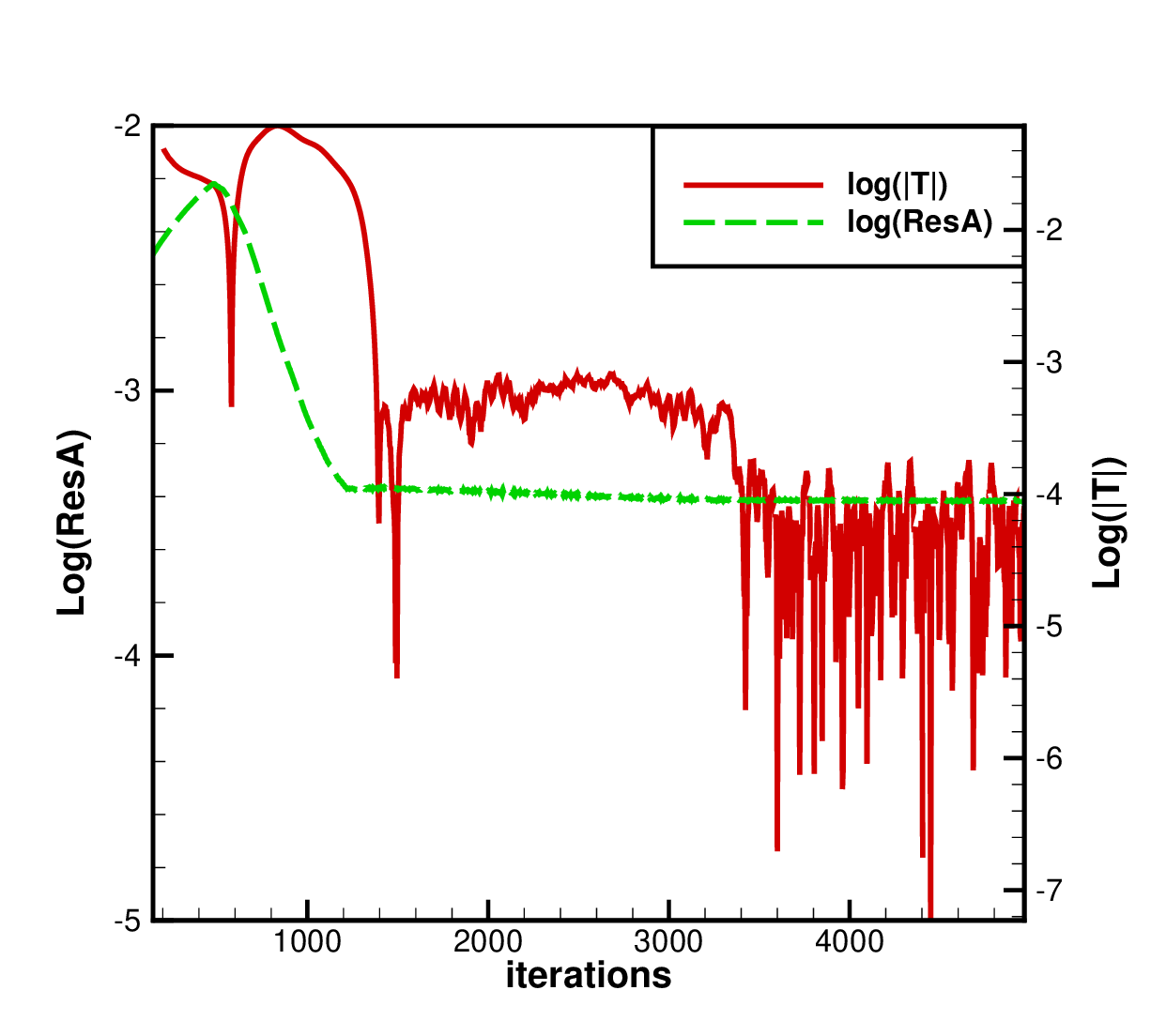}
\end{minipage}}
\subfigure[Example4]
{\begin{minipage}[t]{0.3\linewidth}
\includegraphics[width=2.1in]{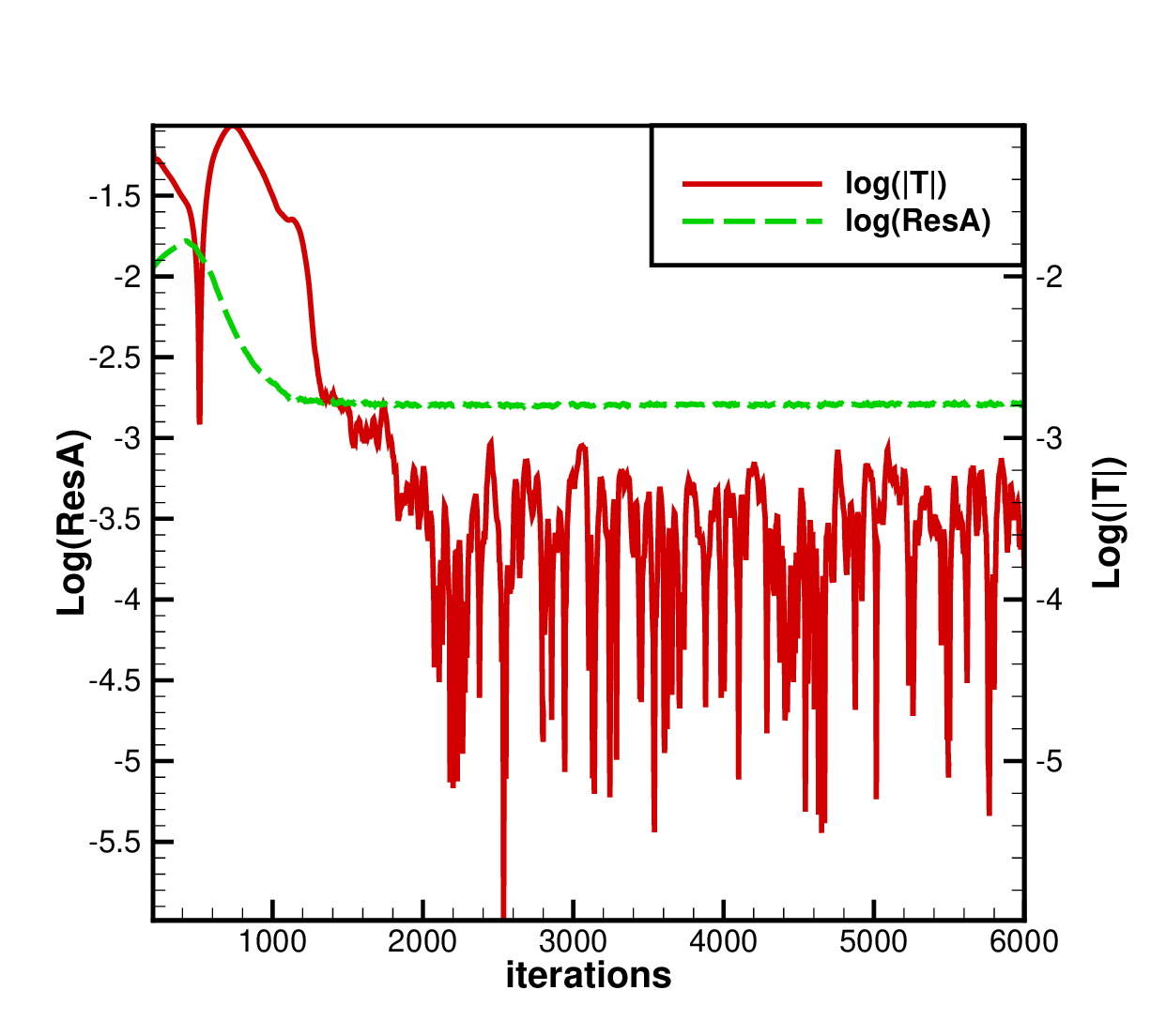}
\end{minipage}}

\subfigure[Example5]
{\begin{minipage}[t]{0.3\linewidth}
\includegraphics[width=2.1in]{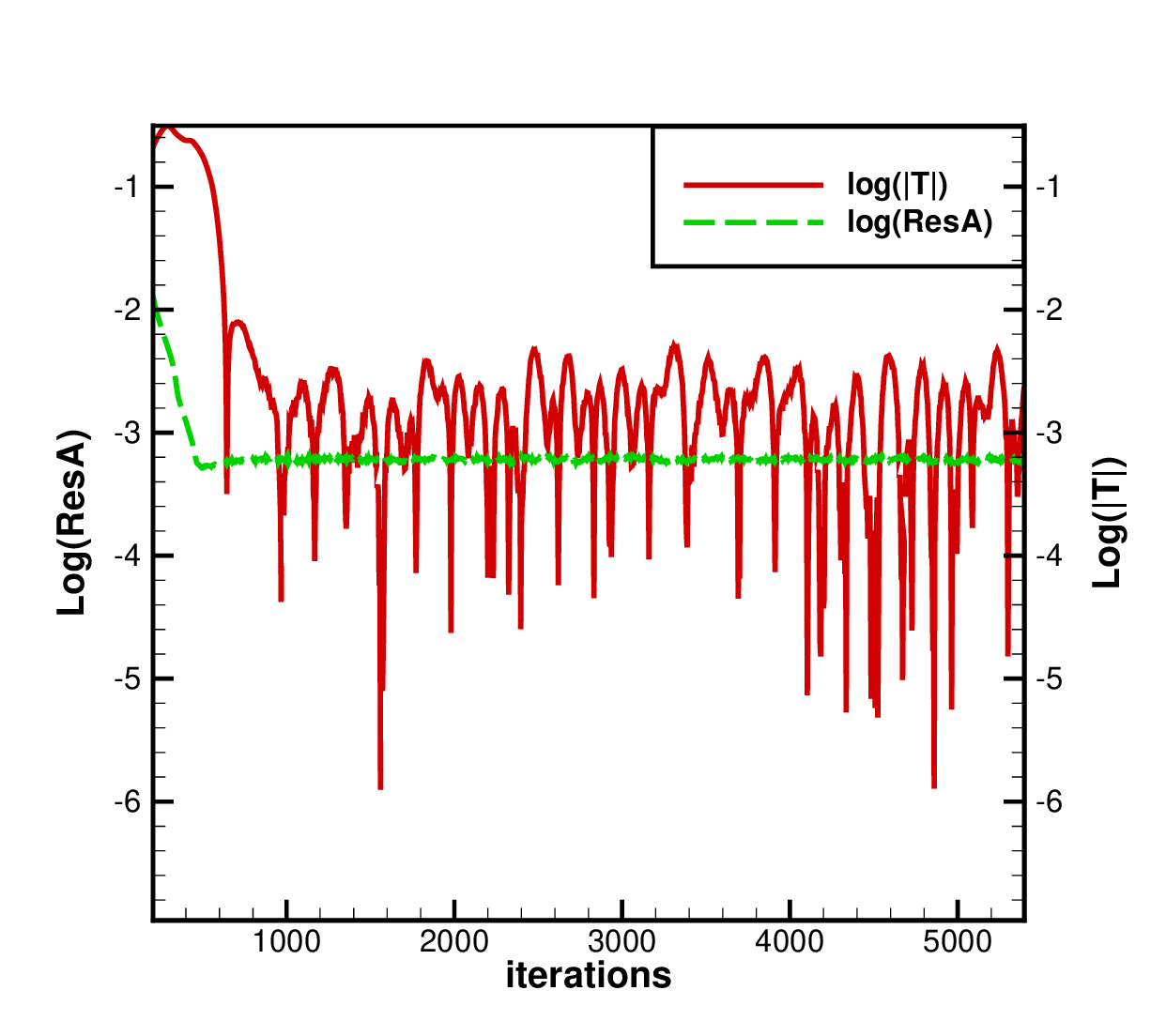}
\end{minipage}}
\subfigure[Example6]
{\begin{minipage}[t]{0.3\linewidth}
\includegraphics[width=2.1in]{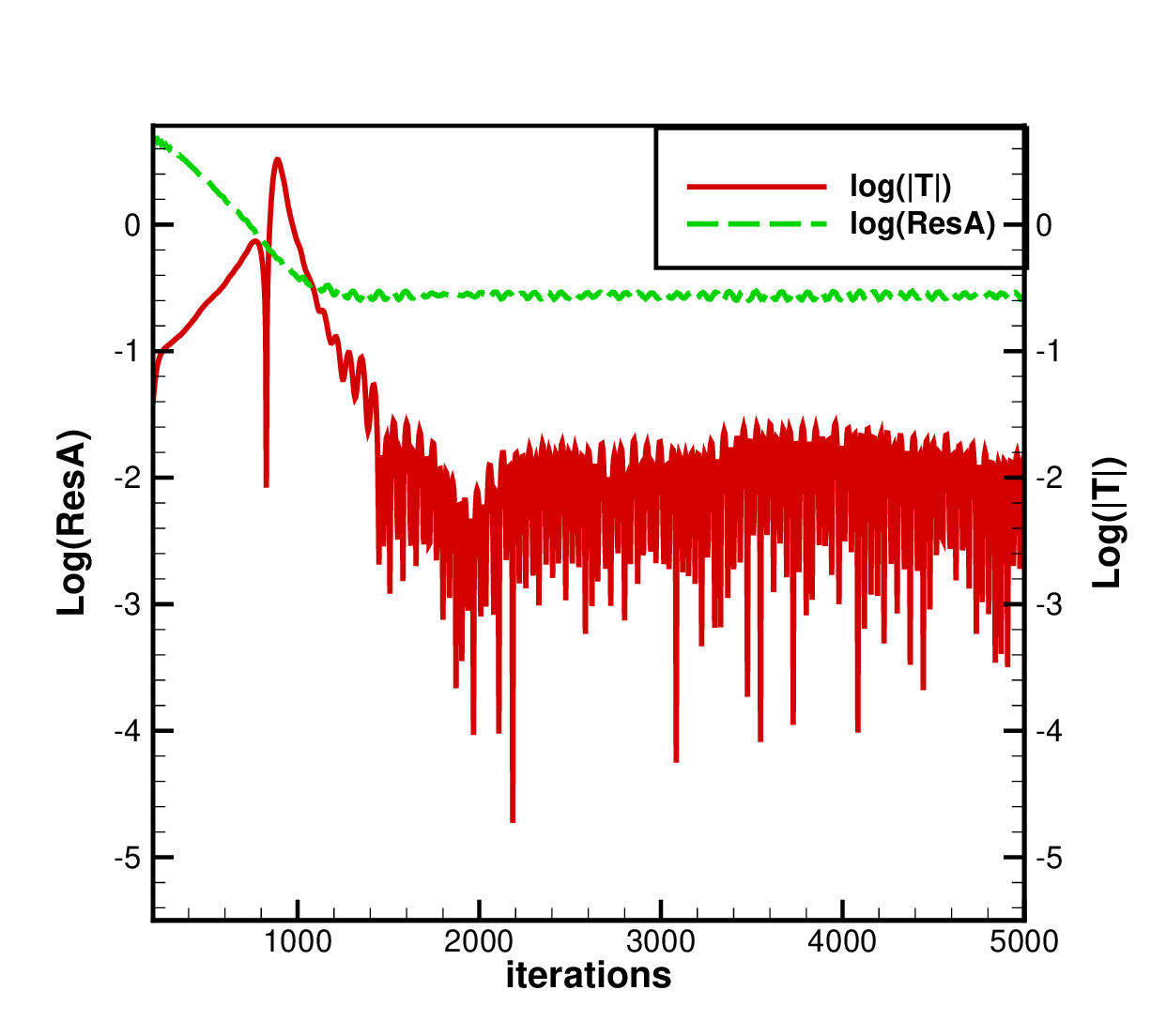}
\end{minipage}}
\subfigure[Example7]
{\begin{minipage}[t]{0.3\linewidth}
\includegraphics[width=2.1in]{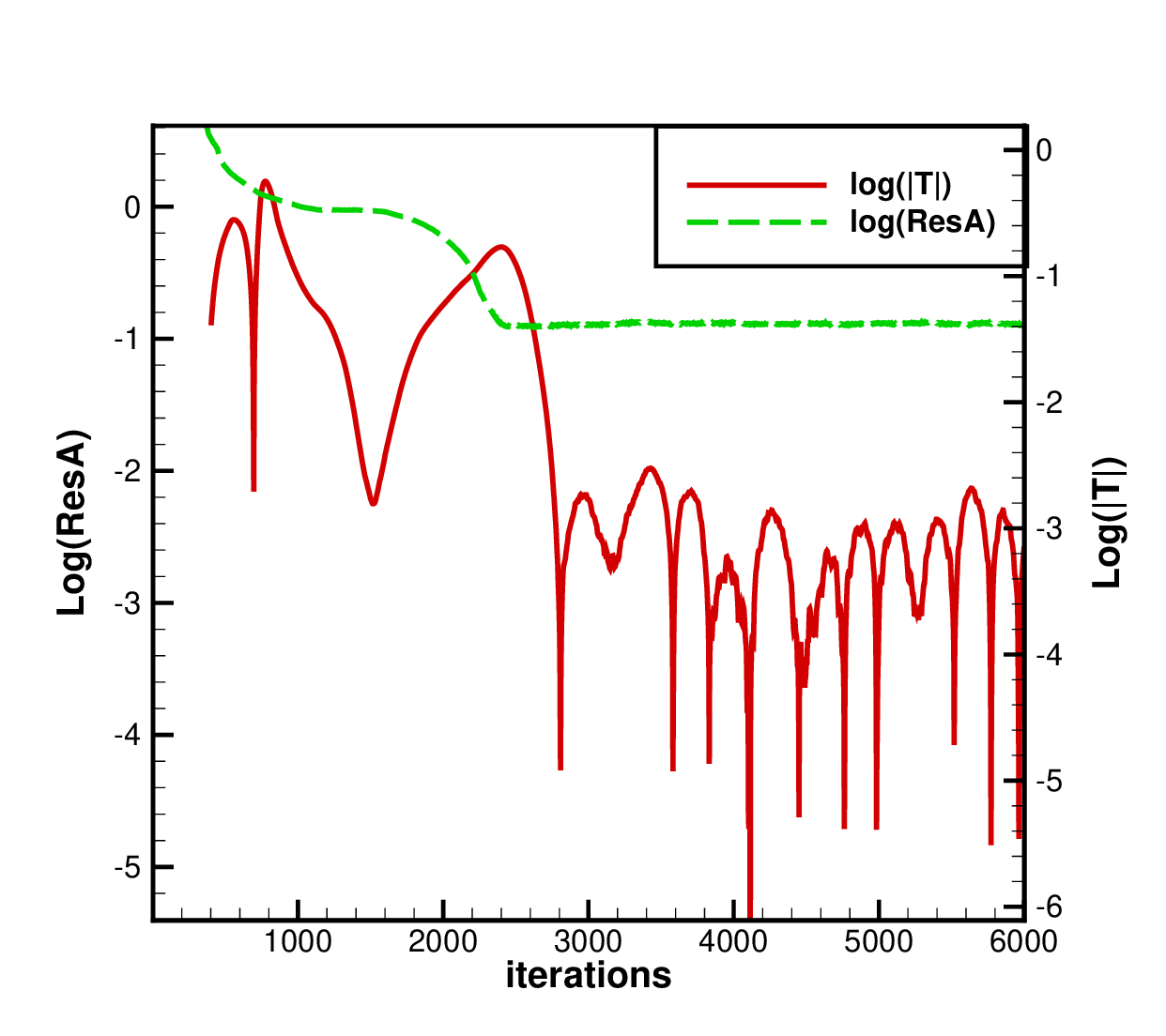}
\end{minipage}}
\caption{\label{bb} Comparison of the residues $Log(ResA)$ and relative error $T(n)$ for FS-WENOJS scheme (Example 2-7, CFL=0.6, m=100).}
\end{figure}

Figure \ref{bb} shows a comparison between the iteration residue $\log(ResA)$ and the relative error $\log{T(n)}$ when $m = 100$. We can see that although the residues for each test case stabilize at different levels, the relative errors all approach zero. Moreover, the residues stabilize as the relative errors vanish. Therefore, it is effective to judge the stabilization of the residues based on the relative error. It is also observed that once the residual stabilizes, the value of $\mid{T(n)}\mid$ remains around $10^{-3}$. Hence, we set the error tolerance to $10^{-3}$. Note that setting error tolerance to $10^{-4}$ is also acceptable, it only delays the freezing of the nonlinear weights.

\begin{table}
		\centering
\begin{tabular}{|c|c|c|c|c|c|c|}\hline
       m     &Example 2&Example 3&Example 4&Example 5&Example 6&Example 7\\\hline
50&1271 &1172	&1072	&644 &	 1840    &1492       \\                           \hline
100&1273 &1344	&1308	&952 &	 1920   & 2838       \\                           \hline
150&1333 &1490	&1448	&808 &	 1840    &  3358    \\ \hline
200&1399 &1616	&1576	&1060 &	 1816    &  3334    \\ \hline
		\end{tabular}
\caption{\label{cc}The number of iterations when freezing the weights with different values of m in Example 2-7 for the FS-WENOJS-AC scheme.}
\end{table}
\subsubsection{Freezing the Nonlinear Weights}

Once the iteration residue has stabilized, the corresponding nonlinear weights are stored. For all subsequent numerical flux computations, the method abandons the nonlinear WENO reconstruction. Instead, the stored nonlinear weights are employed directly as the linear weights. At this point, the scheme has essentially turned into a linear scheme. This approach avoids recomputing the nonlinear weights, thereby significantly saving computational time while allowing the residue to drop to near machine zero. It is important to note that, as can be seen from the figure, $T(n)$ often crosses the x-axis once before reaching a steady state. In such cases, although $|T(n)|<10^{-3}$, the residue has not yet stabilized. To avoid such misjudgments, we need to slightly modify the stopping criterion to require both $|T(n)|<10^{-3}$ and $|T(n-m)|<10^{-3}$. Specifically, once the iteration reaches step $n$, the nonlinear weights are stored and reused for all subsequent steps. Since the fast sweeping method employs four sweeping directions, the corresponding boundary conditions may vary with each sweep such as the boundary conditions on the short plate in Example 2. Therefore, when storing the nonlinear weights, they should be saved separately, with one set stored for each sweeping direction. This approach often leads to faster convergence.

{\bf Remark 1:}
Although the proposed frozen weights strategy is presented and demonstrated within the framework of the classical WENO-JS scheme, its applicability is by no means limited to this particular variant. This strategy can be directly extended to other WENO variants, such as WENO-Z, which improves accuracy near critical points yet still encounters the same stagnation issue in steady state computations. In fact, preliminary tests (not shown here for brevity) confirm that combining the present freezing technique with WENO-Z yields essentially identical convergence behavior and computational savings. The same holds for other variants that employ different smoothness indicators or nonlinear weight mappings, as long as the weights are computed locally from the solution and tend to oscillate after the residual levels off.

{\bf Remark 2:}
A key point to clarify is why freezing the nonlinear weights does not compromise the essentially nonoscillatory (ENO) property. The freezing procedure is activated only after the residual has stabilized, at which stage the numerical solution has already approached the steady state, and the locations of discontinuities such as shocks and contact discontinuities are well resolved and essentially fixed. By this time, the nonlinear weights have already adapted according to the local smoothness of the solution: weights associated with stencils that cross discontinuities are appropriately suppressed, while those lying in smooth regions receive larger values. Since the subsequent iterative updates induce only negligible changes in the solution, the frozen weight distribution remains a valid choice for preventing spurious oscillations. Consequently, fixing the weights does not introduce oscillatory behavior, because the linear combination of substencils at the frozen state still constitutes a convex combination that heavily favors smooth substencils near shocks exactly as in the original nonlinear reconstruction. In essence, the ENO property is already encoded in the weights at the moment of freezing, and the subsequent linear evolution remains stable due to the minimal variation in the solution. This explains why the frozen weight scheme retains sharp shock capturing without producing oscillatory artifacts.
\section{Numerical experiments}

In this section, we test several challenging cases for which numerous papers have shown that the conventional TVD-RK3 WENOJS scheme cannot converge the residues to machine zero. We refer to the two types of WENO schemes based on equal-sized spatial stencils constructed in this paper as the absolutely convergent fast sweeping WENO-JS scheme and the absolutely convergent Runge-Kutta WENO-JS scheme, denoted as FS-WENOJS-AC and RK-WENOJS-AC, respectively. These difficult cases will be used to evaluate the performance of the new scheme. To demonstrate the excellent convergence and high computational efficiency of this scheme, we compare it with the fast sweeping multi-resolution WENO scheme, which has been proven to have good convergence in $\cite{LZZ}$. In the following context, this scheme is referred to as the FS-MRWENO scheme. Unless otherwise stated, the convergence criterion for all numerical examples is $10^{-12}$, and the CPU time is measured when the residual falls below this tolerance. To more comprehensively evaluate the performance of the new scheme, we also tested a case in which the iteration residue of the FS-WENOJS scheme can also settle down to machine zero. The purpose of testing this case is to examine whether the new algorithm would prematurely freeze the weights when the residues are still decreasing, and to verify that the weight-freezing operation based on computing relative errors does not require excessive additional time.

\bigskip
\noindent{\bf Example 1. 2D accuracy test}

\noindent We first use the two-dimensional Euler equations to test the accuracy and computational efficiency of those numerical schemes.
\[
\frac{\partial}{\partial x}
\begin{pmatrix}
\rho u \\
\rho u^2 + p \\
\rho uv \\
u(E + p)
\end{pmatrix}
+
\frac{\partial}{\partial y}
\begin{pmatrix}
\rho v \\
\rho uv \\
\rho v^2 + p \\
v(E + p)
\end{pmatrix}
= \begin{pmatrix}
0.4\cos(x+y) \\
0.6\cos(x+y) \\
0.6\cos(x+y) \\
1.8\cos(x+y)
\end{pmatrix}.
\]
The exact steady-state solution \(\rho(x, y) = 1 + 0.2 \sin(x + y), u(x, y) = 1, v(x, y) = 1\), and \(p(x, y) = 1 + 0.2 \sin(x + y)\) on the domain \((x, y) \in [0,2\pi]^2 \) is used as the initial guess. With the exact solution prescribed as boundary conditions, the iterative schemes drive this non-satisfying solution to their numerical steady states. The convergence criterion is \(10^{-12}\). Table \ref{a1} reports the \(L_1\) and \(L_\infty\) numerical errors and accuracy orders, and CPU times from a mesh refinement study. From the table it can be seen that all schemes achieve the designed order of accuracy. For this problem, the residuals of all schemes can be reduced smoothly to $10^{-12}$, so neither the RK-WENOJS-AC nor the FS-WENOJS-AC scheme freezes the nonlinear weights. As shown in the table, the FS-WENOJS-AC scheme requires slightly more computational time than the FS-WENOJS scheme; however, the extra cost incurred by evaluating $T(n)$ is negligible compared with the overall cost of the whole problem. The table also reveals that the MRWENO scheme is approximately half as efficient as the WENOJS scheme, due to the simpler construction of the latter. Furthermore, the fast sweeping method saves about 70\% of the computational time relative to the Runge-Kutta method, which is consistent with the conclusions reported in many relevant studies. From this numerical example, it is evident that the FS-WENOJS-AC scheme is significantly more efficient than both the FS-MRWENO and RK-WENOJS-AC schemes, and this conclusion will be further confirmed in the subsequent examples.
\begin{table}
	\centering
	\caption{Example 1. Accuracy, $L^{1}$ and $ L^{\infty}$ error and order of four schemes. CFL=1.0}\label{a1}
\begin{tabular}{|c|c|c|c|c|c|c|}\hline
			&\multicolumn{6}{|c|}{FS-MRWENO }\\    \hline
			N$\times$N  &  L$_{1}$ error  &  L$_{1}$ order&  L$_{\infty}$ error  &  L$_{\infty} $order  & iteration number&  CPU time\\\hline
   $40\times40$ & 4.49E-07 &  & 1.13E-06 &  & 1010 & 8.08\\
\hline
$50\times50$ & 1.50E-07 & 4.92 & 3.74E-07 & 4.96 & 1213 & 17.25\\
\hline
$60\times60$ & 6.08E-08 & 4.94 & 1.51E-07 & 4.95 & 1421 & 27.69\\
\hline
$70\times70$ & 2.83E-08 & 4.95 & 7.04E-08 & 4.97 & 1622 & 44.03\\
\hline
$80\times80$ & 1.46E-08 & 4.96 & 3.62E-08 & 4.99 & 1814 & 60.97\\
\hline
   			&\multicolumn{6}{|c|}{RK-WENOJS-AC }\\    \hline
  $40\times40$ & 2.77E-06 &  & 7.54E-06 &  & 3513 & 22.94\\
\hline
$50\times50$ & 9.18E-07 & 4.95 & 2.43E-06 & 5.08 & 4071 & 31.97\\
\hline
$60\times60$ & 3.69E-07 & 5.01 & 9.63E-07 & 5.08 & 4737 & 49.63\\
\hline
$70\times70$ & 1.69E-07 & 5.04 & 4.40E-07 & 5.09 & 5391 & 73.50\\
\hline
$80\times80$ & 8.61E-08 & 5.06 & 2.22E-07 & 5.12 & 6045 & 104.83\\
\hline
      &\multicolumn{6}{|c|}{FS-WENOJS }\\    \hline
  $40\times40$ & 2.76E-06 &  & 7.54E-06 & & 1125 & 5.39\\
\hline
$50\times50$ & 9.14E-07 & 4.95 & 2.43E-06 & 5.08 & 1297 & 9.38\\
\hline
$60\times60$ & 3.67E-07 & 5.00 & 9.63E-07 & 5.08 & 1505 & 16.22\\
\hline
$70\times70$ & 1.69E-07 & 5.04 & 4.40E-07 & 5.09 & 1709 & 22.95\\
\hline
$80\times80$ & 8.59E-08 & 5.06 & 2.22E-07 & 5.12 & 1913 & 32.78\\
\hline
   &\multicolumn{6}{|c|}{FS-WENOJS-AC }\\    \hline
$40\times40$ & 2.76E-06 &  & 7.54E-06 &  & 1125 & 5.52\\
\hline
$50\times50$ & 9.14E-07 & 4.95 & 2.43E-06 & 5.08 & 1297 & 11.23\\
\hline
$60\times60$ & 3.67E-07 & 5.00 & 9.63E-07 & 5.08 & 1505 & 17.16\\
\hline
$70\times70$ & 1.69E-07 & 5.04 & 4.40E-07 & 5.09 & 1709 & 23.05\\
\hline
$80\times80$ & 8.59E-08 & 5.06 & 2.22E-07 & 5.12 & 1913 & 33.20\\
\hline
   	\end{tabular}
	\end{table}

\bigskip
\noindent{\bf Example 2. Regular shock reflection}

\begin{figure}
		\centering
\subfigure[FS-MRWENO]
{\begin{minipage}[t]{0.3\linewidth}
\includegraphics[width=2.0in]{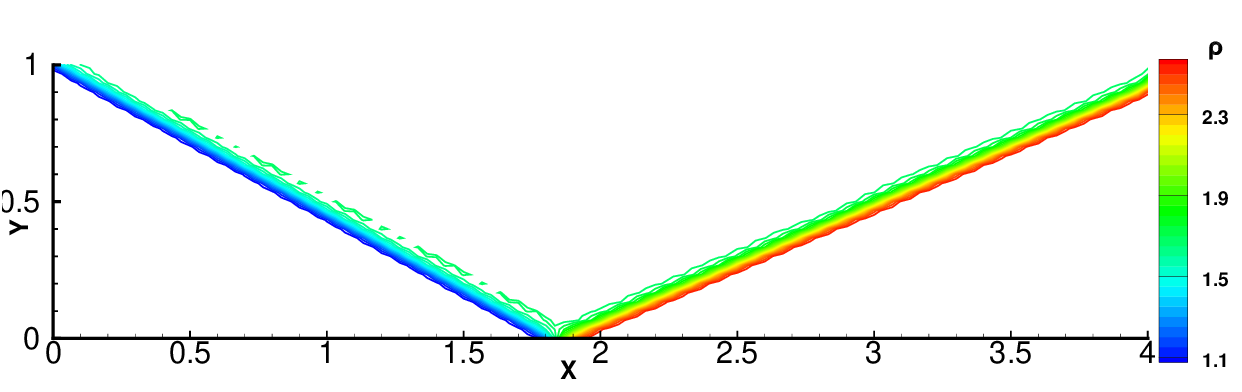}
\end{minipage}}
\subfigure[RK-WENOJS-AC]
{\begin{minipage}[t]{0.3\linewidth}
\includegraphics[width=2.0in]{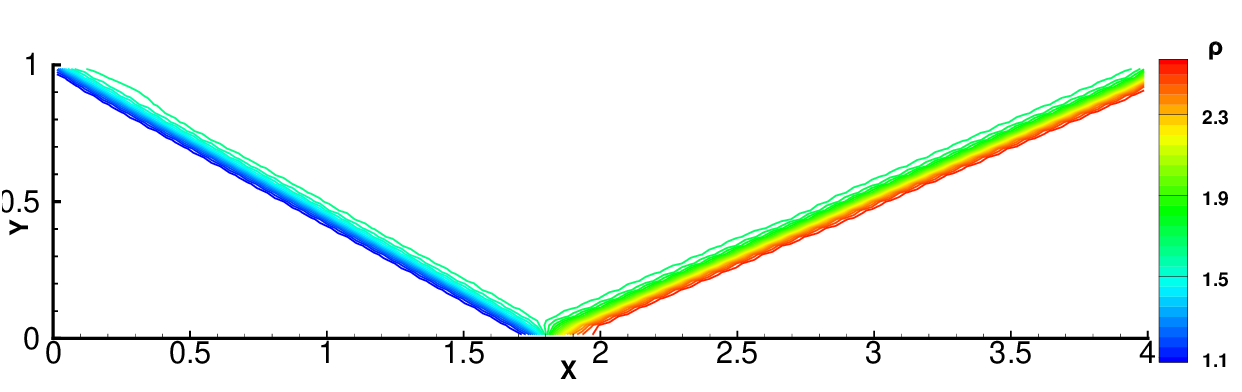}
\end{minipage}}
\subfigure[FS-WENOJS-AC]
{\begin{minipage}[t]{0.3\linewidth}
\includegraphics[width=2.0in]{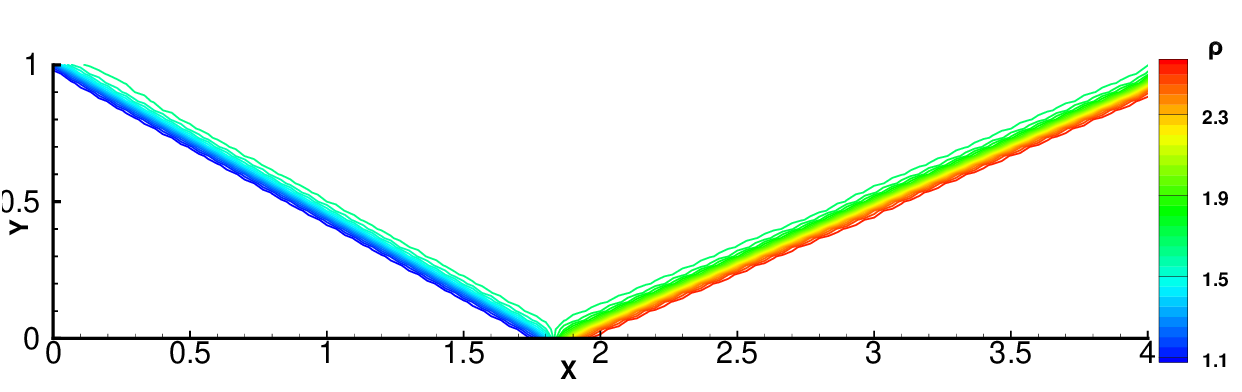}
\end{minipage}}
\caption{\label{1.2}Example 2: Thirty equally spaced density contours from 1.1 to 2.6 of the converged steady states of numerical solutions by three WENO schemes.}
\end{figure}

\begin{figure}
		\centering
\subfigure[FS-MRWENO]
{\begin{minipage}[t]{0.28\linewidth}
\includegraphics[width=2in]{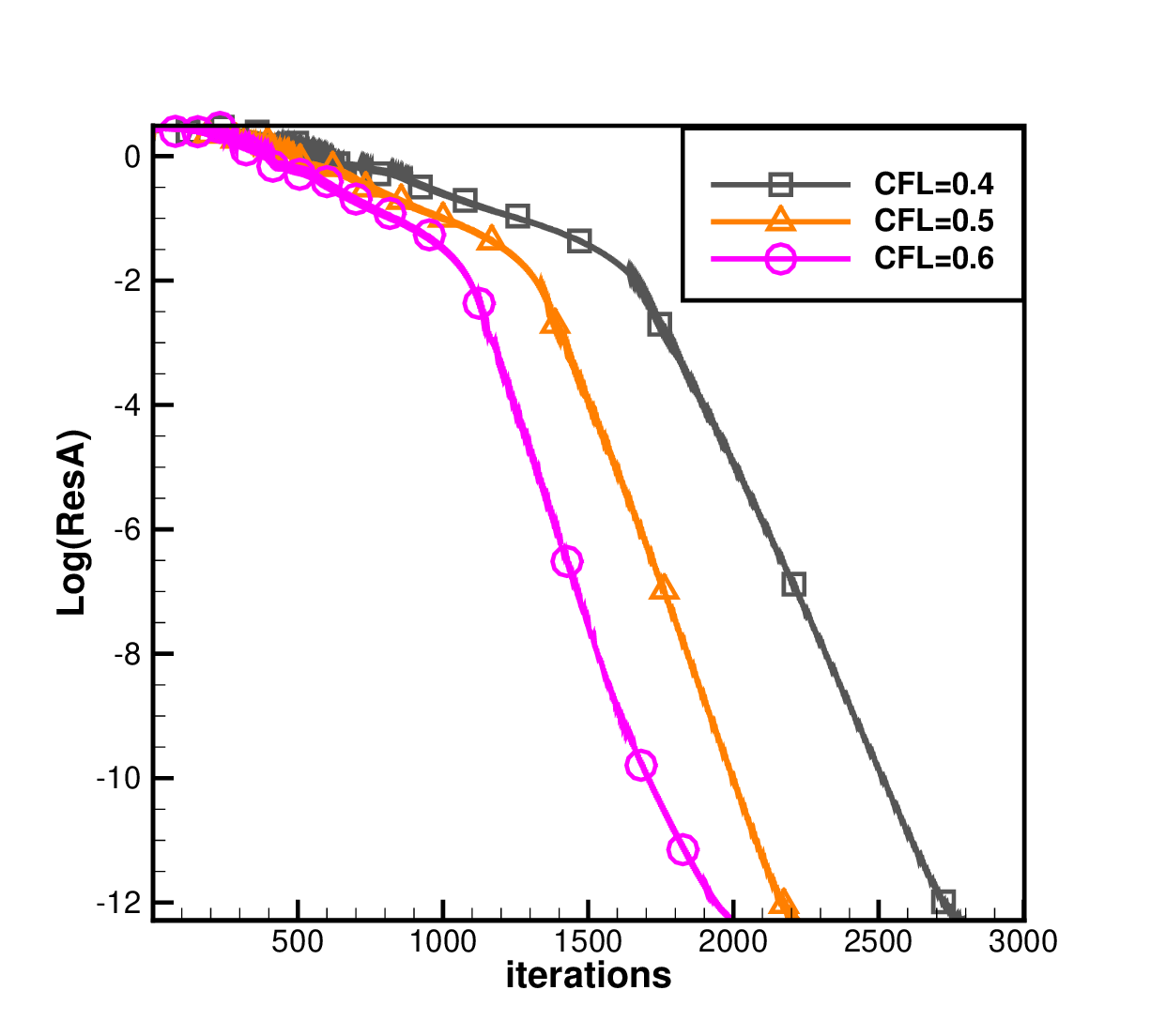}
\end{minipage}}
\subfigure[FS-WENOJS]
{\begin{minipage}[t]{0.28\linewidth}
\includegraphics[width=2in]{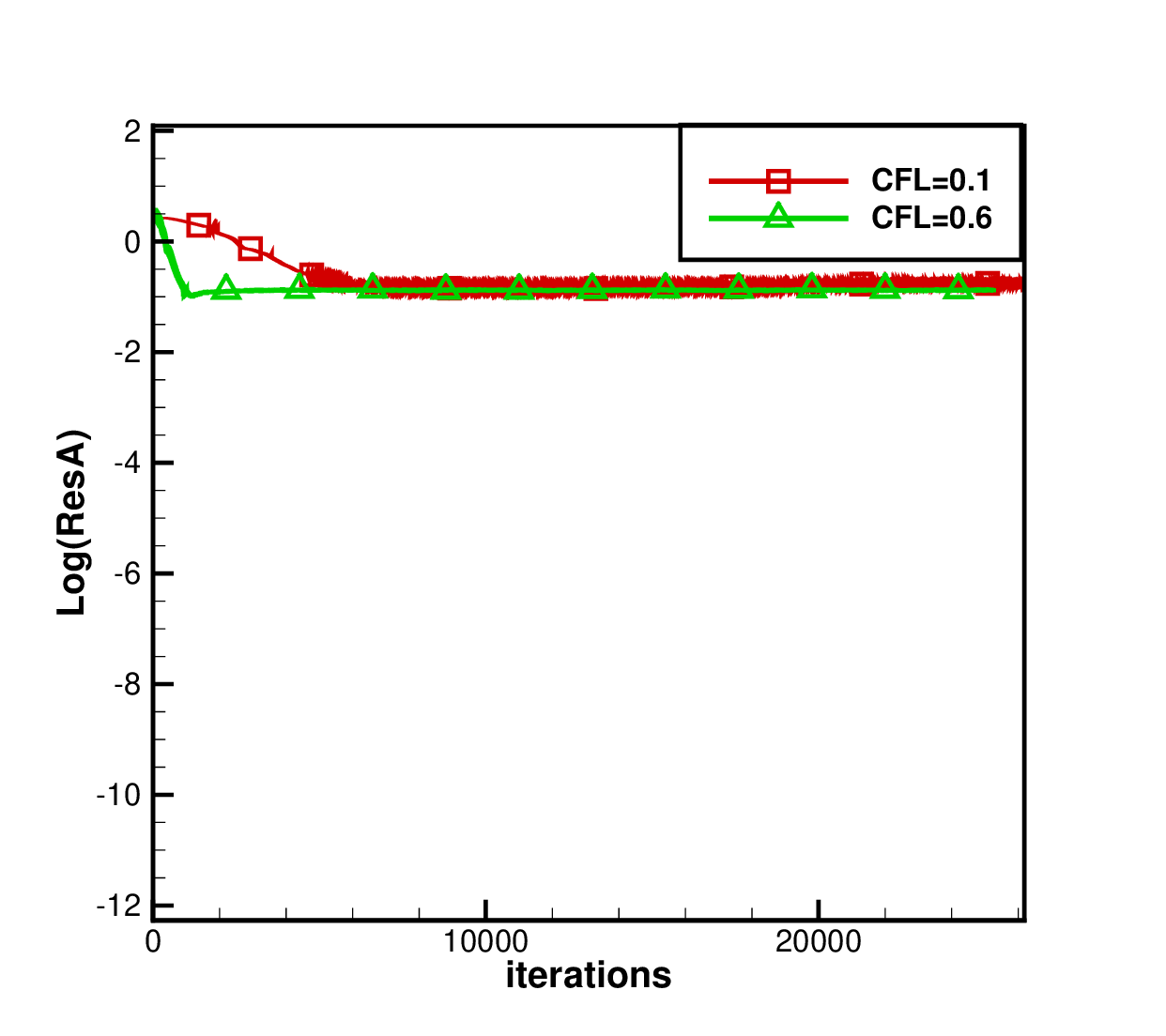}
\end{minipage}}
\subfigure[FS-WENOJS-AC]
{\begin{minipage}[t]{0.28\linewidth}
\includegraphics[width=2in]{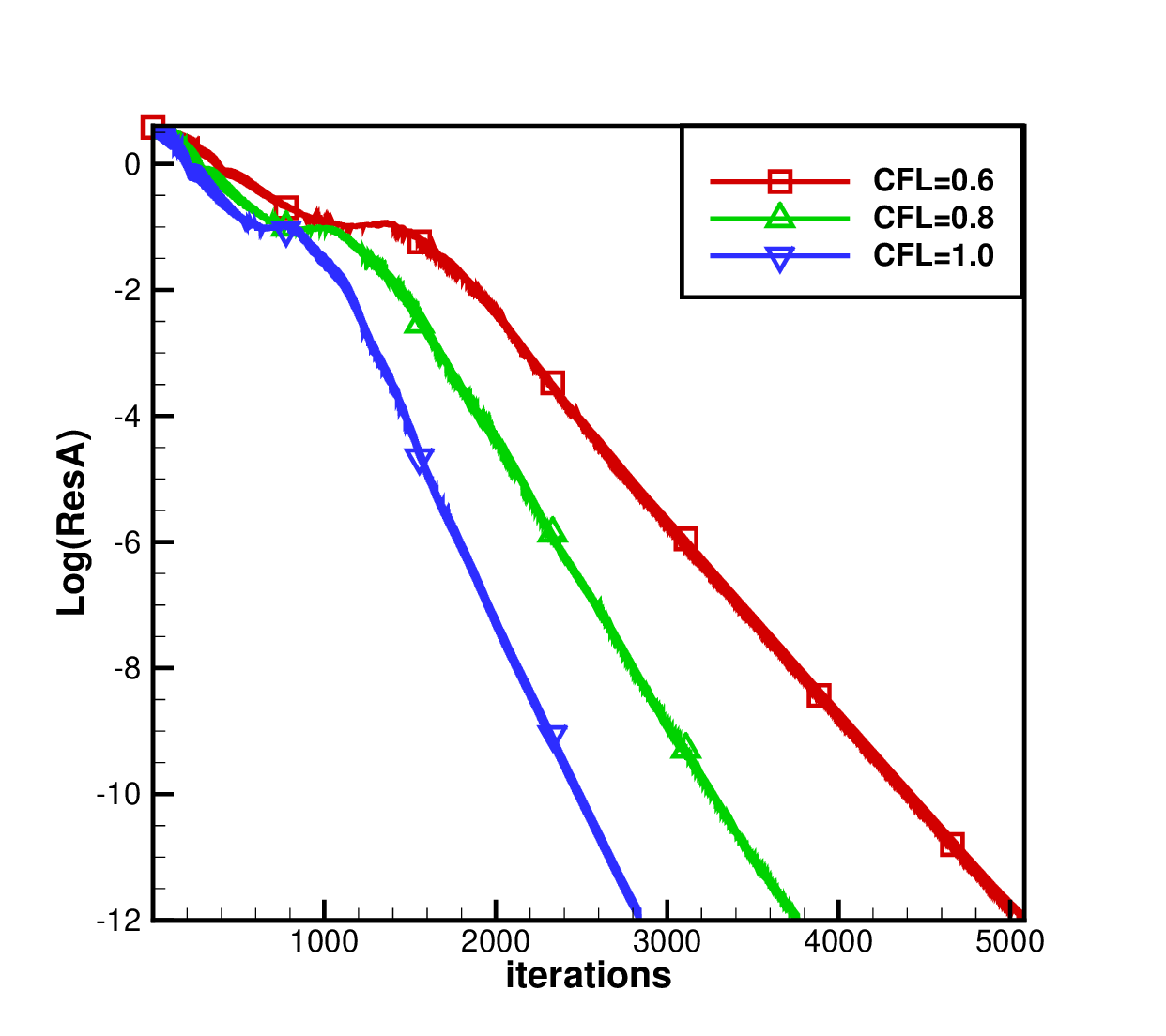}
\end{minipage}}

\subfigure[RK-WENOJS-AC]
{\begin{minipage}[t]{0.28\linewidth}
\includegraphics[width=2in]{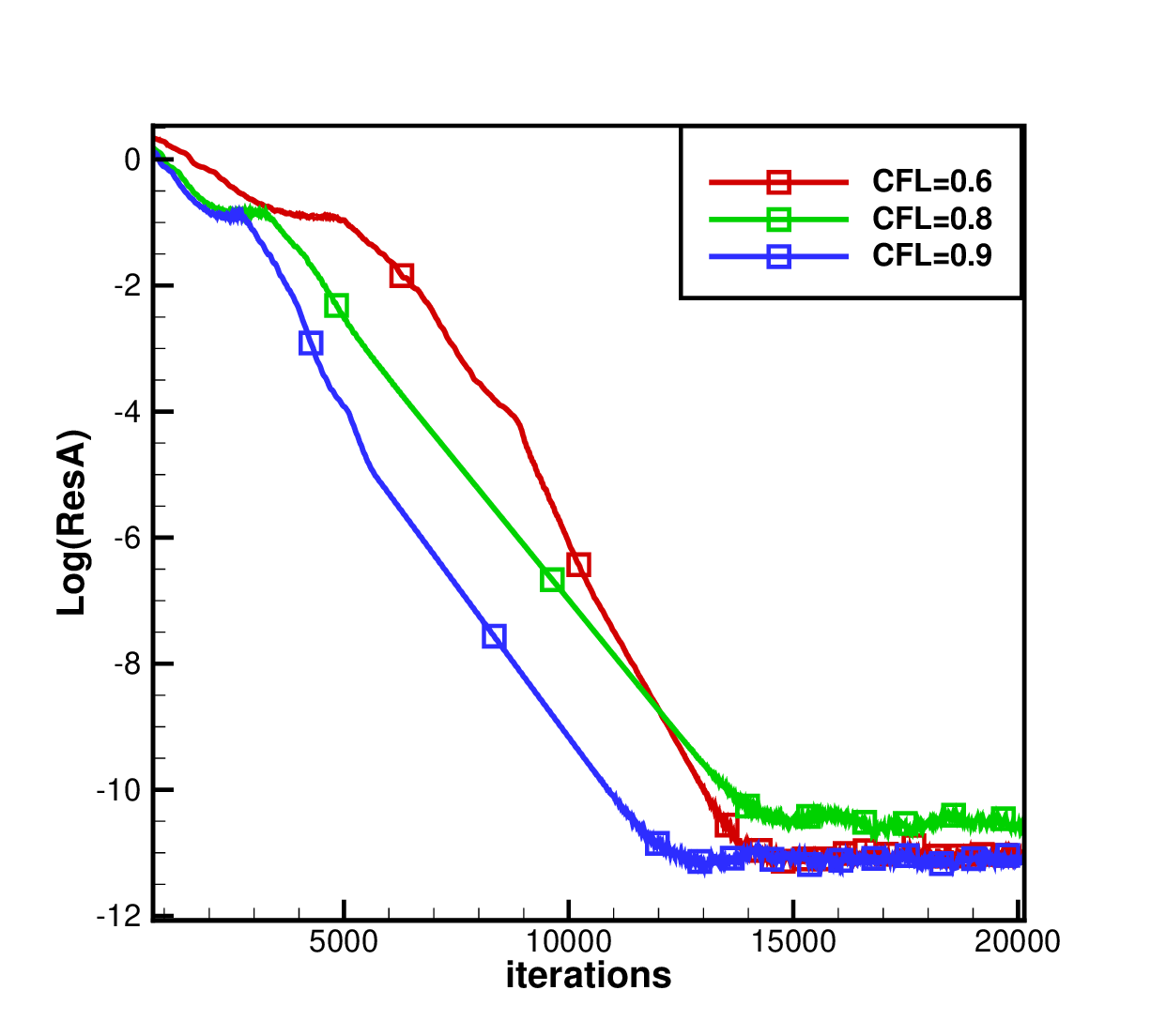}
\end{minipage}}
\subfigure[Density]
{\begin{minipage}[t]{0.28\linewidth}
\includegraphics[width=2in]{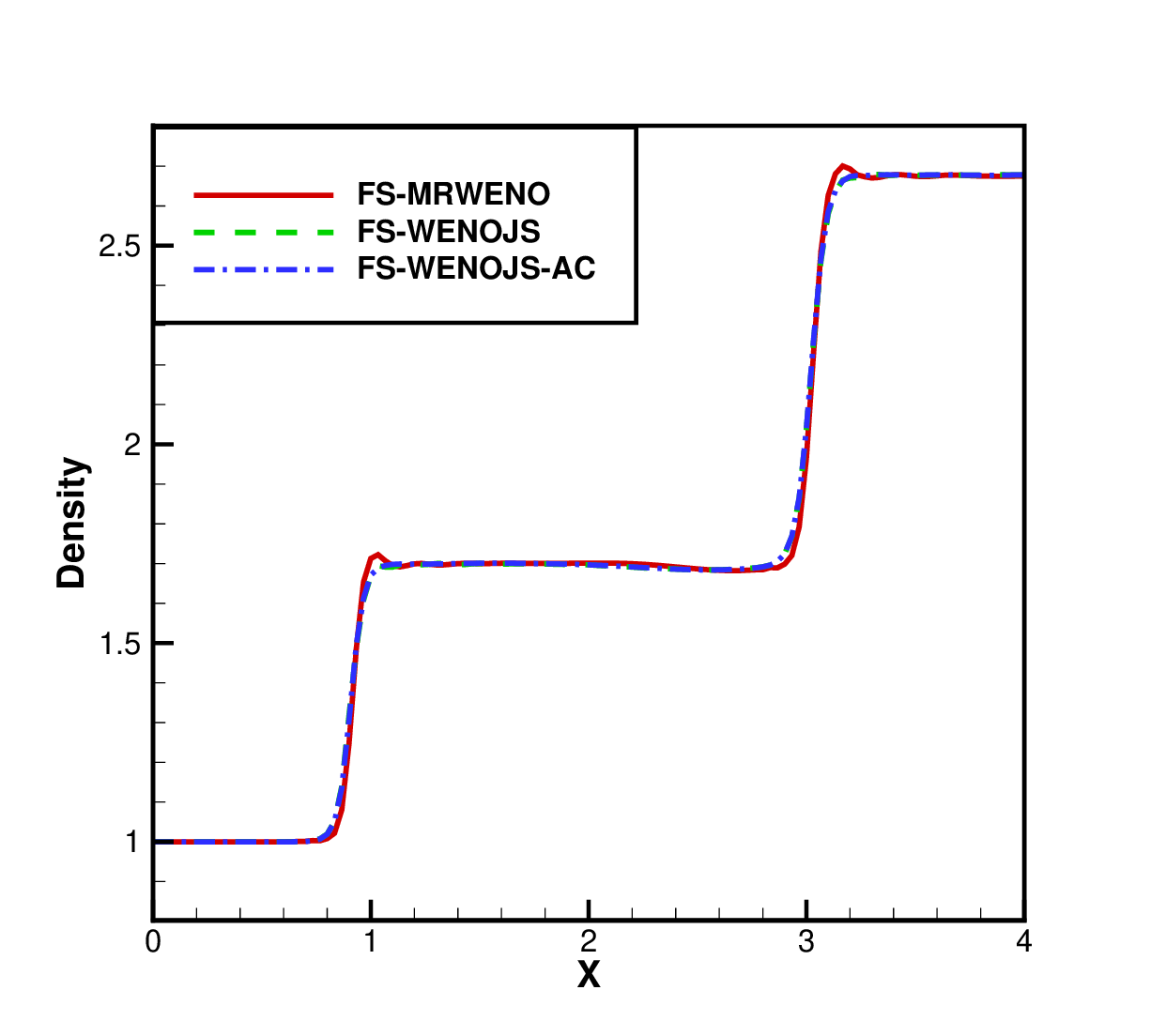}
\end{minipage}}
\caption{\label{1.3}Example 2: The convergence history of the residue as a function of number of iterations for four schemes. Density distributions along the line $y = 0.5$.}
\end{figure}

\noindent This is a classic case often used to test the convergence of a numerical scheme to a steady state. The computational domain is $[0,4]\times[0,1]$ and is discretized using a $120\times30$ grid. The boundary conditions are applied as follows. The bottom edge is set to a reflective wall condition, and the right boundary to supersonic outflow. Dirichlet conditions are imposed on the left boundary with $(\rho,u,v,p)=(1.0,2.9,0,5/7)$, and on the top boundary with $(\rho,u,v,p)=(1.69997,2.61934,-0.50632,1.52819)$. Initially, we set the solution in the entire domain to be that at the left boundary. It is difficult for the residue of many high order WENO schemes to converge to the level of round off errors \cite{WuLiang, SSCW, SCW}. Figure \ref{1.3} displays the residue evolution plots for four different WENO schemes applied to this problem. It can be observed that the FS-WENOJS scheme fails to achieve full convergence even with a very small CFL number, whereas both the FS-MRWENO and FS-WENOJS-AC schemes successfully attain full convergence. Notably, the FS-WENOJS-AC scheme allows for a larger CFL condition number compared to MR-WENO. Figure \ref{1.3}(e) presents the density profile along $y=0.5$. Although the FS-MRWENO scheme achieves full convergence in Figure \ref{1.3} and has similar contour plots in Figure \ref{1.2}, it can be seen from Figure \ref{1.3}(e) that there are weak numerical oscillations near the shock wave. In contrast, the newly proposed scheme does not exhibit this issue. Table \ref{1.1} presents the number of iterations and CPU time required for both schemes under different CFL numbers, up to the one that yields the fastest convergence. It can be observed that for this problem, the newly proposed scheme can achieve convergence with a CFL number as high as 1, whereas the MRWENO scheme only reaches a maximum CFL number of 0.6. A comparison of CPU time shows that the FS-WENOJS-AC scheme requires only 15.42 seconds to achieve the fastest convergence, reducing computational time by more than half compared to the FS-MRWENO scheme. It is worth noting that although the FS-WENOJS-AC scheme has a faster computational rate, it requires more iterations. This is because, after fixing the nonlinear weights, the computational cost per iteration drops significantly, which ultimately results in less total CPU time. Figure \ref{1.2} presents the density contours for the three numerical schemes. Although their final residues and computational times differ, the contour plots closely resemble one another. It should be additionally noted that the residue of the RK-WENOJS-AC scheme can only drop to a level of $10^{-11}$. Therefore, in this example, its computational efficiency is no longer compared with that of the other schemes.
\begin{table}
		\centering
\begin{tabular}{|c|c|c|}\hline
			\multicolumn{3}{|c|}{FS-MRWENO}\\\hline
            $\gamma:$ CFL number & iteration number & CPU time \\\hline
0.5	&2170	&40.22\\\hline
0.6	&1934	&34.14\\\hline
			\multicolumn{3}{|c|}{FS-WENOJS-AC}\\\hline
            $\gamma:$ CFL number & iteration number & CPU time \\\hline
0.6	&5082		&29.92\\\hline
0.8	&3762		&19.89\\\hline
1.0	&2838		&15.42\\\hline
		\end{tabular}
		\caption{\label{1.1}Example 2: Number of iterations, the final time, and the total CPU time of two different iterative schemes when convergence is obtained. Convergence
criterion threshold value is $10^{-12}$. CPU time unit: second.}
	\end{table}

\bigskip
\noindent{\bf Example 3. Supersonic flow past one short plate with an attack angle}

\begin{figure}
		\centering
\subfigure[FS-MRWENO]
{\begin{minipage}[t]{0.3\linewidth}
\includegraphics[width=2.0in]{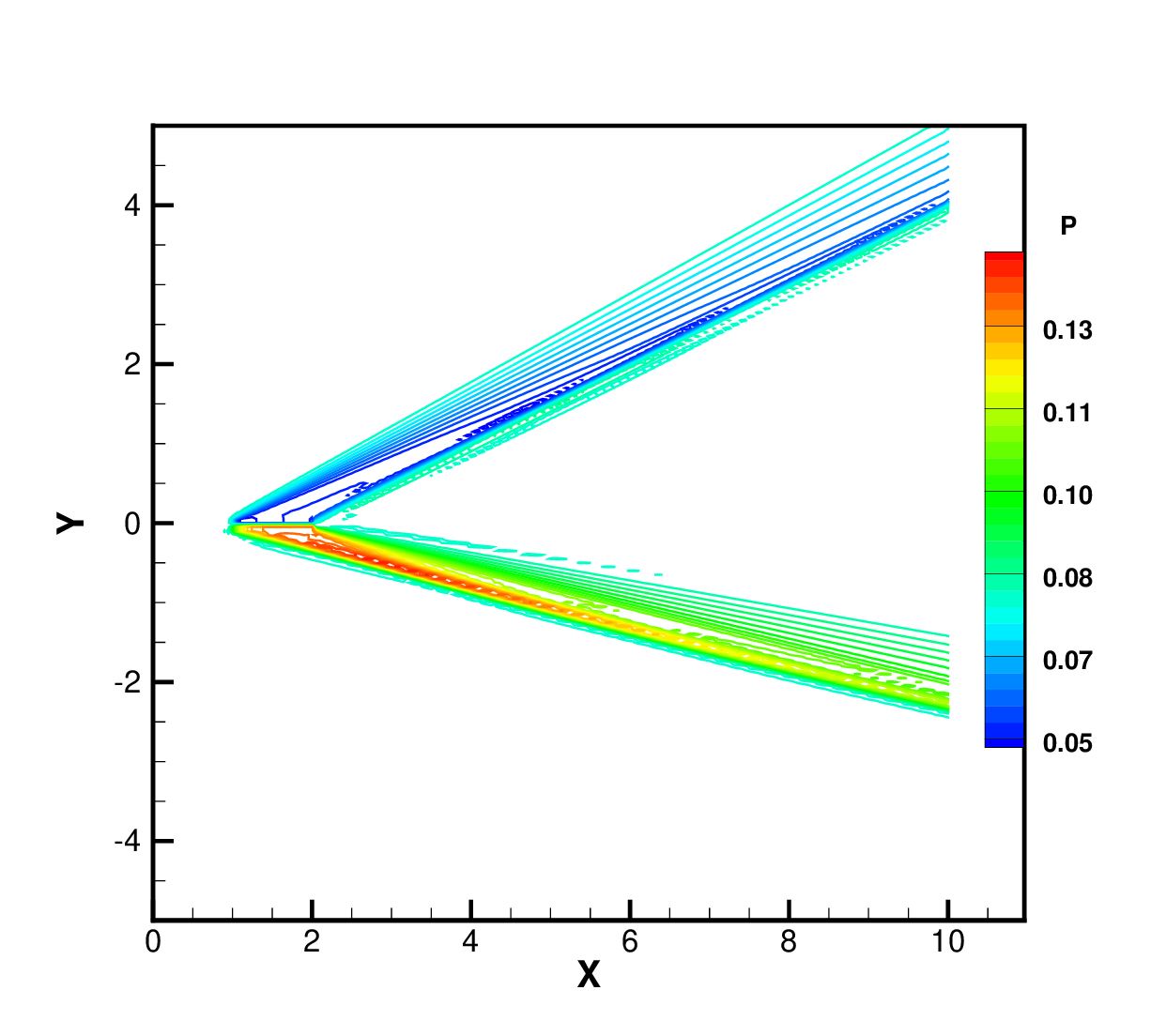}
\end{minipage}}
\subfigure[RK-WENOJS-AC]
{\begin{minipage}[t]{0.3\linewidth}
\includegraphics[width=2.0in]{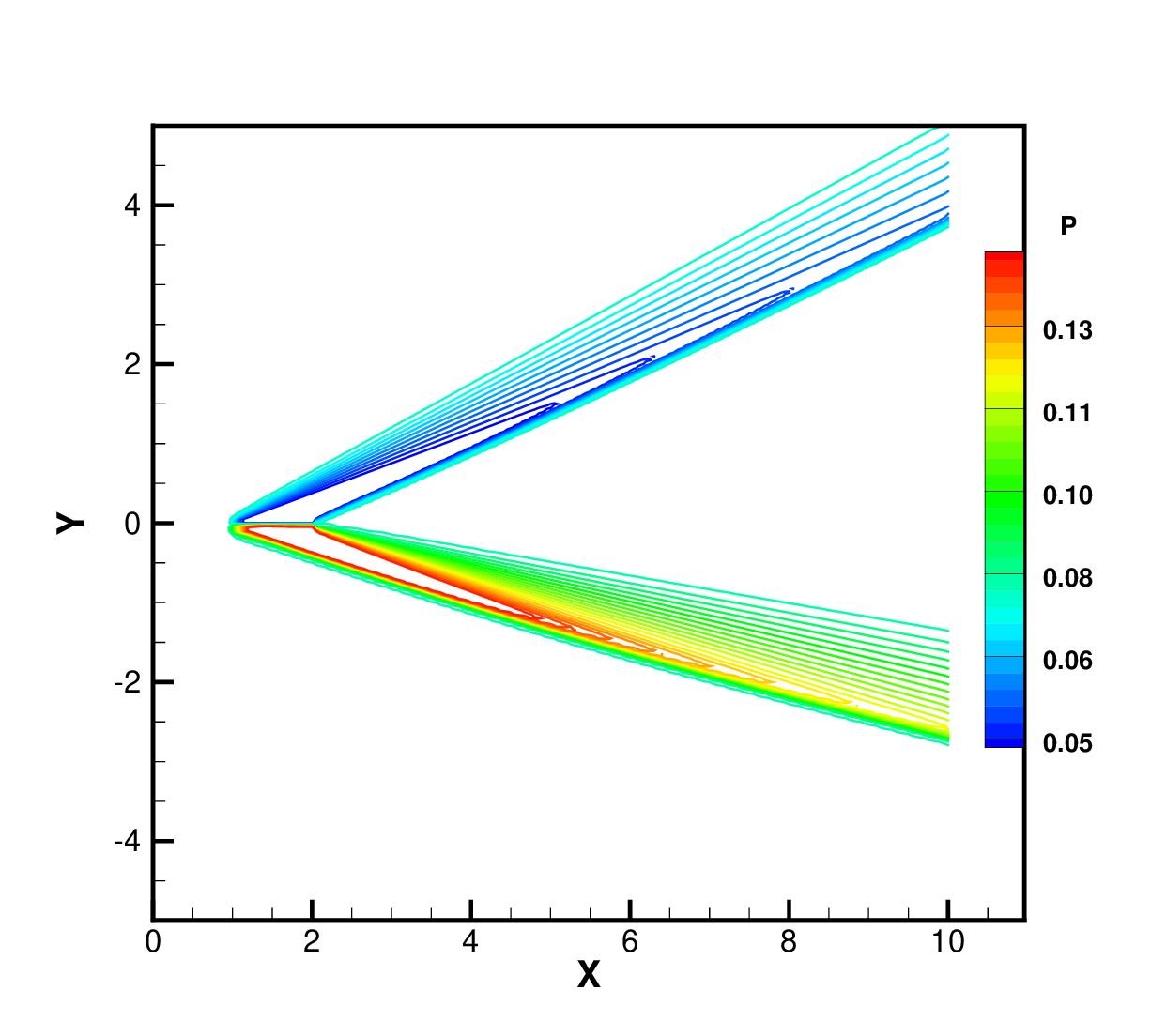}
\end{minipage}}
\subfigure[FS-WENOJS-AC]
{\begin{minipage}[t]{0.3\linewidth}
\includegraphics[width=2.0in]{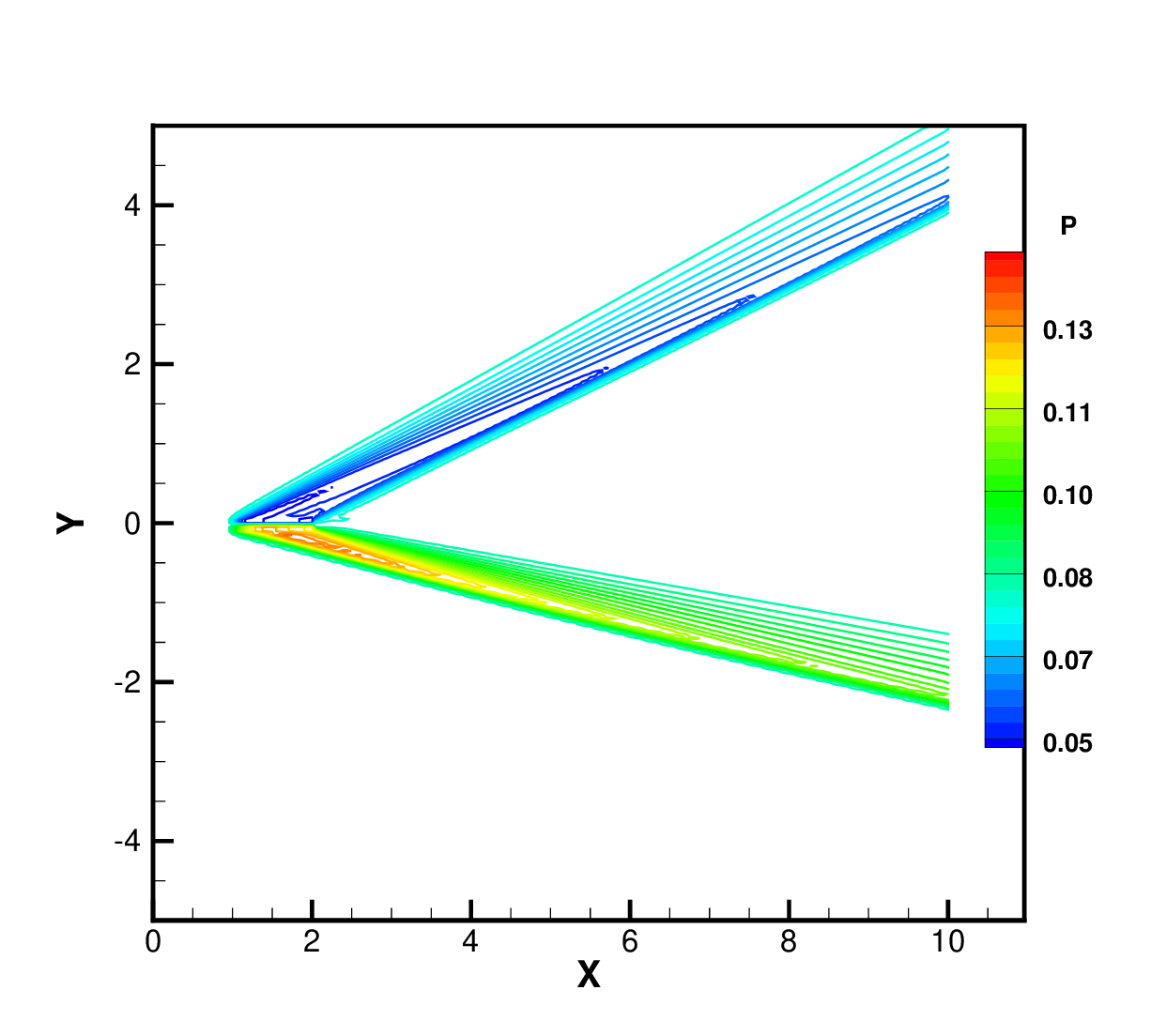}
\end{minipage}}
\caption{\label{2.2}Example 3: Thirty equally spaced pressure contours from 0.05 to 0.14 of the converged steady states of numerical solutions by three WENO schemes.}
\end{figure}

\begin{figure}
		\centering
\subfigure[FS-MRWENO]
{\begin{minipage}[t]{0.23\linewidth}
\includegraphics[width=1.7in]{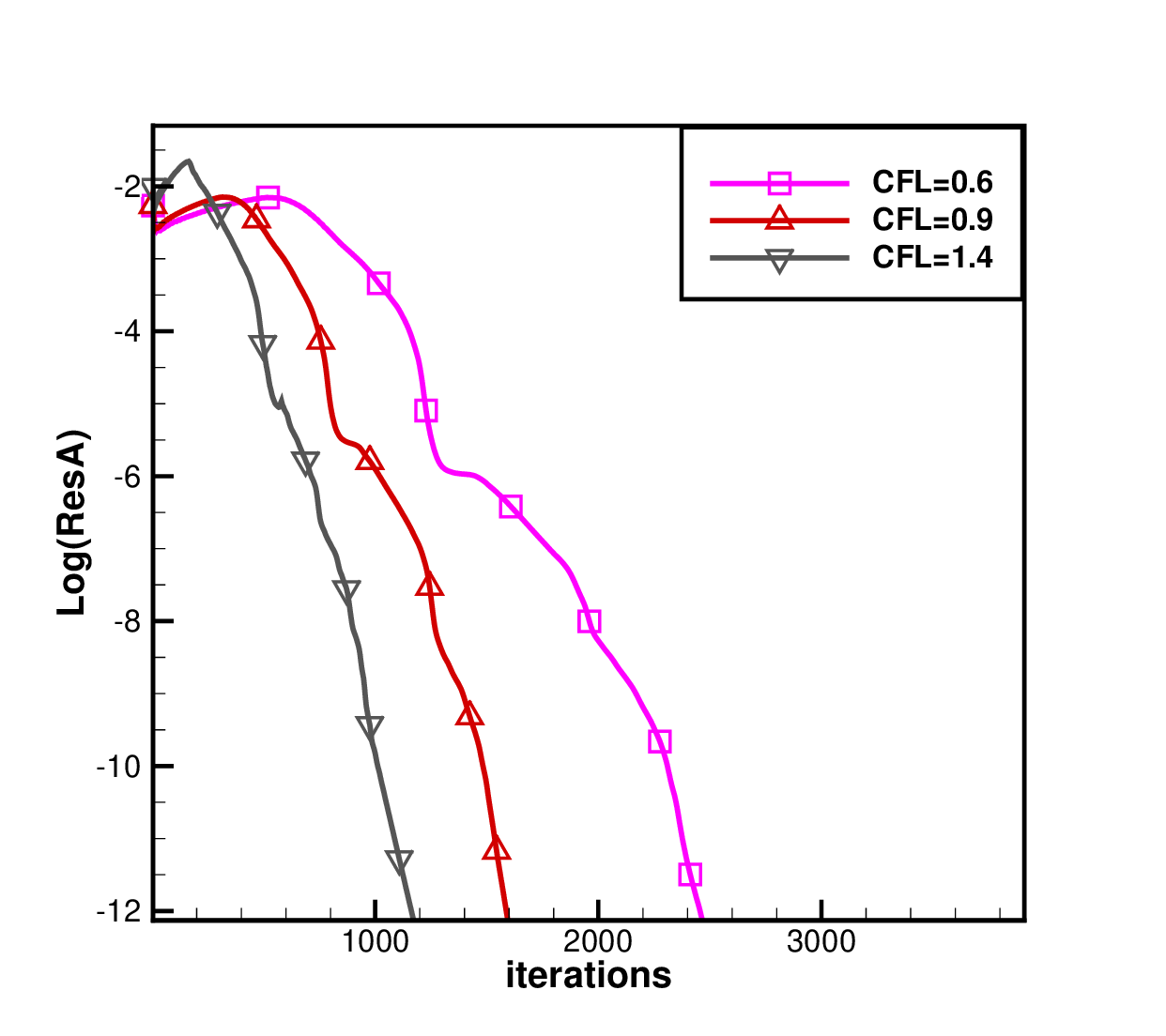}
\end{minipage}}
\subfigure[FS-WENOJS]
{\begin{minipage}[t]{0.23\linewidth}
\includegraphics[width=1.7in]{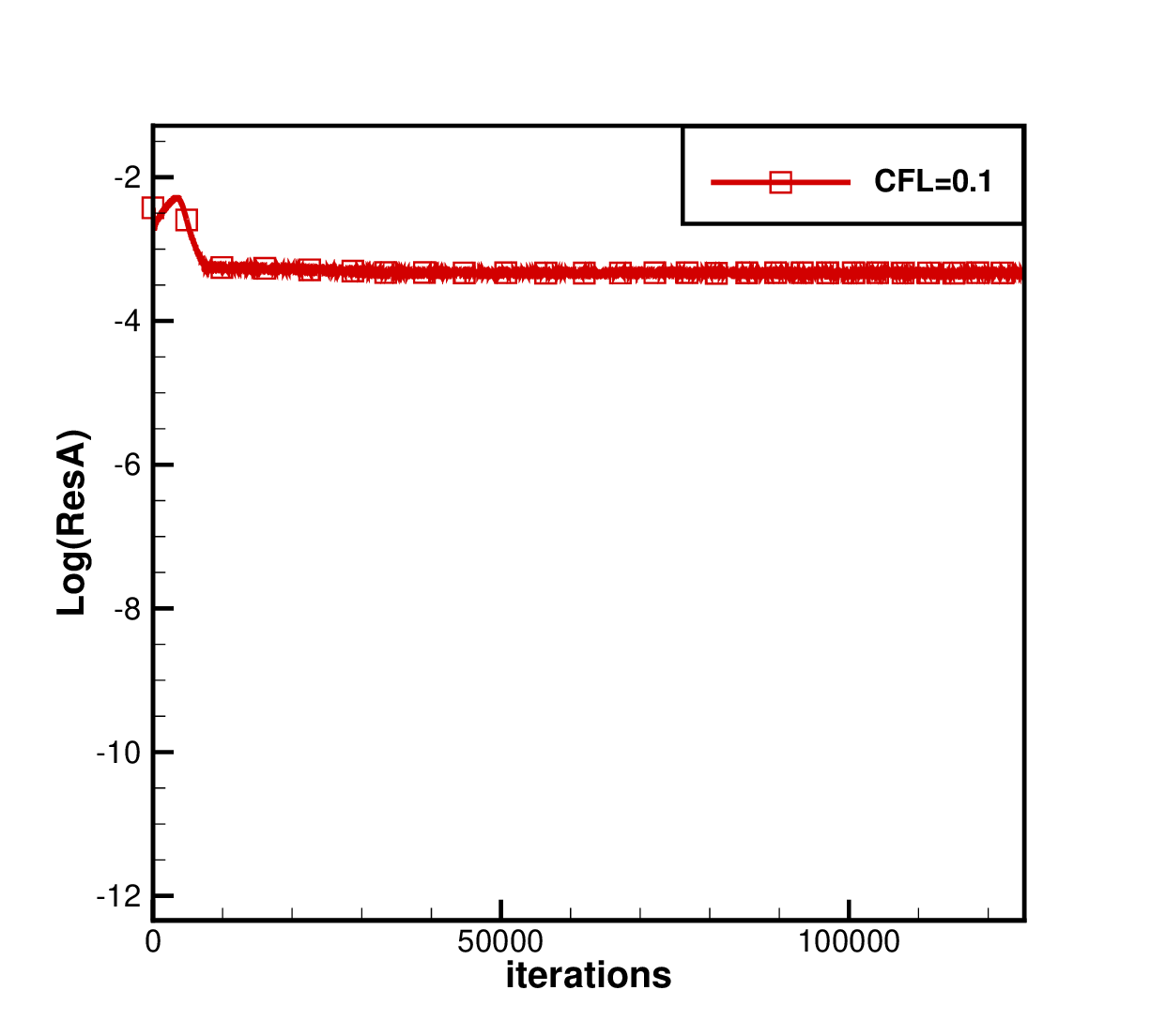}
\end{minipage}}
\subfigure[RK-WENOJS-AC]
{\begin{minipage}[t]{0.23\linewidth}
\includegraphics[width=1.7in]{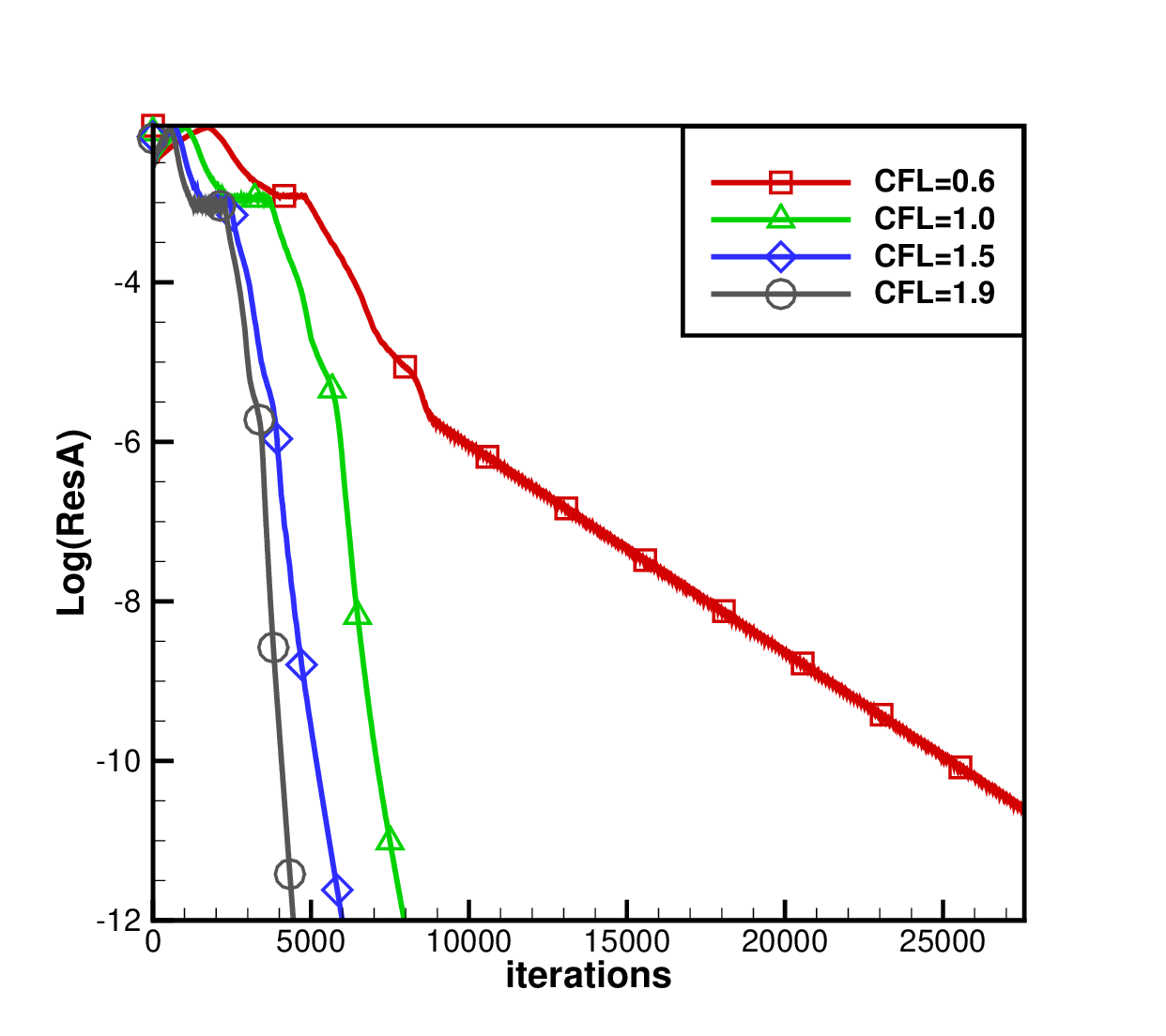}
\end{minipage}}
\subfigure[FS-WENOJS-AC]
{\begin{minipage}[t]{0.23\linewidth}
\includegraphics[width=1.7in]{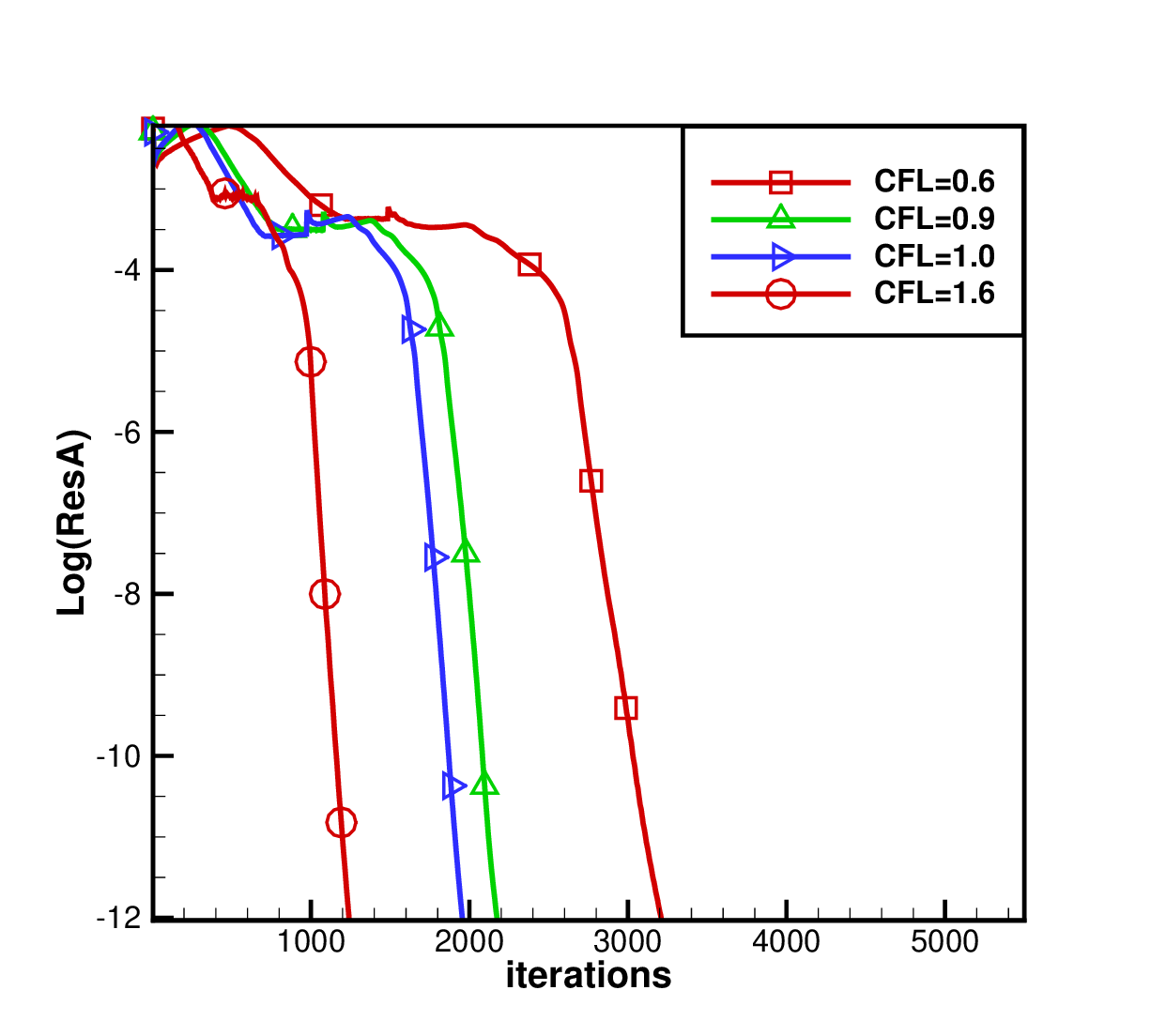}
\end{minipage}}
\caption{Example 3: The convergence history of the residue as a function of number of iterations for four schemes.}\label{2.3}
\end{figure}

\noindent The example \cite{JUNZ3} presented illustrates a scenario where supersonic flow passes over a short plate set at an angle of attack, $\alpha=10^{\circ}$. The free-stream Mach number is specified as $M_{\infty}$=3. An ideal gas approaches the plate, which is situated within the range $x \in [1,2]$ along the line $y=0$, flowing from left to right. The initial conditions are defined by $p=\frac{1}{\gamma M_{\infty}^{2}}$, $\rho$=1, $u=\cos(\alpha)$ and $v=\sin(\alpha)$ with $\gamma=1.4$. The computational domain spans  $[0,10]\times [-5,5]$, discretized into a 200$\times$200 grid. A slip boundary condition is applied along the plate, while appropriate physical values are enforced for inflow and outflow boundaries. As observed in Figure \ref{2.3}, when the FS-WENOJS scheme is applied to this problem, its residue still fails to decrease to the level of machine zero. Even with a very small CFL number of 0.1, the iteration residue stagnates at around $10^{-3.3}$ and cannot decrease further. In contrast, the other three schemes successfully achieve convergence. The computational performance of those schemes under varying CFL numbers is summarized in Table \ref{2.1}. It can be observed that the RK-WENOJS-AC scheme can tolerate a larger CFL number to achieve convergence, yet its fastest computational speed remains twice as slow as that of the FS-WENOJS-AC scheme. In contrast, the FS-MRWENO scheme requires a smaller CFL number and a longer convergence time, with its convergence time being approximately three times that of FS-WENOJS-AC scheme. The pressure contours of the three schemes are shown in Figure \ref{2.2}, and they exhibit similar behavior.

\begin{table}
		\centering
\begin{tabular}{|c|c|c|}\hline
			\multicolumn{3}{|c|}{FS-MRWENO}\\\hline
            $\gamma:$ CFL number & iteration number  & CPU time \\\hline
0.6	&2460		&432.94\\\hline
0.9	&1592		&319.33\\\hline
1.4	&1164		&232.50\\\hline
			\multicolumn{3}{|c|}{RK-WENOJS-AC}\\\hline
            $\gamma:$ CFL number & iteration number  & CPU time \\\hline
1.0	&7953		&230.66\\\hline
1.5	&5997		&181.38\\\hline
1.9	&4446		&137.09\\\hline
			\multicolumn{3}{|c|}{FS-WENOJS-AC}\\\hline
            $\gamma:$ CFL number & iteration number  & CPU time \\\hline
0.6	&3212		&238.38\\\hline
0.9	&2176		&187.77\\\hline
1.0	&1956		&165.89\\\hline
1.5	&1408		&88.61\\\hline
1.6	&1242		&71.81\\\hline
		\end{tabular}
		\caption{\label{2.1}Example 3: Number of iterations, the final time, and the total CPU time of three different iterative schemes when convergence is obtained. Convergence
criterion threshold value is $10^{-12}$. CPU time unit: second.}
	\end{table}

\bigskip
\noindent{\bf Example 4. Supersonic flow past two short plates with an attack angle}

\begin{figure}
		\centering
\subfigure[FS-MRWENO]
{\begin{minipage}[t]{0.3\linewidth}
\includegraphics[width=2.0in]{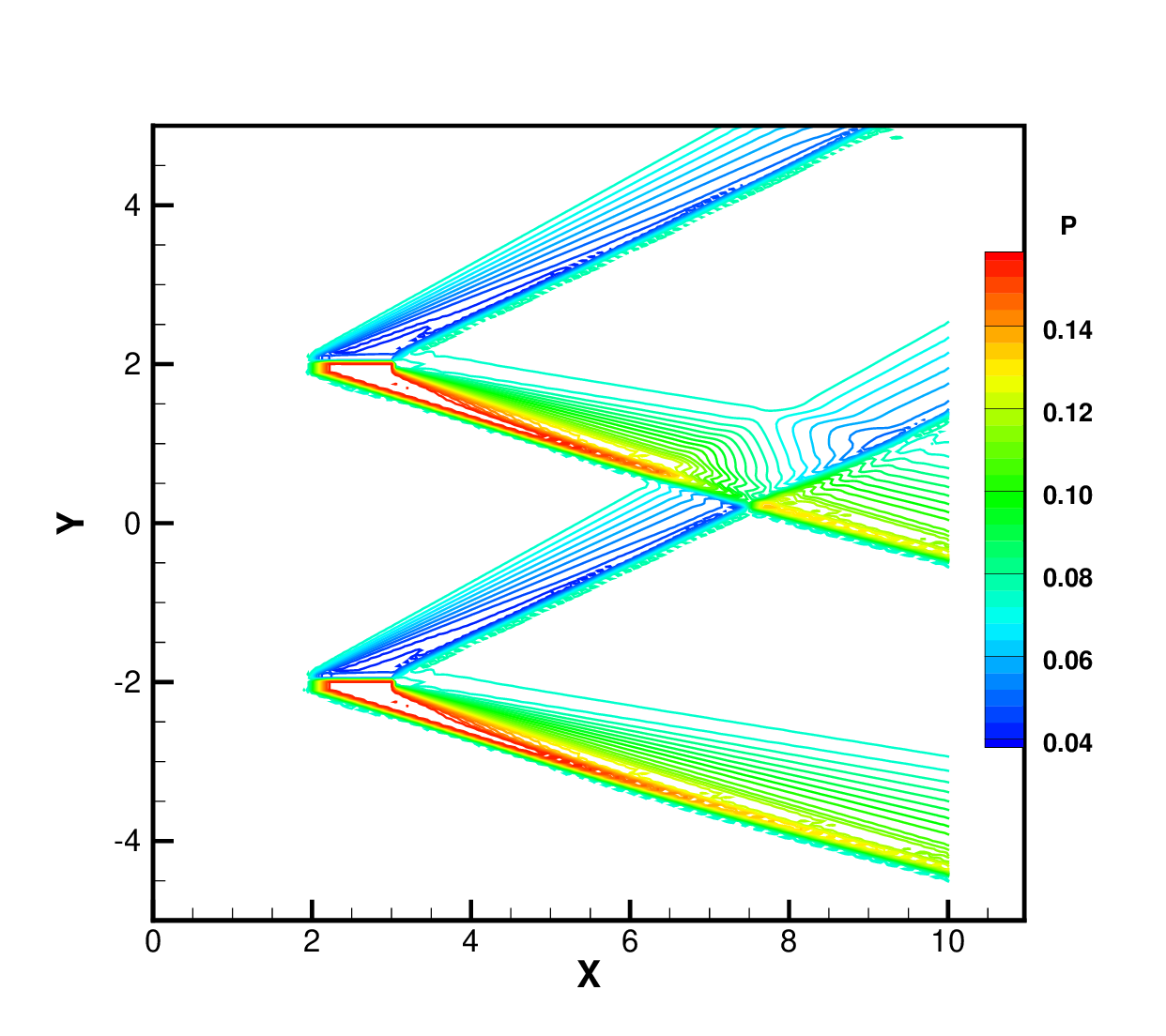}
\end{minipage}}
\subfigure[RK-WENOJS-AC]
{\begin{minipage}[t]{0.3\linewidth}
\includegraphics[width=2.0in]{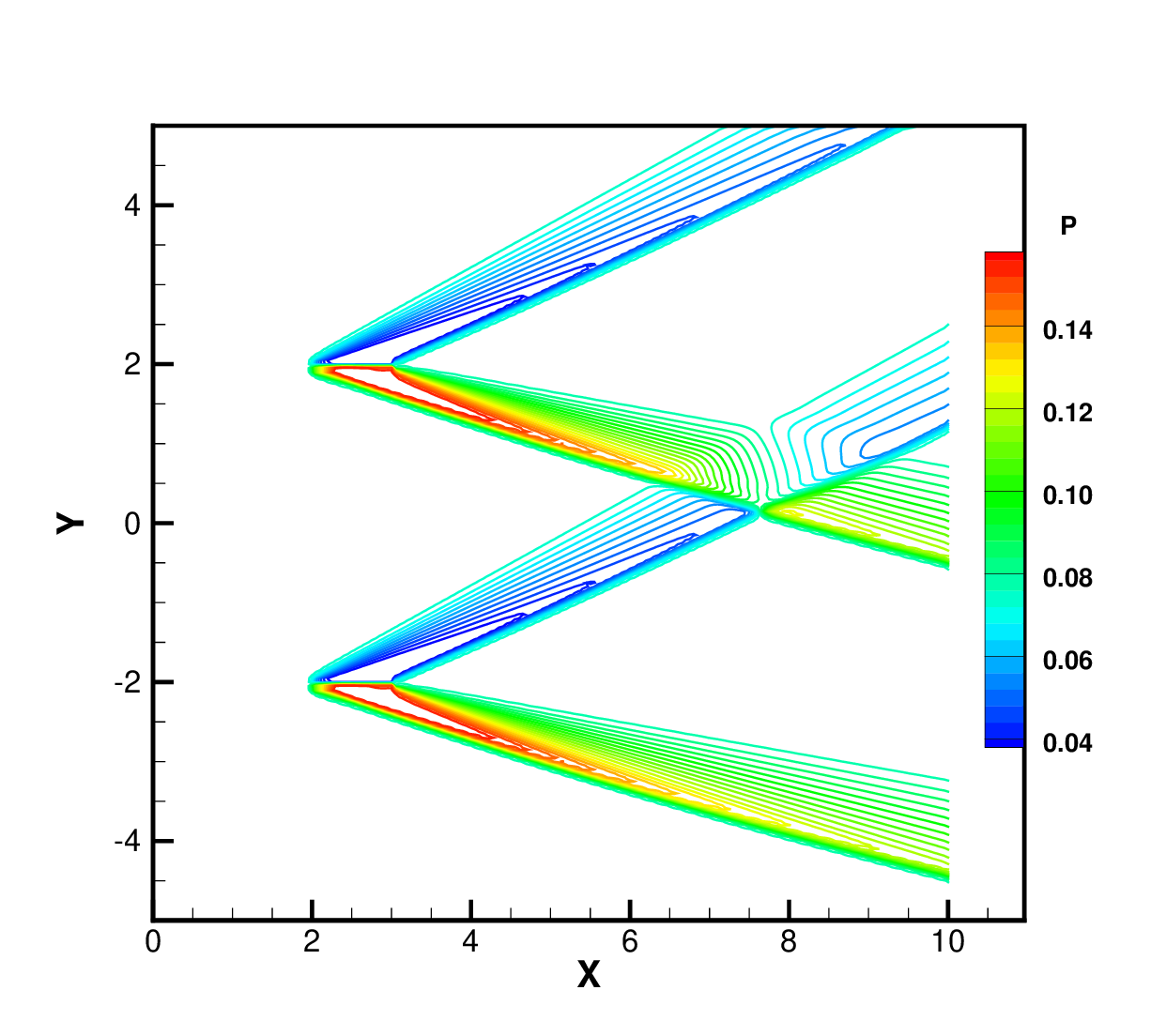}
\end{minipage}}
\subfigure[FS-WENOJS-AC]
{\begin{minipage}[t]{0.3\linewidth}
\includegraphics[width=2.0in]{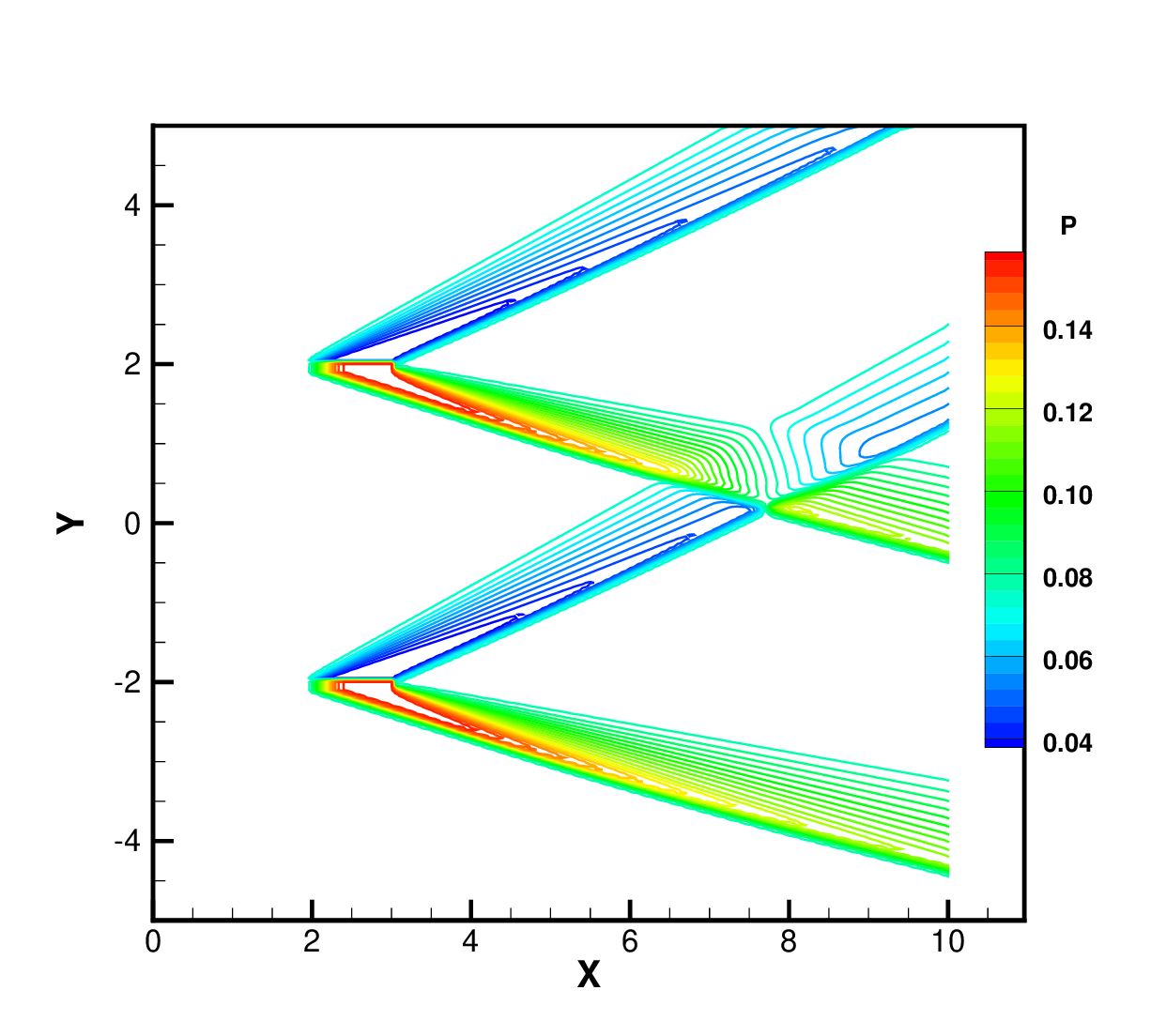}
\end{minipage}}
\caption{\label{3.2}Example 4: Thirty equally spaced pressure contours from 0.04 to 0.16 of the converged steady states of numerical solutions by three different iterative schemes.}
\end{figure}

\begin{figure}
		\centering
\subfigure[FS-MRWENO]
{\begin{minipage}[t]{0.23\linewidth}
\includegraphics[width=1.7in]{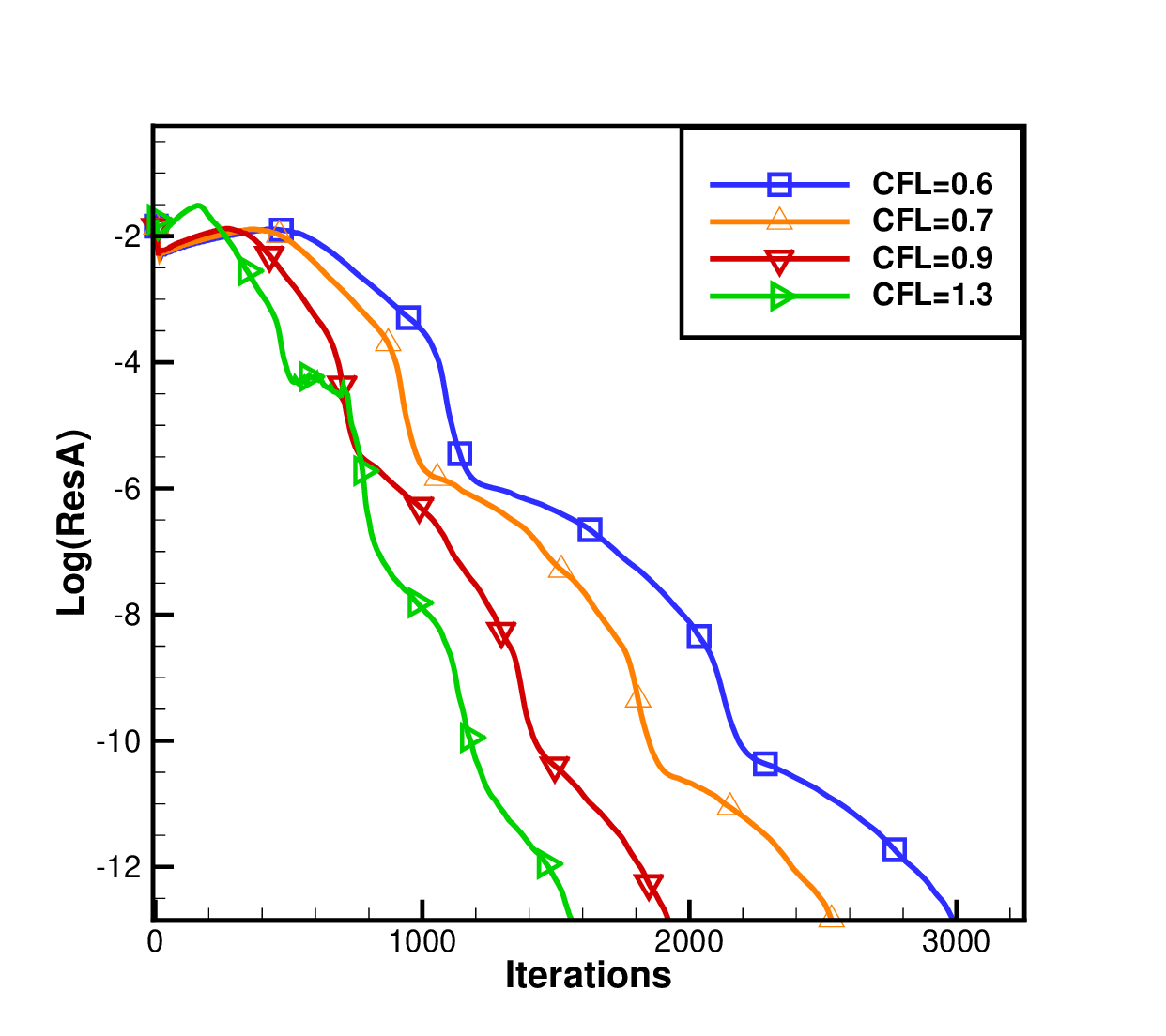}
\end{minipage}}
\subfigure[FS-WENOJS]
{\begin{minipage}[t]{0.23\linewidth}
\includegraphics[width=1.7in]{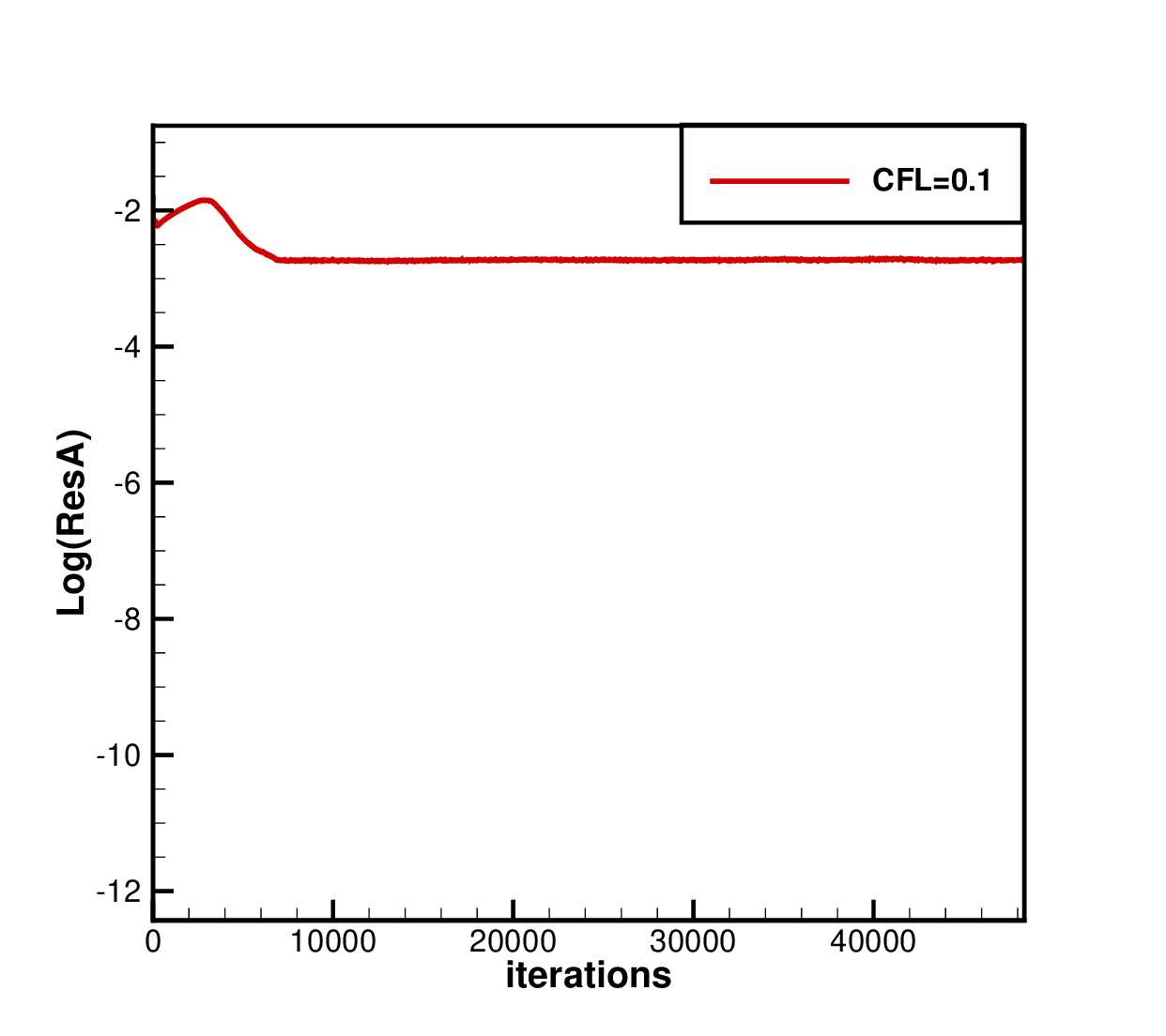}
\end{minipage}}
\subfigure[RK-WENOJS-AC]
{\begin{minipage}[t]{0.23\linewidth}
\includegraphics[width=1.7in]{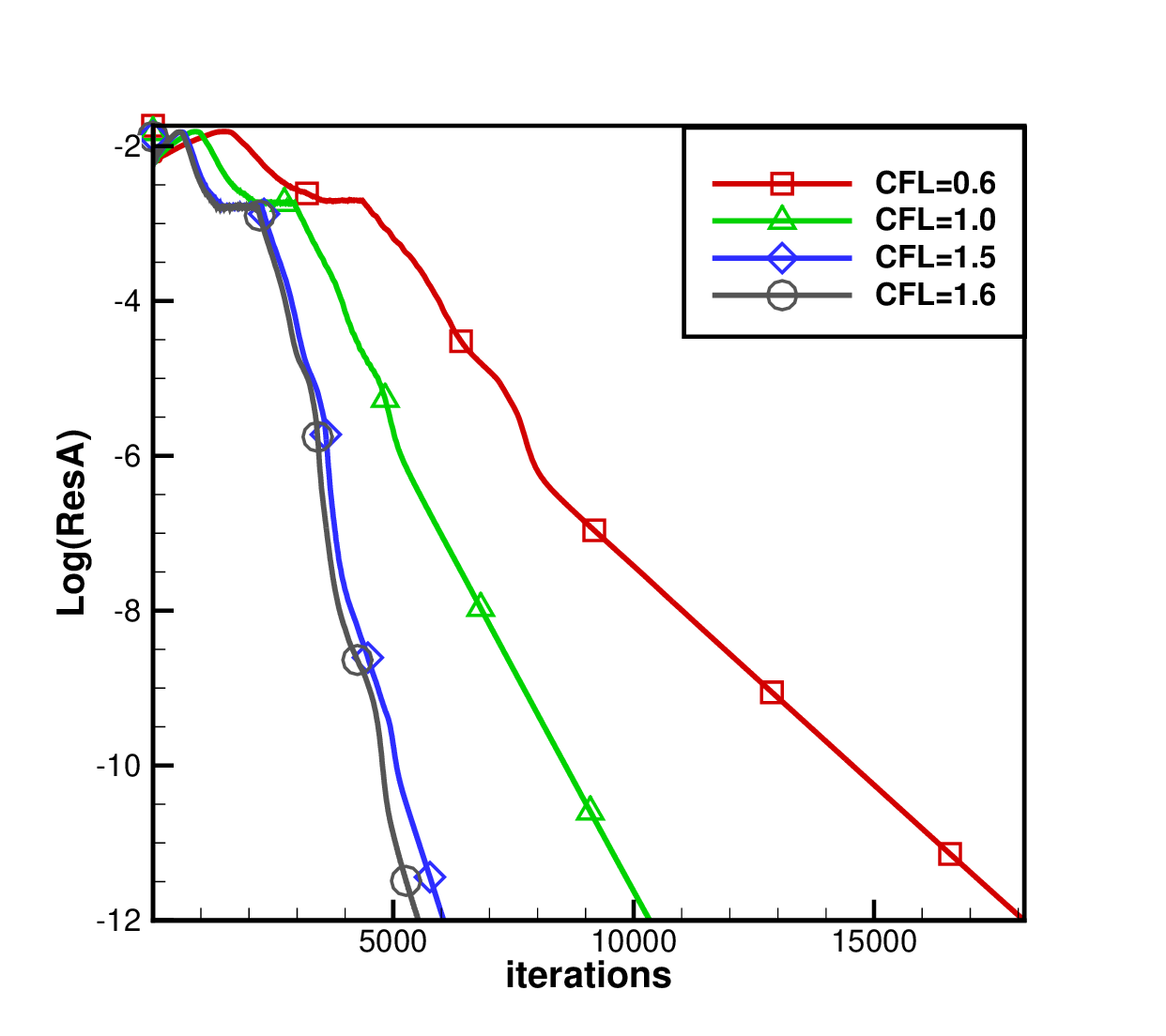}
\end{minipage}}
\subfigure[FS-WENOJS-AC]
{\begin{minipage}[t]{0.23\linewidth}
\includegraphics[width=1.7in]{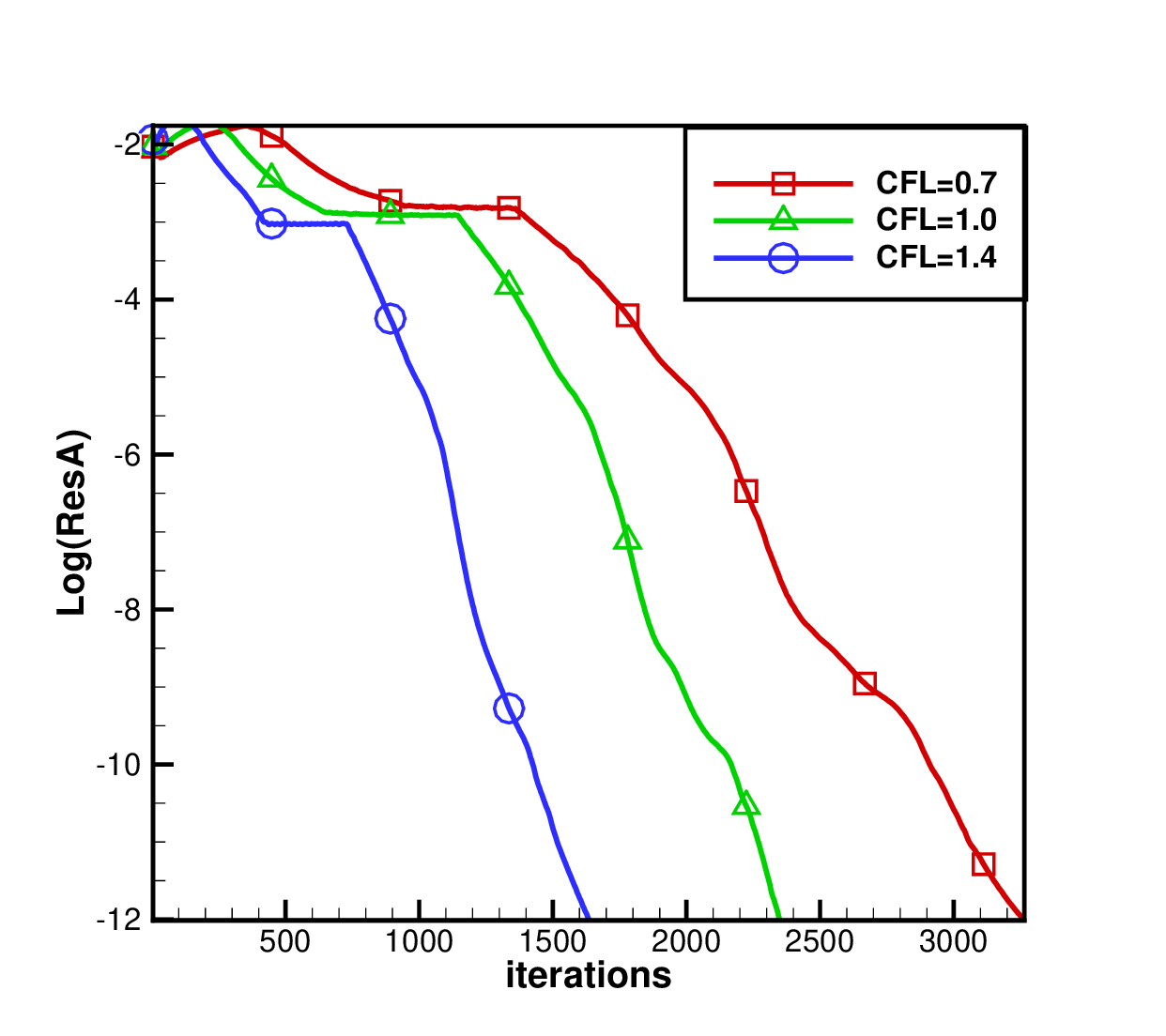}
\end{minipage}}
\caption{\label{3.3}Example 4: The convergence history of the residue as a function of number of iterations for four schemes.}
\end{figure}

\noindent This test case \cite{taijie2} is similar to the previous one, with the only difference being that two short plates are placed within the flow field: one located at $x\in[2, 3], y = -2$, and the other at $x\in [2, 3], y = 2$. The information of the flow field and the boundary conditions are the same as in the previous calculation case. This case involves shock wave interactions, making it more challenging for high-order WENO schemes to converge to machine zero. The residues of the four numerical schemes for this problem are shown in Figure \ref{3.3}. Similarly, and in contrast to the other schemes, the residual for the FS-WENOJS scheme fails to achieve convergence to machine zero. In table \ref{3.1}, number of iterations required to reach the convergence threshold value $10^{-12}$ and total CPU time when convergence is obtained for these three iterative schemes with different CFL numbers $\gamma$ are reported. All three schemes can achieve a large CFL number. The FS-WENOJS-AC scheme remains the most efficient, with its CPU time being approximately $60\%$ lower than that of the FS-MRWENO scheme. The fast sweeping method also saves $30\%$ of the time compared to the Runge-Kutta algorithm. The pressure contour plots of the three schemes are shown in Figure \ref{3.2}, and they also exhibit nearly identical contours.

\begin{table}
		\centering
\begin{tabular}{|c|c|c|}\hline
			\multicolumn{3}{|c|}{FS-MRWENO}\\\hline
            $\gamma:$ CFL number & iteration number  & CPU time \\\hline
0.6	&2836		&491.92\\\hline
0.9	&1820		&368.34\\\hline
1.3	&1476		&299.64\\\hline
			\multicolumn{3}{|c|}{RK-WENOJS-AC}\\\hline
            $\gamma:$ CFL number & iteration number  & CPU time \\\hline
0.6	&18123		&560.41\\\hline
1.0	&10338	&317.20\\\hline
1.5 &6051  &186.27\\\hline
1.6 &5526   &166.39\\\hline
			\multicolumn{3}{|c|}{FS-WENOJS-AC}\\\hline
            $\gamma:$ CFL number & iteration number & CPU time \\\hline
0.7	&3264	&236.83\\\hline
1.0	&2352	&169.17\\\hline
1.4 &1636   &119.05\\\hline
		\end{tabular}
		\caption{\label{3.1}Example 4: Number of iterations, the final time, and the total CPU time of three different iterative schemes when convergence is obtained. Convergence
criterion threshold value is $10^{-12}$. CPU time unit: second.}
	\end{table}

\bigskip
\noindent{\bf Example 5.  Supersonic flow past three long plates with an attack angle}

\begin{figure}
		\centering
\subfigure[FS-MRWENO]
{\begin{minipage}[t]{0.3\linewidth}
\includegraphics[width=2.0in]{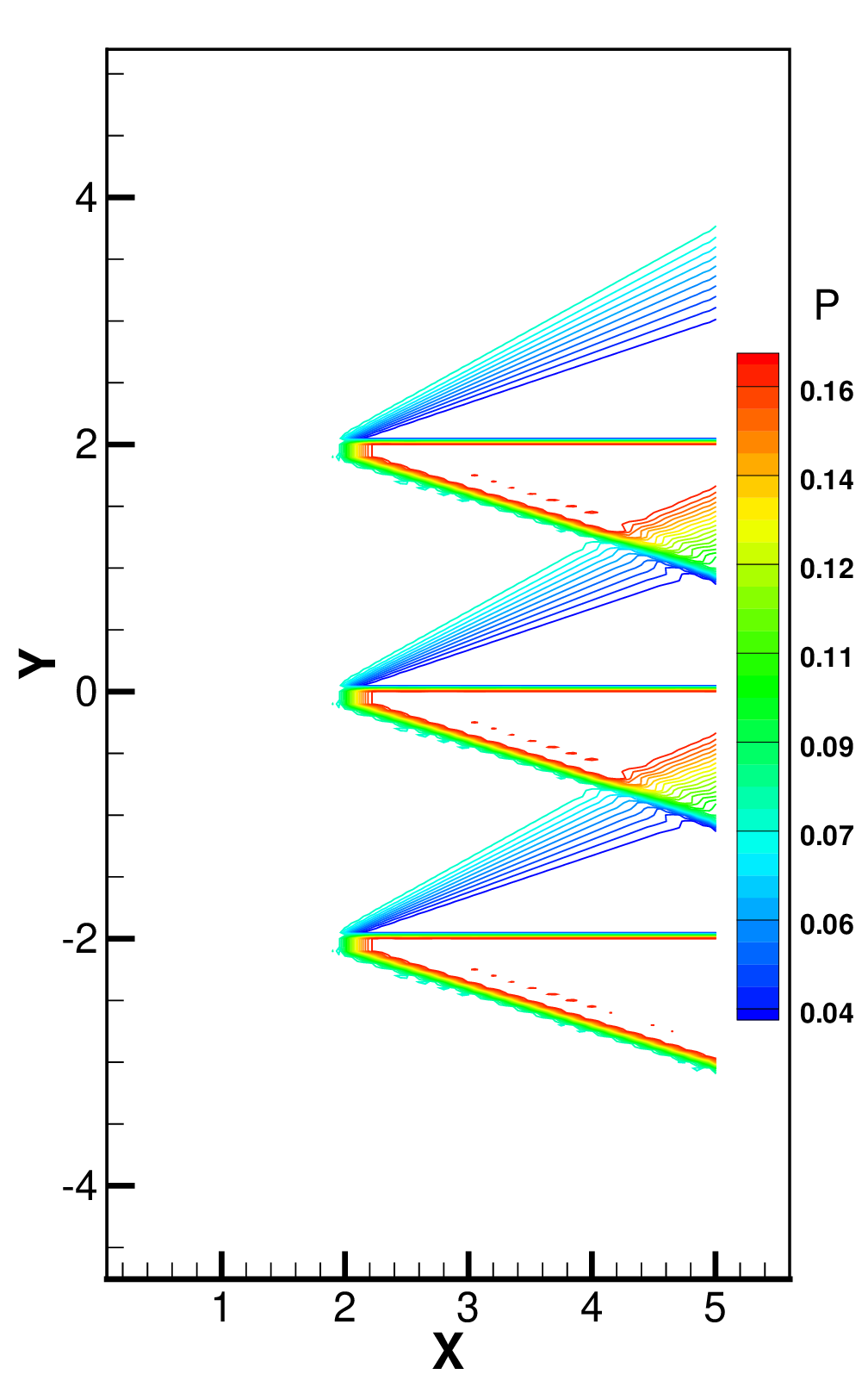}
\end{minipage}}
\subfigure[RK-WENOJS-AC]
{\begin{minipage}[t]{0.3\linewidth}
\includegraphics[width=2.0in]{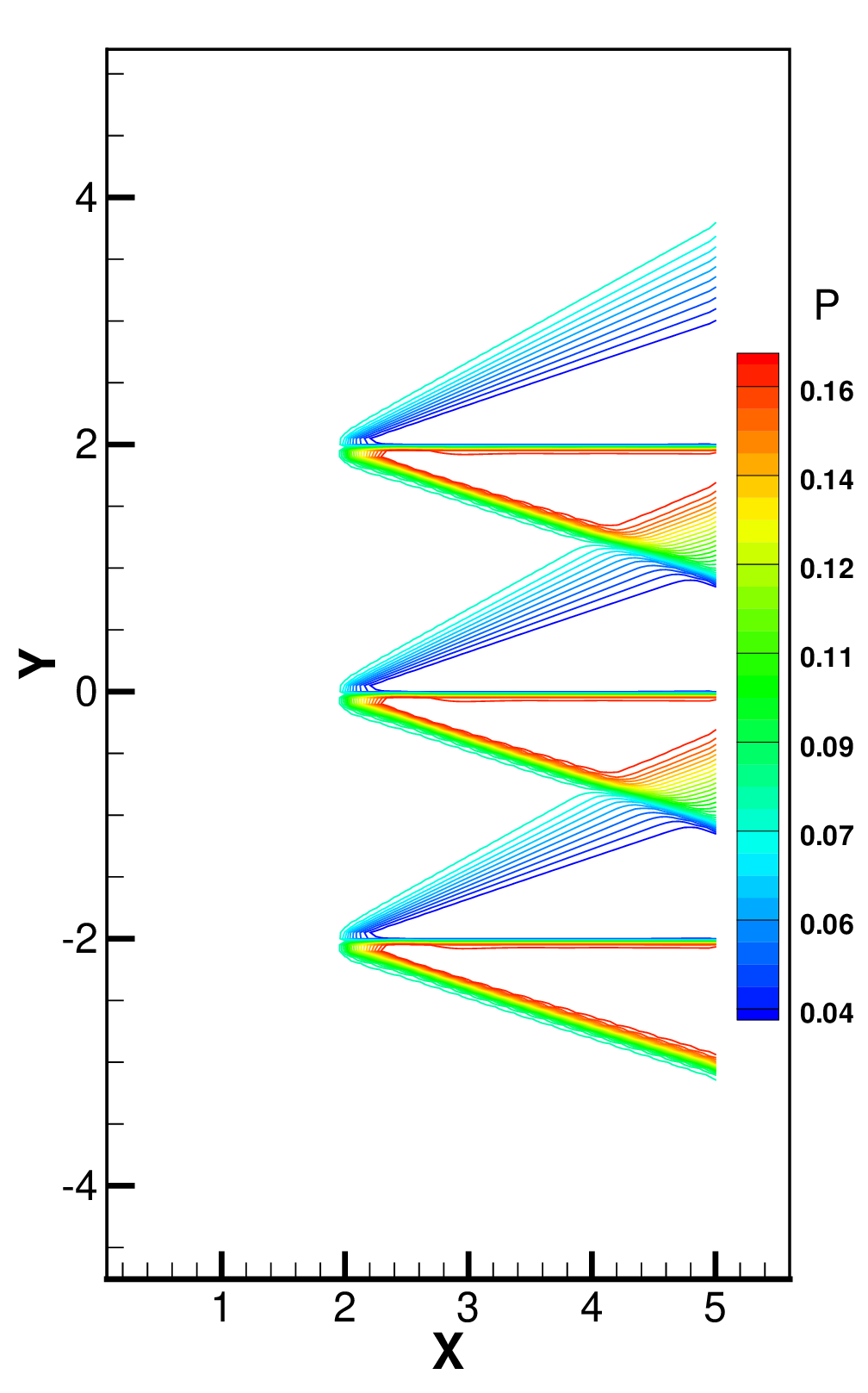}
\end{minipage}}
\subfigure[FS-WENOJS-AC]
{\begin{minipage}[t]{0.3\linewidth}
\includegraphics[width=2.0in]{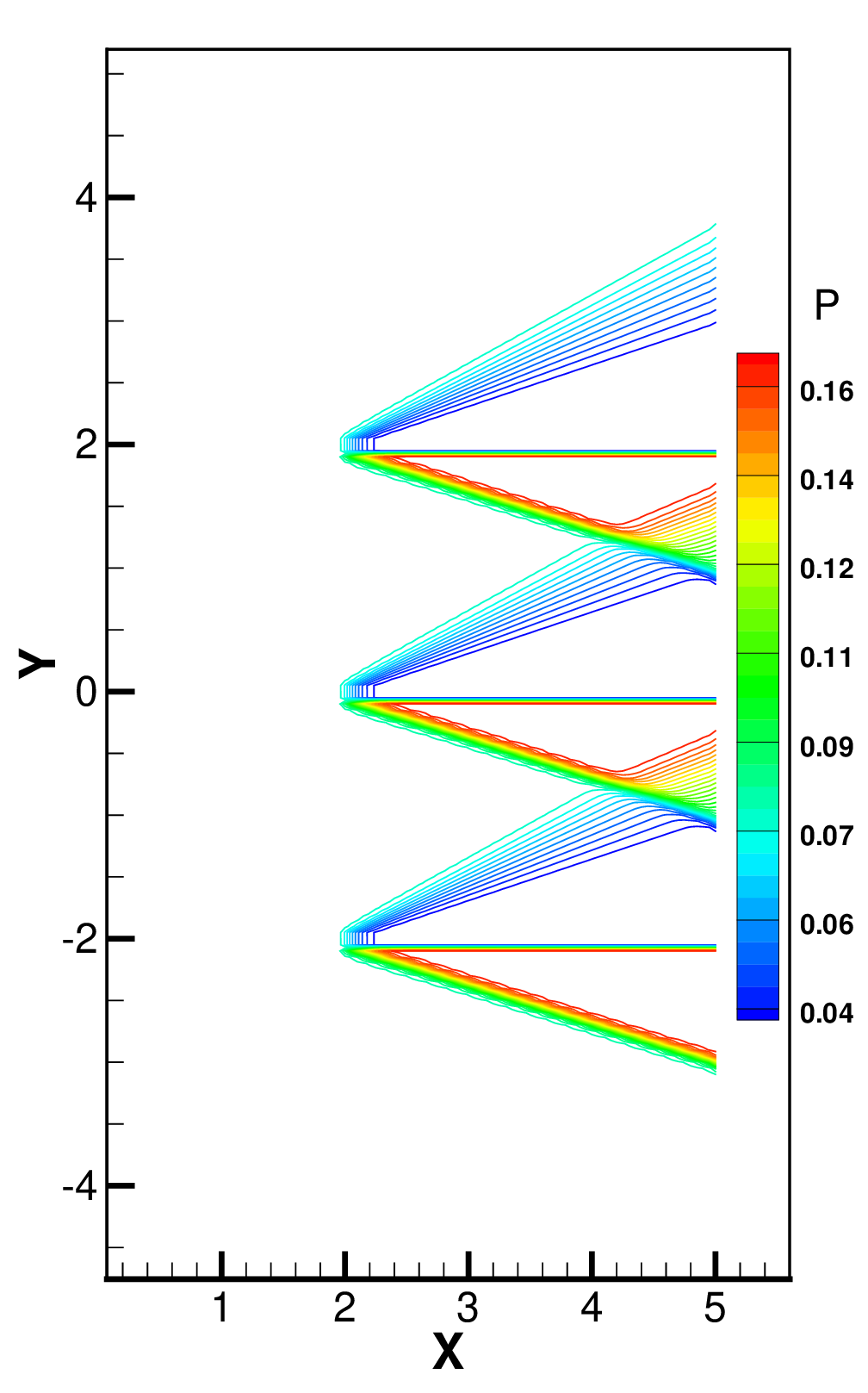}
\end{minipage}}
\caption{\label{4.2}Example 5: Thirty equally spaced pressure contours from 0.04 to 0.17 of the converged steady states of numerical solutions by three different iterative schemes.}
\end{figure}

\begin{figure}
		\centering
\subfigure[FS-MRWENO]
{\begin{minipage}[t]{0.23\linewidth}
\includegraphics[width=1.7in]{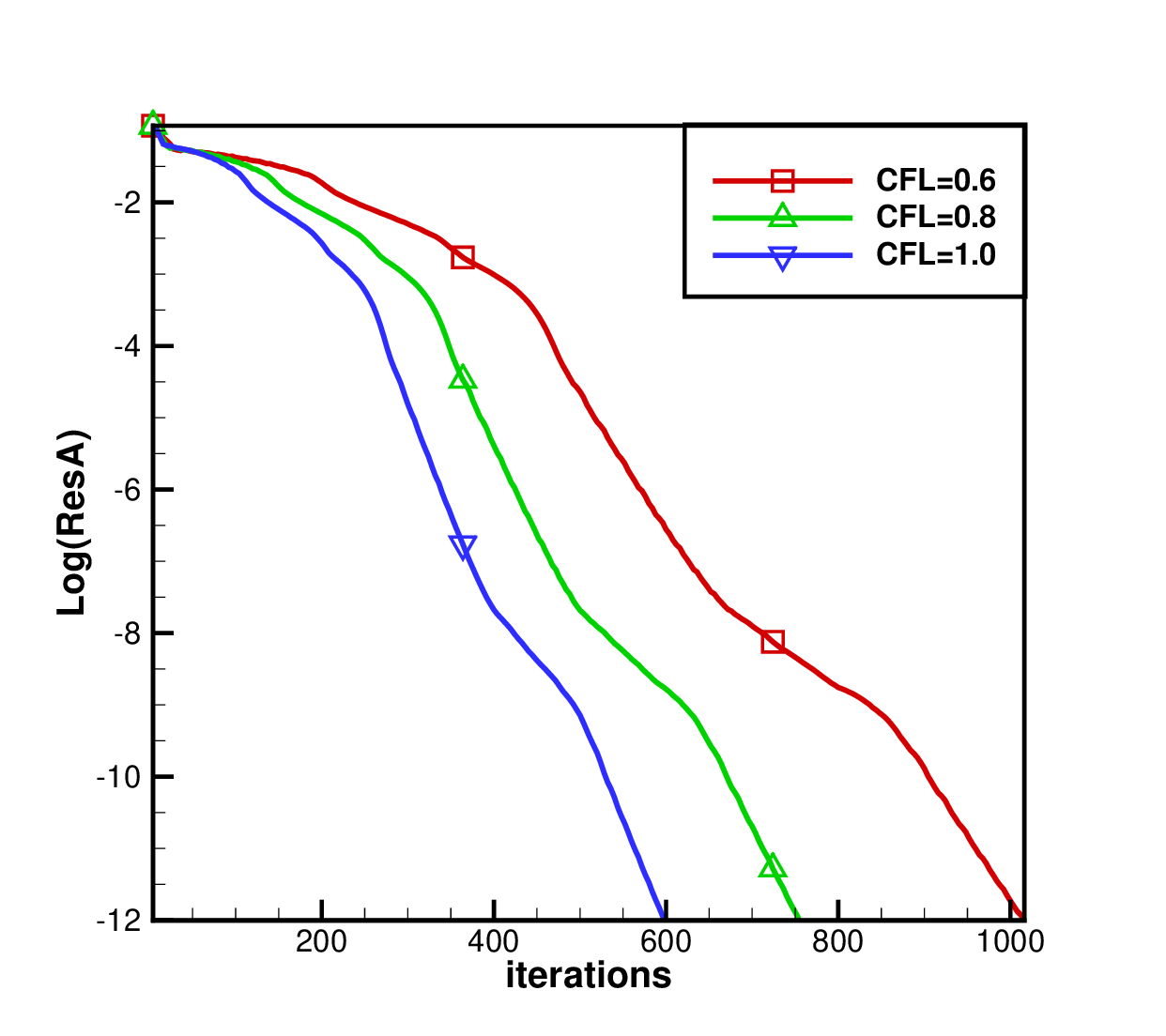}
\end{minipage}}
\subfigure[FS-WENOJS]
{\begin{minipage}[t]{0.23\linewidth}
\includegraphics[width=1.7in]{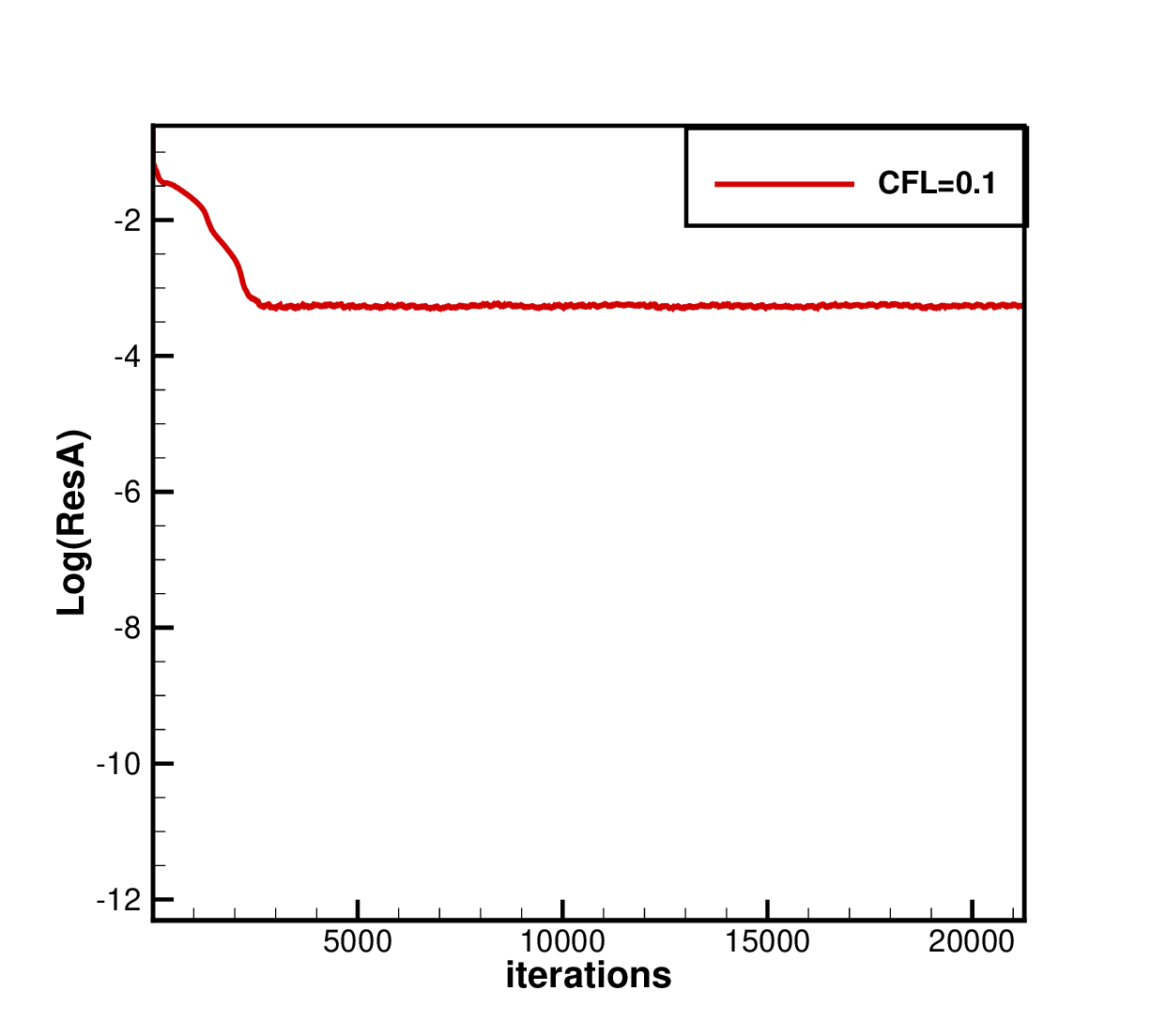}
\end{minipage}}
\subfigure[RK-WENOJS-AC]
{\begin{minipage}[t]{0.23\linewidth}
\includegraphics[width=1.7in]{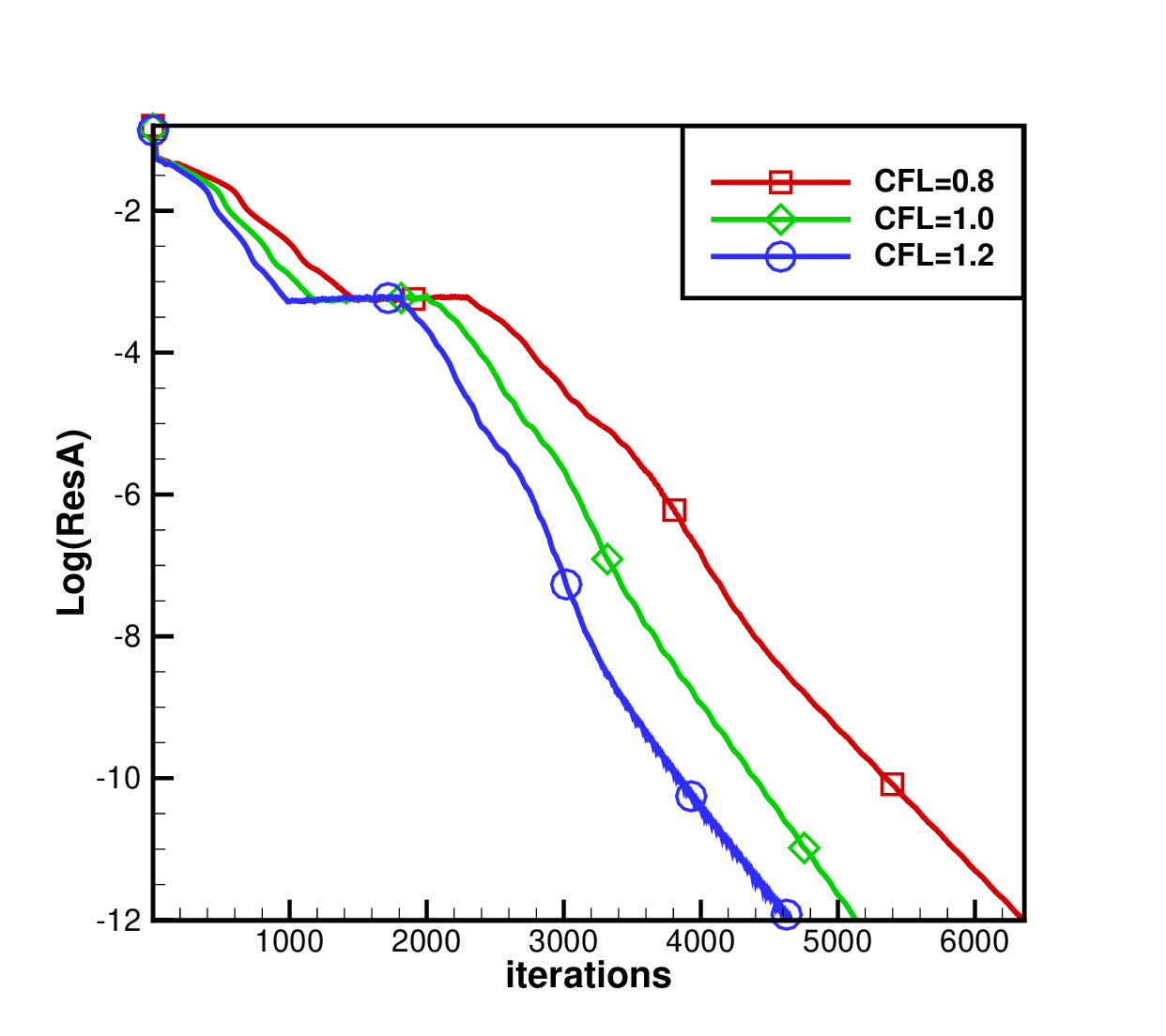}
\end{minipage}}
\subfigure[FS-WENOJS-AC]
{\begin{minipage}[t]{0.23\linewidth}
\includegraphics[width=1.7in]{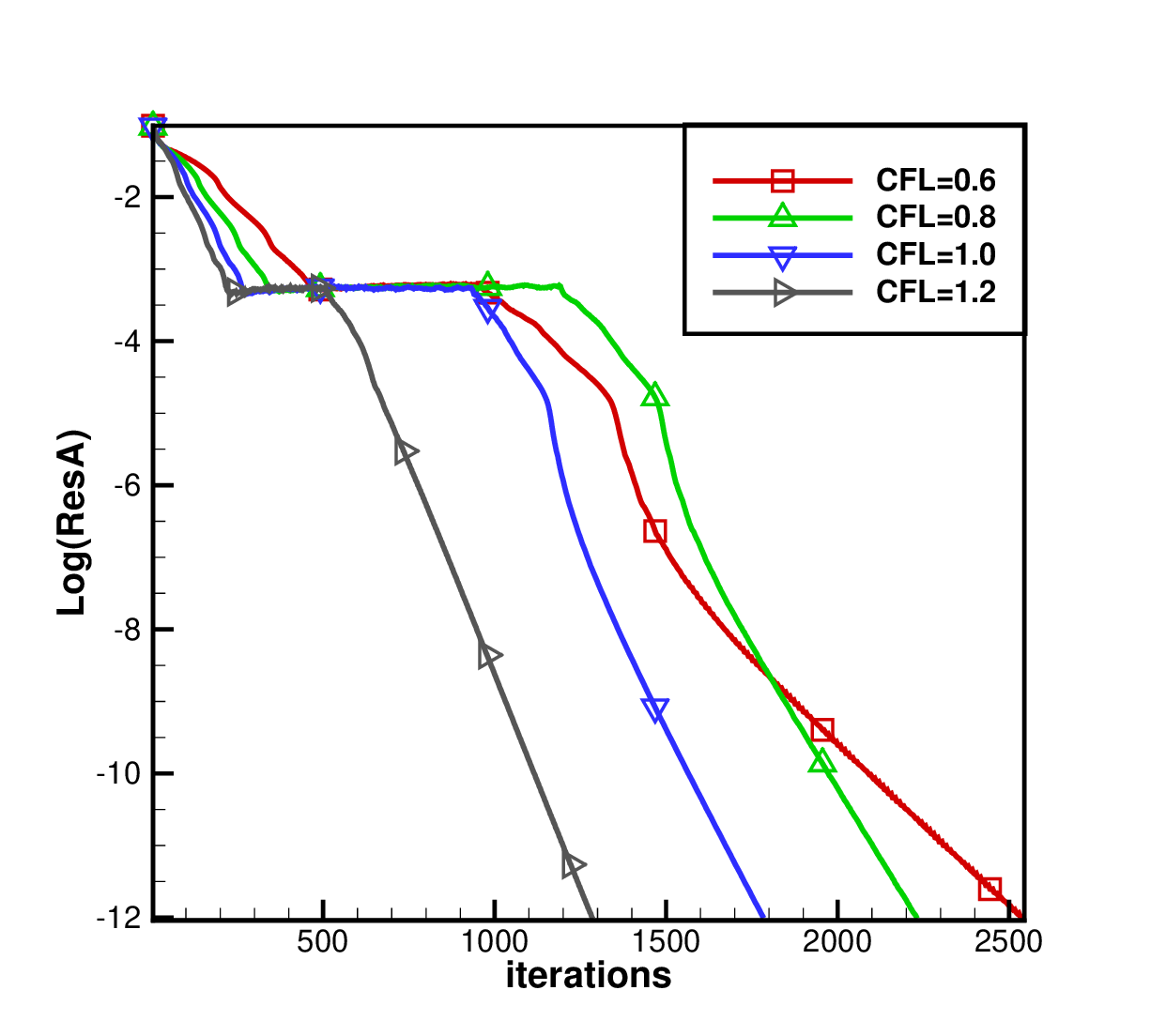}
\end{minipage}}
\caption{\label{4.3}Example 5: The convergence history of the residue as a function of number of iterations for four schemes.}
\end{figure}

\noindent This example has the same flow field conditions as the previous two, but the plates within the flow field are more numerous and longer, with their lengths extending all the way to the right boundary. The computational domain is $[0,5]\times[-5,5]$. The positions of the three long plates are: $x \in [2,5], y=2; x \in [2,5], y=0; x \in [2,5], y=-2$. A slip boundary condition is applied along the plates, while appropriate physical values are enforced at inflow and outflow boundaries. The domain is discretized using a $100 \times 200$ grids. In this case, the residue of the FS-WENOJS scheme still fails to converge to machine zero. The convergence history of the corresponding residue as a function of number of iterations for those schemes are presented in Figure \ref{4.3}. Table \ref{4.1} summarizes the computational performance of three WENO schemes under various CFL numbers. The FS-WENOJS-AC scheme demonstrates the fastest convergence, with CPU times lower than those of FS-MRWENO and RK-WENOJS-AC scheme. The CPU time of the FS-WENOJS-AC scheme is 80 percent of that of the FS-MRWENO scheme and 45 percent of that of the RK-WENOJS-AC scheme. Figure \ref{4.2} displays the pressure contours of the numerical solutions by three schemes, which show excellent agreement in structure.

\begin{table}
		\centering
\begin{tabular}{|c|c|c|}\hline
			\multicolumn{3}{|c|}{FS-MRWENO}\\\hline
            $\gamma:$ CFL number & iteration number  & CPU time \\\hline
0.6& 1016   &100.86 \\ \hline
0.8& 756   &75.83 \\ \hline
1.0& 600   &61.16 \\ \hline
			\multicolumn{3}{|c|}{RK-WENOJS-AC}\\\hline
            $\gamma:$ CFL number & iteration number  & CPU time \\\hline
0.8   &6360   & 148.84\\\hline
1.0   &5136   & 119.09\\\hline
1.2   &4653       &108.31\\\hline
			\multicolumn{3}{|c|}{FS-WENOJS-AC}\\\hline
            $\gamma:$ CFL number & iteration number  & CPU time \\\hline
0.6  &  2544     & 99.18\\\hline
0.8  &  2232     & 89.08\\\hline
1.0  & 1784     &76.69\\\hline
1.2  &  1286    & 48.72\\\hline
		\end{tabular}
		\caption{\label{4.1}Example 5: Number of iterations, the final time, and the total CPU time of three different iterative schemes when convergence is obtained. Convergence
criterion threshold value is $10^{-12}$. CPU time unit: second.}
	\end{table}

\bigskip
\noindent{\bf Example 6. $30^\circ$ steady oblique shock problem}

\begin{figure}
		\centering
\subfigure[FS-MRWENO]
{\begin{minipage}[t]{0.3\linewidth}
\includegraphics[width=2.0in]{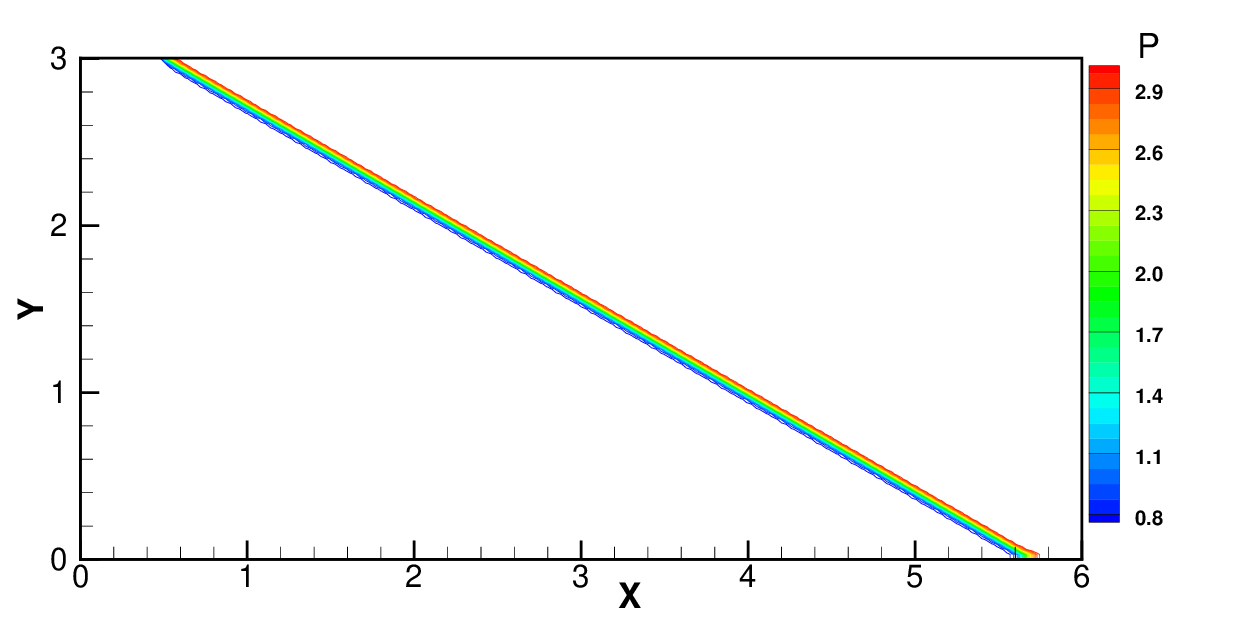}
\end{minipage}}
\subfigure[RK-WENOJS-AC]
{\begin{minipage}[t]{0.3\linewidth}
\includegraphics[width=2.0in]{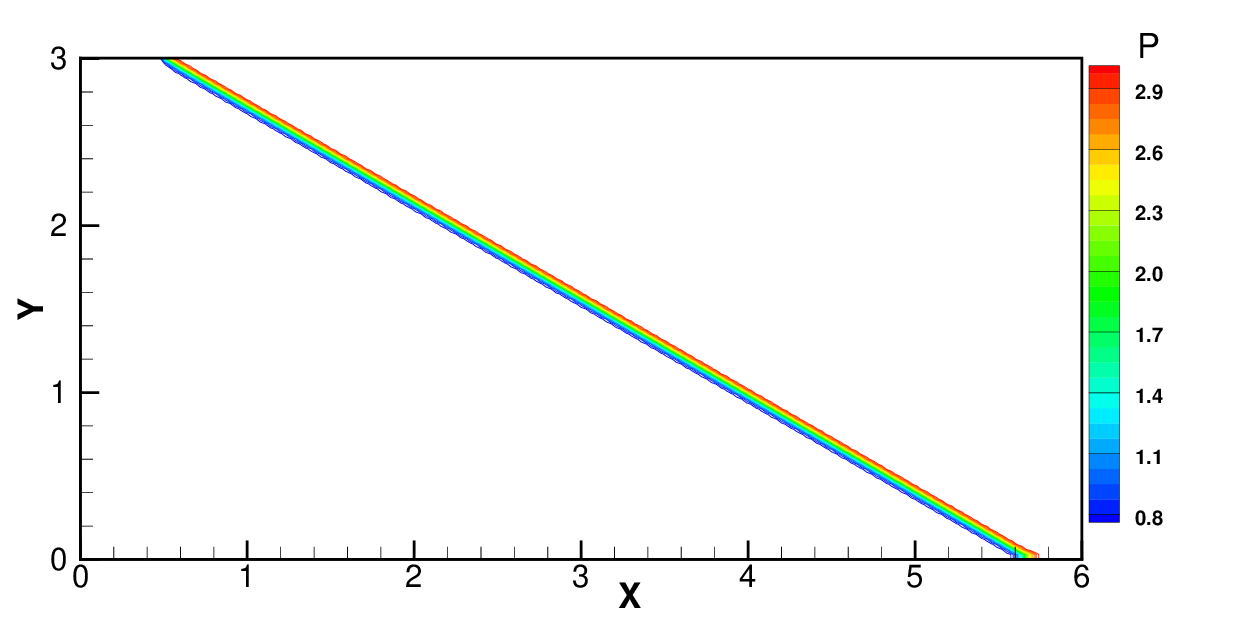}
\end{minipage}}
\subfigure[FS-WENOJS-AC]
{\begin{minipage}[t]{0.3\linewidth}
\includegraphics[width=2.0in]{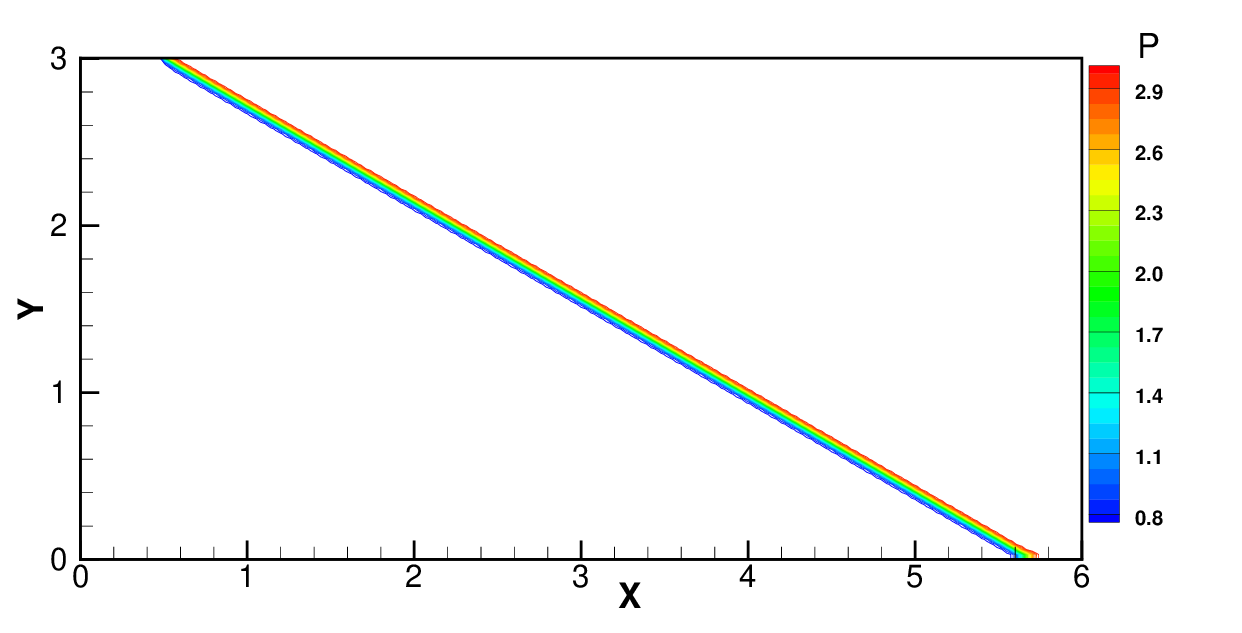}
\end{minipage}}
\caption{\label{5.2}Example 6: Thirty equally spaced pressure contours from 0.8 to 3.0 of the converged steady states of numerical solutions by three different iterative schemes.}
\end{figure}

\begin{figure}
		\centering
\subfigure[FS-MR-WENO]
{\begin{minipage}[t]{0.23\linewidth}
\includegraphics[width=1.7in]{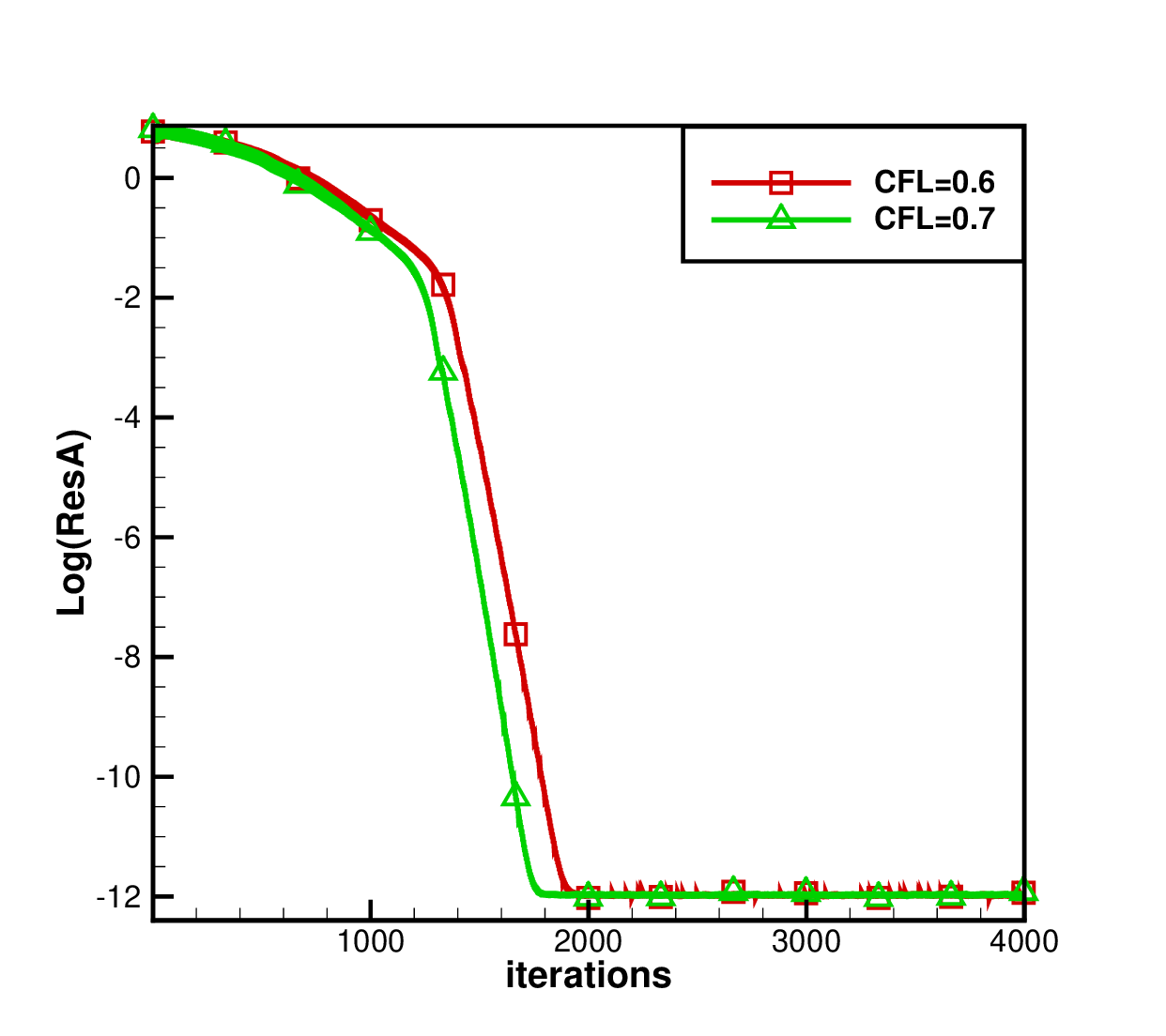}
\end{minipage}}
\subfigure[FS-WENOJS]
{\begin{minipage}[t]{0.23\linewidth}
\includegraphics[width=1.7in]{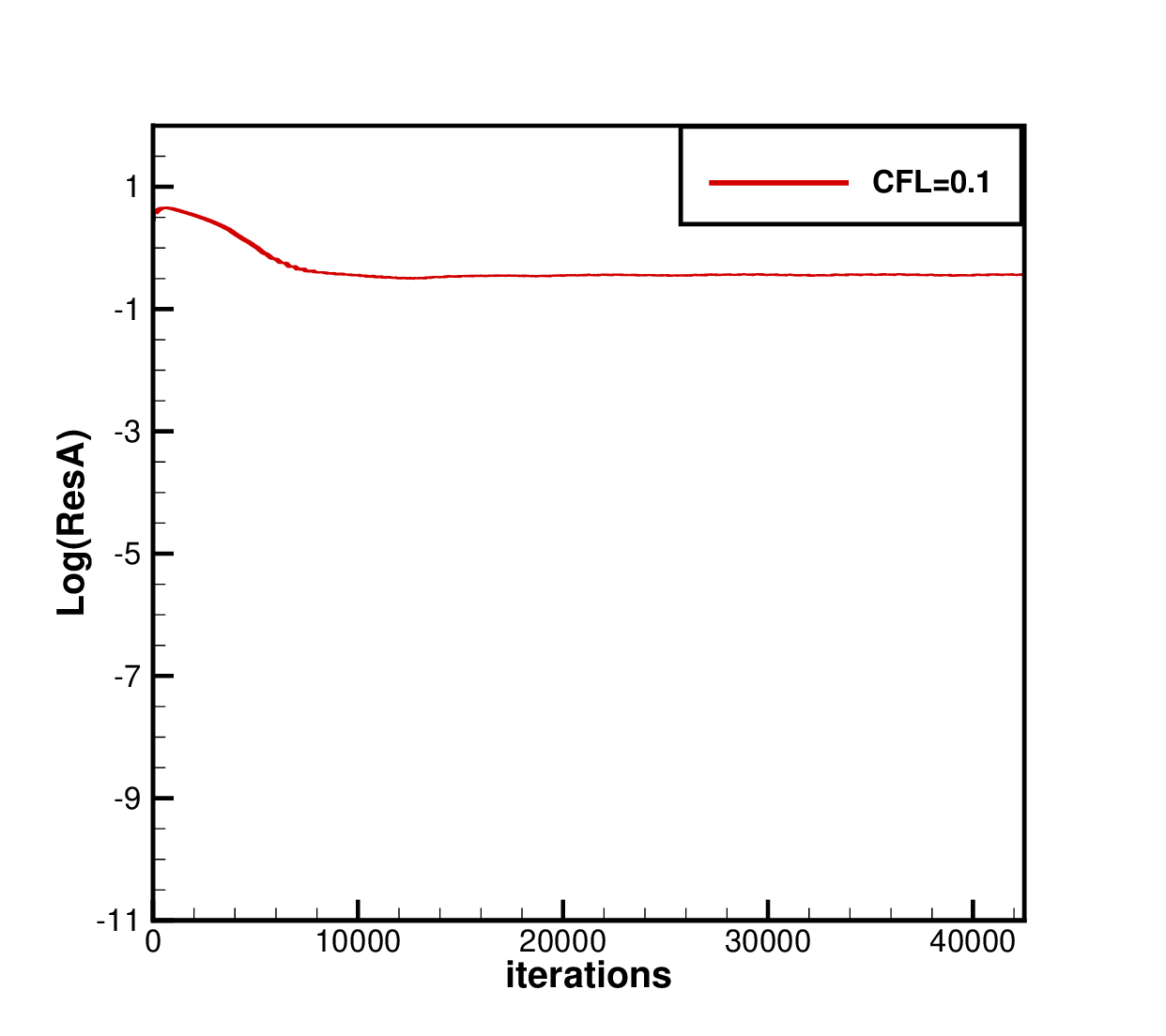}
\end{minipage}}
\subfigure[RK-WENOJS-AC]
{\begin{minipage}[t]{0.23\linewidth}
\includegraphics[width=1.7in]{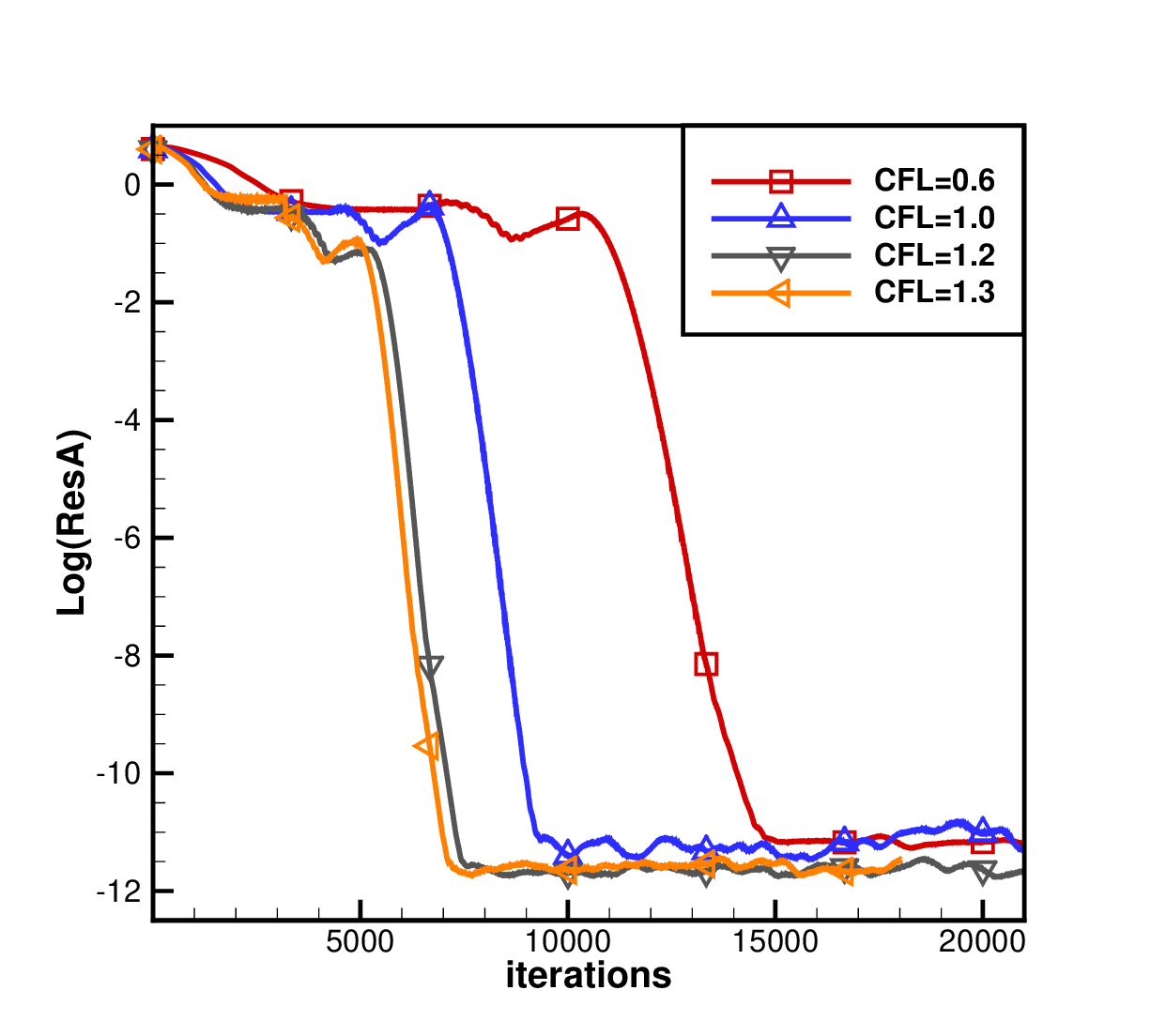}
\end{minipage}}
\subfigure[FS-WENOJS-AC]
{\begin{minipage}[t]{0.23\linewidth}
\includegraphics[width=1.7in]{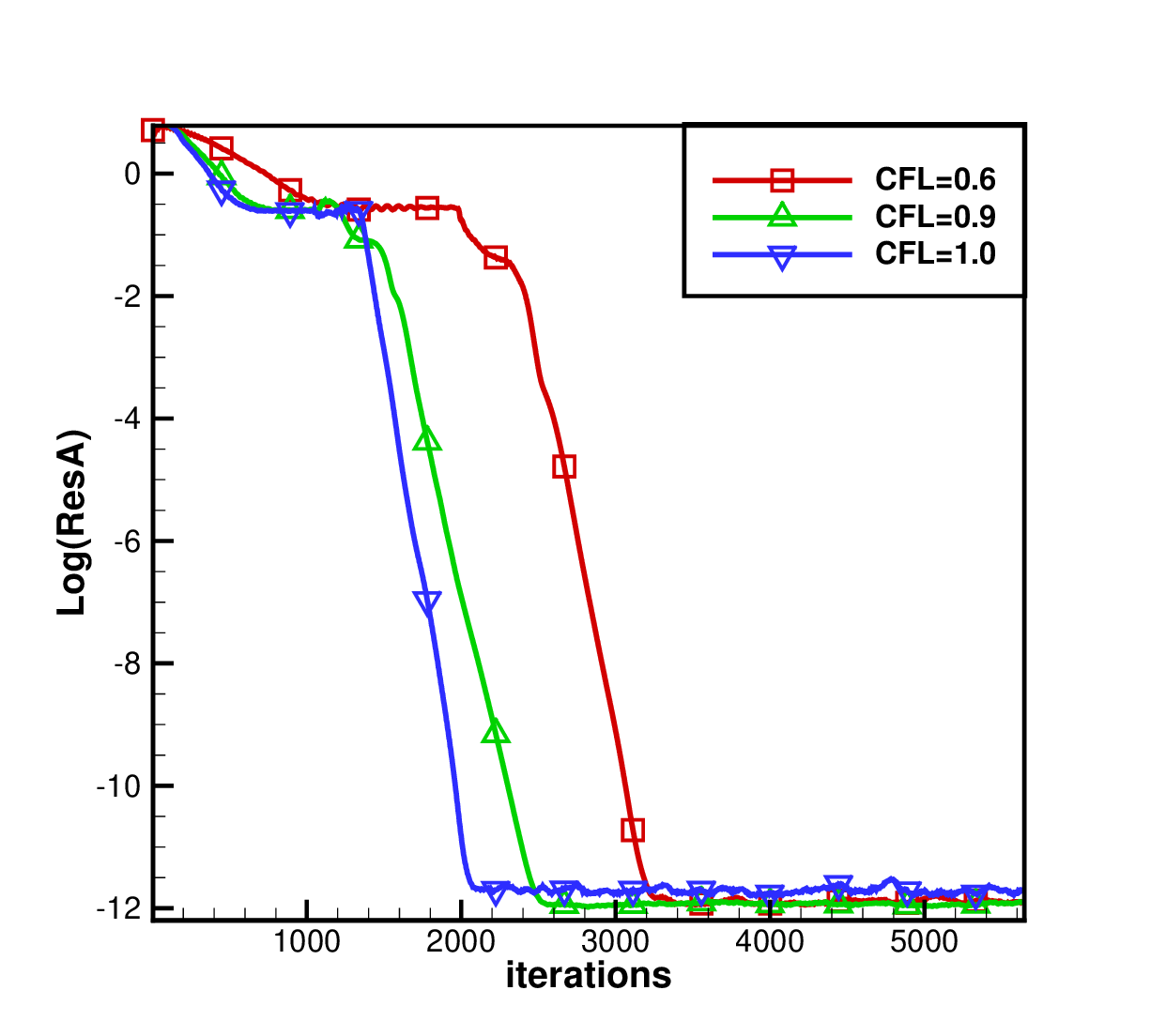}
\end{minipage}}
\caption{\label{5.3}Example 6: The convergence history of the residue as a function of number of iterations for three schemes.}
\end{figure}

\noindent This example investigates a two-dimensional steady oblique shock problem \cite{liangfu}, where a stationary shock wave forms at a $30^\circ$ angle relative to the $x$-axis. The computational domain covers $(x,y)\in[0,6]\times[0,3]$, with outflow conditions imposed on the right and bottom boundaries. The left boundary is assigned $(\rho,u,v,p)=(1,4,0,1/1.4)$. A uniform grid of 200 $\times$ 100 is used for the computation. On the top boundary, a piecewise condition is applied: for $x$ in $[0,0.5]$, the state is $(\rho,u,v,p)=(1,4,0,1/1.4)$; for $x$ in $(0.5,6]$, it becomes $(2.6667,3.3750,-1.0825,3.2143)$. The initial flow field is set equal to the values prescribed at the left boundary. In this numerical example, the residual of the scheme can only decrease to between $10^{-11}$ and $10^{-12}$, so we set the convergence criterion to $10^{-11}$. Table \ref{5.1} shows the computational efficiency of three schemes. We observe that the RK-WENOJS-AC scheme can still achieve a larger CFL number, but the FS-WENOJS-AC scheme remains the most computationally efficient, saving nearly half of the CPU time. Figure \ref{5.3} presents the evolution of residuals over time for four schemes. The FS-WENOJS scheme still fails to achieve full convergence, while the other three schemes all absolutely converge. Figure \ref{5.2} also shows the pressure contour plots of the three schemes that achieve full convergence, and their numerical results are similar.

\begin{table}
		\centering
\begin{tabular}{|c|c|c|}\hline
			\multicolumn{3}{|c|}{FS-MRWENO}\\\hline
            $\gamma:$ CFL number & iteration number  & CPU time \\\hline
0.6	&1835		&176.13\\\hline
0.7	&1698		&160.91\\\hline
			\multicolumn{3}{|c|}{RK-WENOJS-AC}\\\hline
            $\gamma:$ CFL number & iteration number & CPU time \\\hline
0.6	&14673		&321.30\\\hline
0.8	&14247		&302.97\\\hline
1.0	&9222		&200.56\\\hline
1.2	&7296		&157.05\\\hline
1.3	&6972		&151.53\\\hline
			\multicolumn{3}{|c|}{FS-WENOJS-AC}\\\hline
            $\gamma:$ CFL number & iteration number  & CPU time \\\hline
0.6	&3136		&122.08\\\hline
0.9	&2396		&87.91\\\hline
1.0	&2012		&76.67\\\hline
		\end{tabular}
		\caption{\label{5.1}Example 6: Number of iterations and the total CPU time of three different iterative schemes when convergence is obtained. Convergence criterion threshold value is $10^{-11}$. CPU time unit: second.}
	\end{table}

\bigskip
\noindent{\bf Example 7. Supersonic flows past an airfoil}

\begin{figure}
		\centering
\subfigure[FS-MRWENO]
{\begin{minipage}[t]{0.3\linewidth}
\includegraphics[width=2.0in]{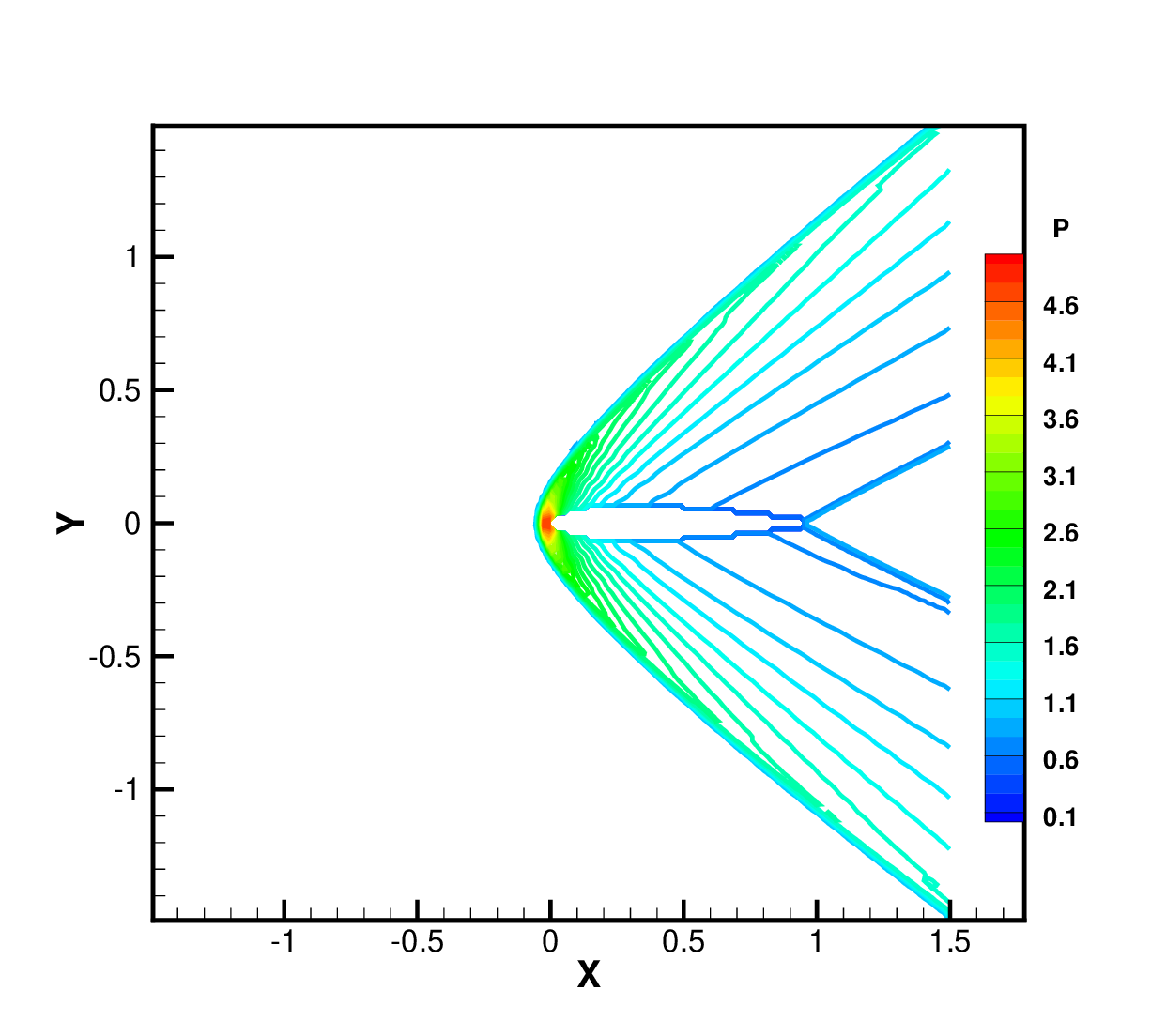}
\end{minipage}}
\subfigure[RK-WENOJS-AC]
{\begin{minipage}[t]{0.3\linewidth}
\includegraphics[width=2.0in]{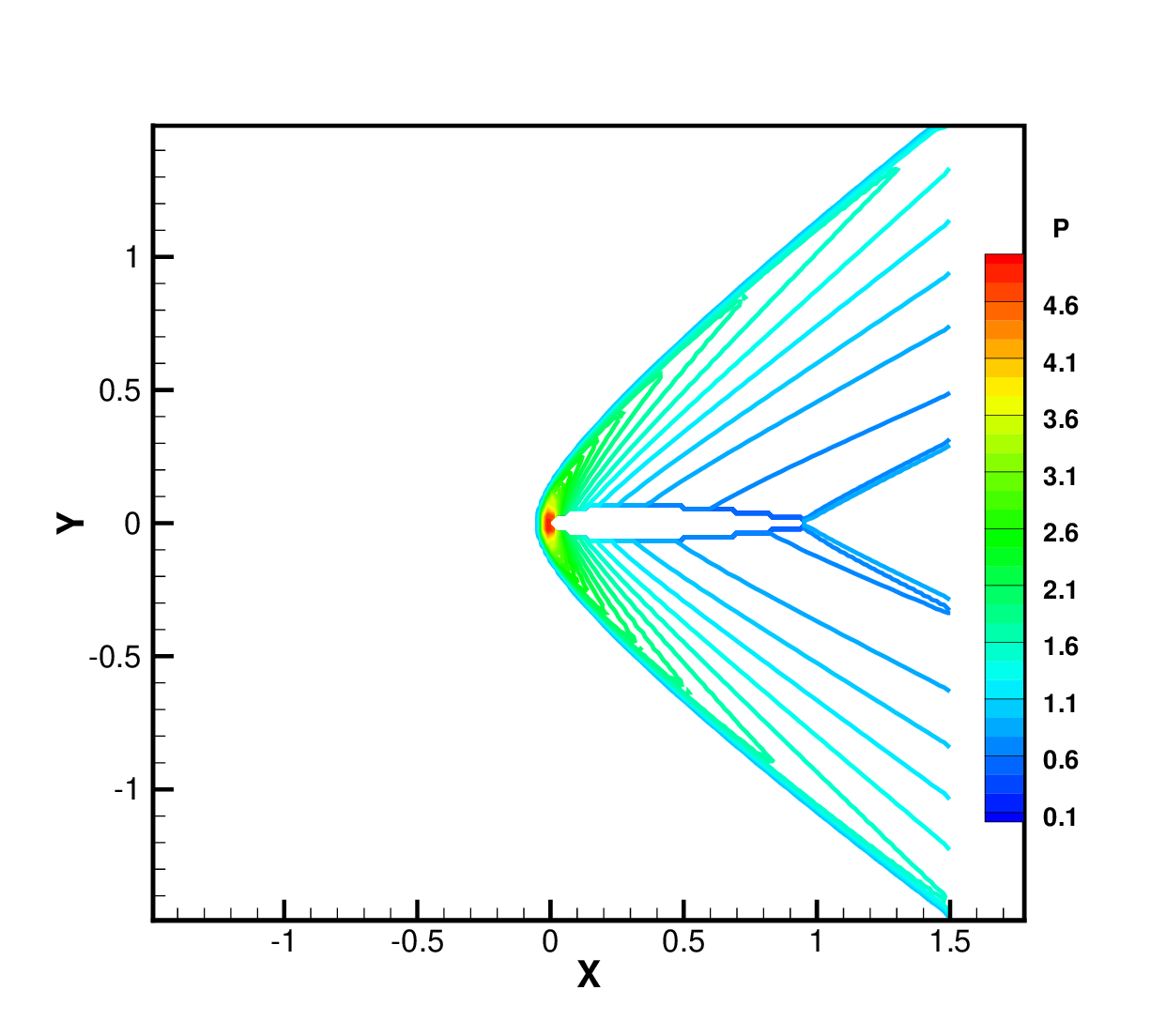}
\end{minipage}}
\subfigure[FS-WENOJS-AC]
{\begin{minipage}[t]{0.3\linewidth}
\includegraphics[width=2.0in]{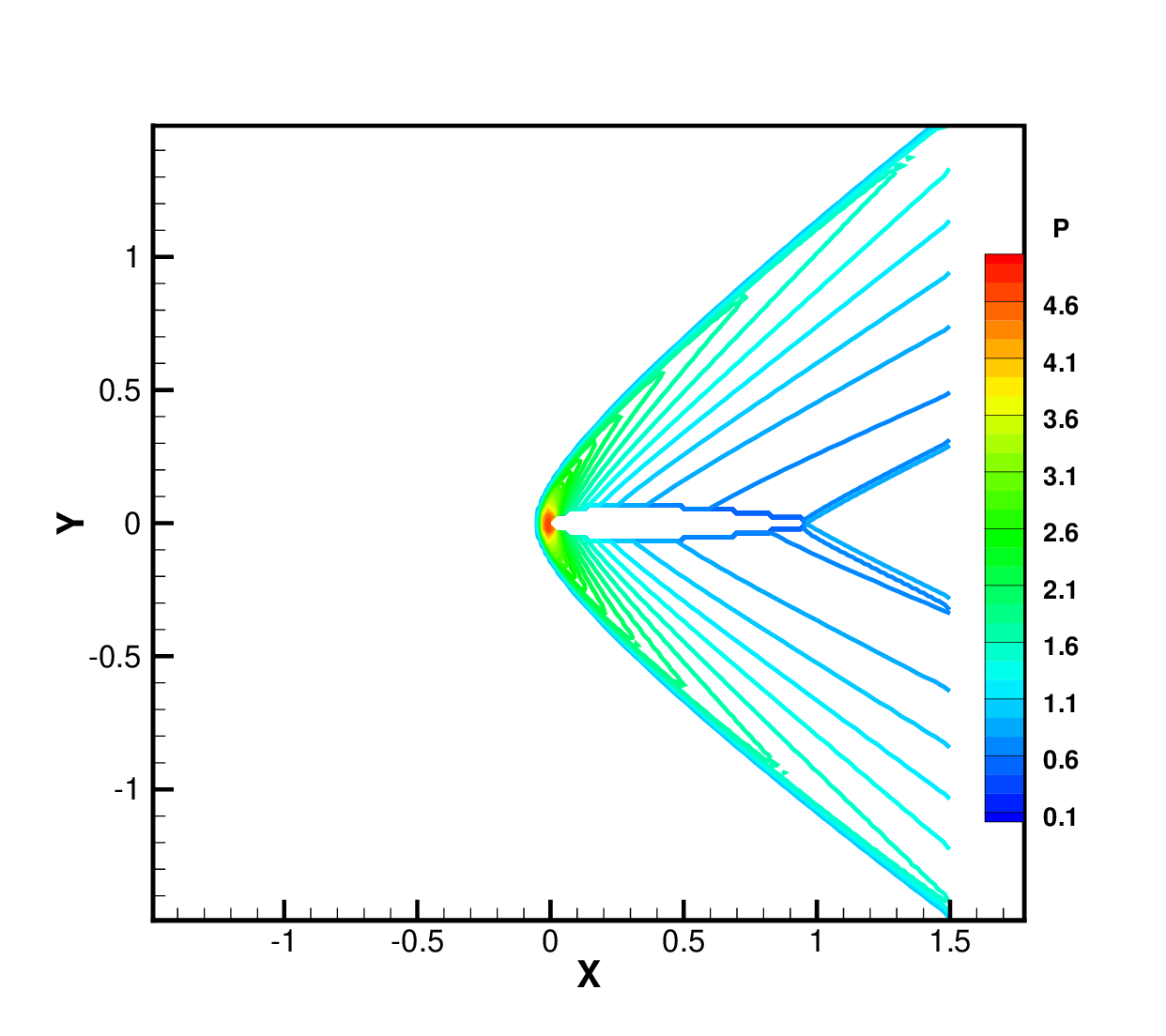}
\end{minipage}}
\caption{\label{6.2}Example 7: Thirty equally spaced pressure contours from 0.1 to 4.9 of the converged steady states of numerical solutions by three different iterative schemes.}
\end{figure}

\begin{figure}
		\centering
\subfigure[FS-MRWENO]
{\begin{minipage}[t]{0.23\linewidth}
\includegraphics[width=1.7in]{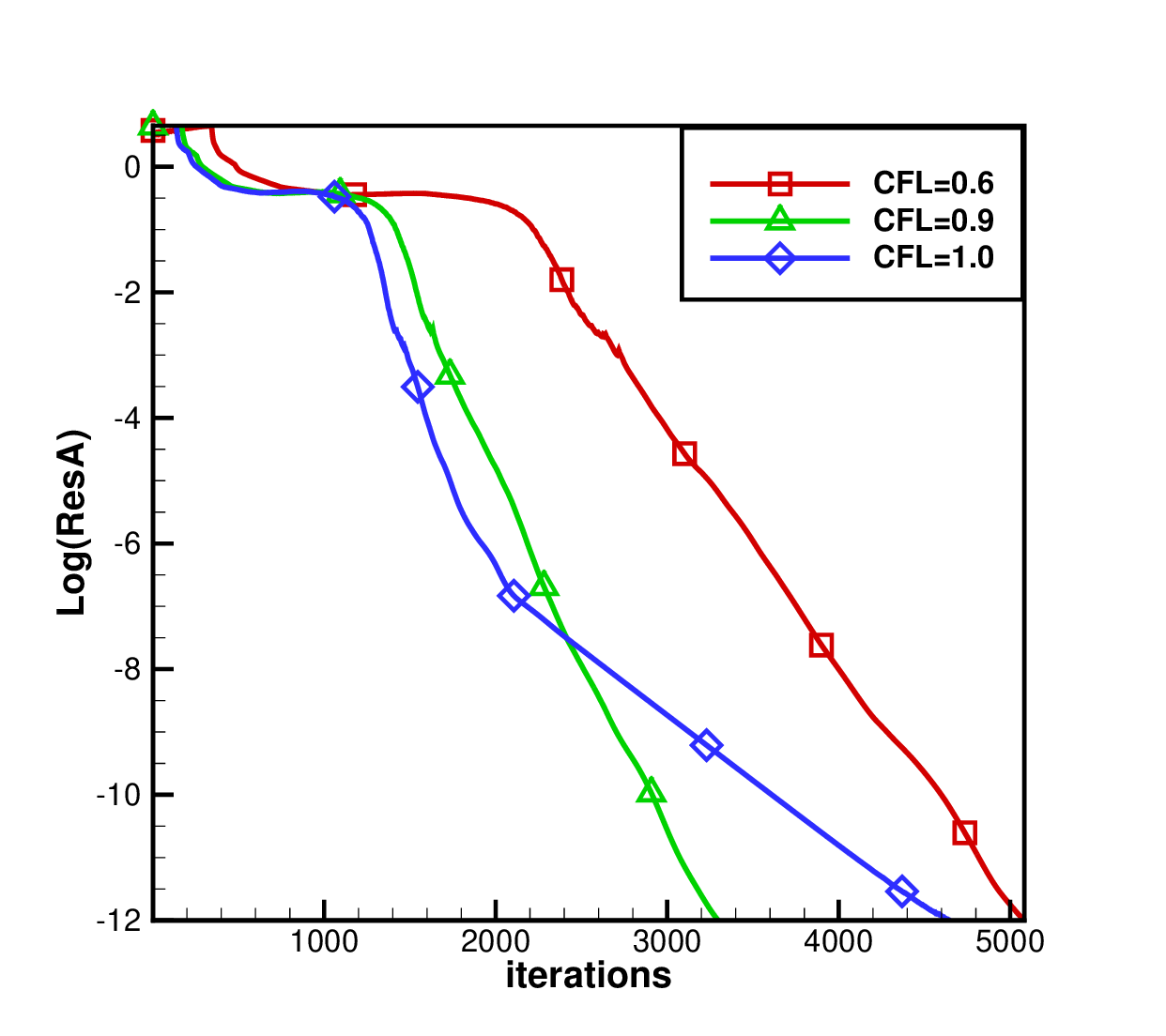}
\end{minipage}}
\subfigure[FS-WENOJS]
{\begin{minipage}[t]{0.23\linewidth}
\includegraphics[width=1.7in]{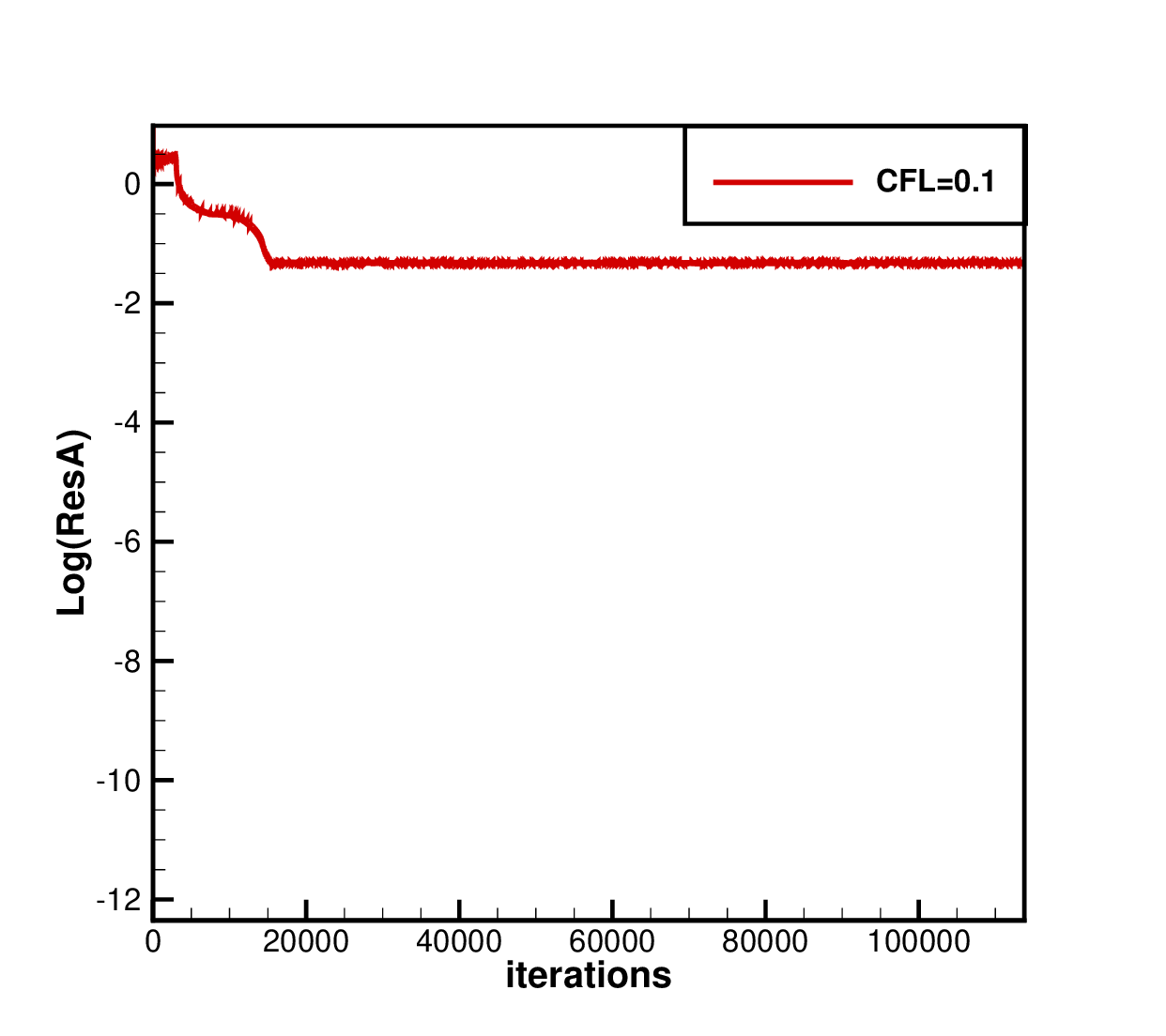}
\end{minipage}}
\subfigure[RK-WENOJS-AC]
{\begin{minipage}[t]{0.23\linewidth}
\includegraphics[width=1.7in]{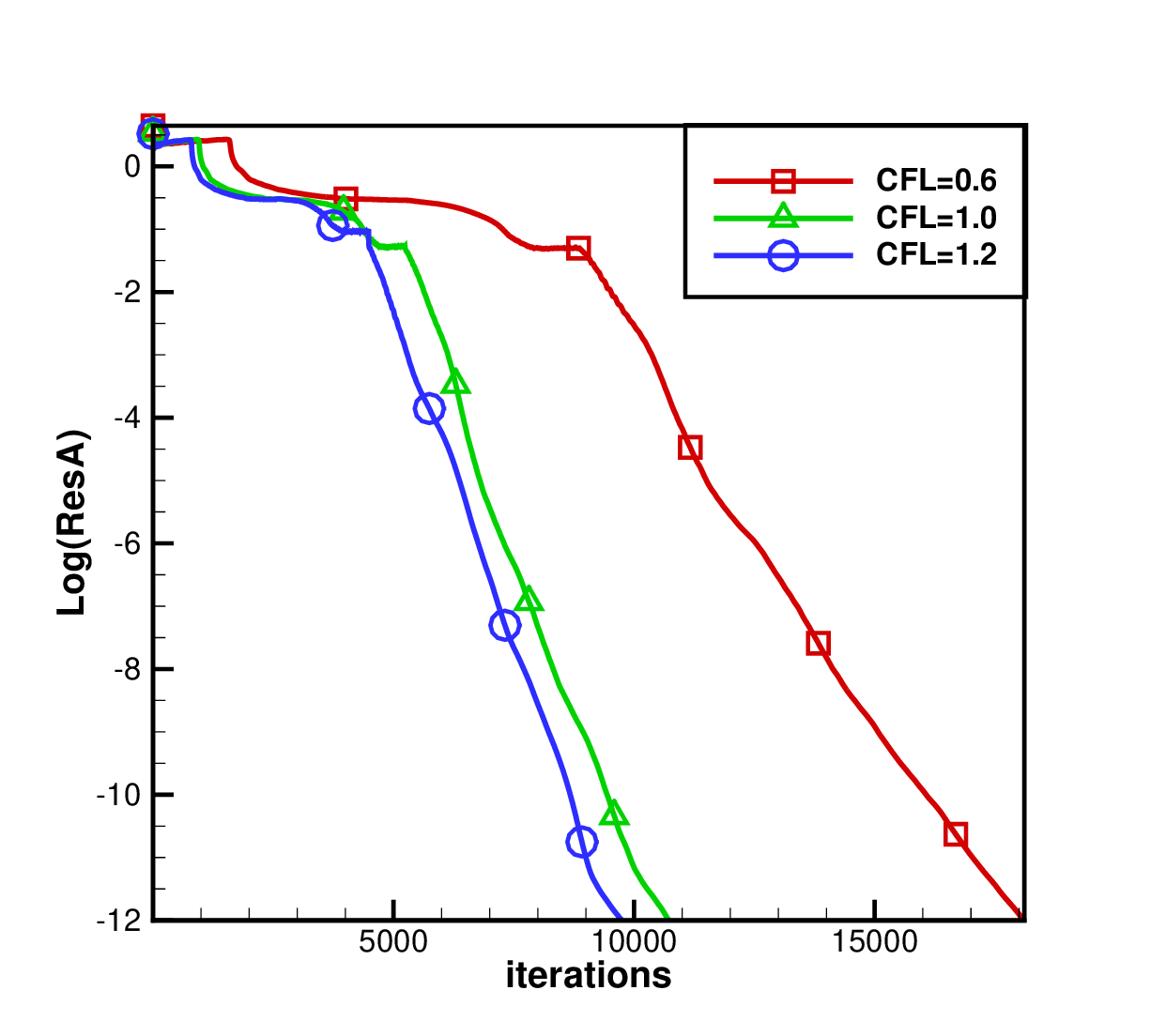}
\end{minipage}}
\subfigure[FS-WENOJS-AC]
{\begin{minipage}[t]{0.23\linewidth}
\includegraphics[width=1.7in]{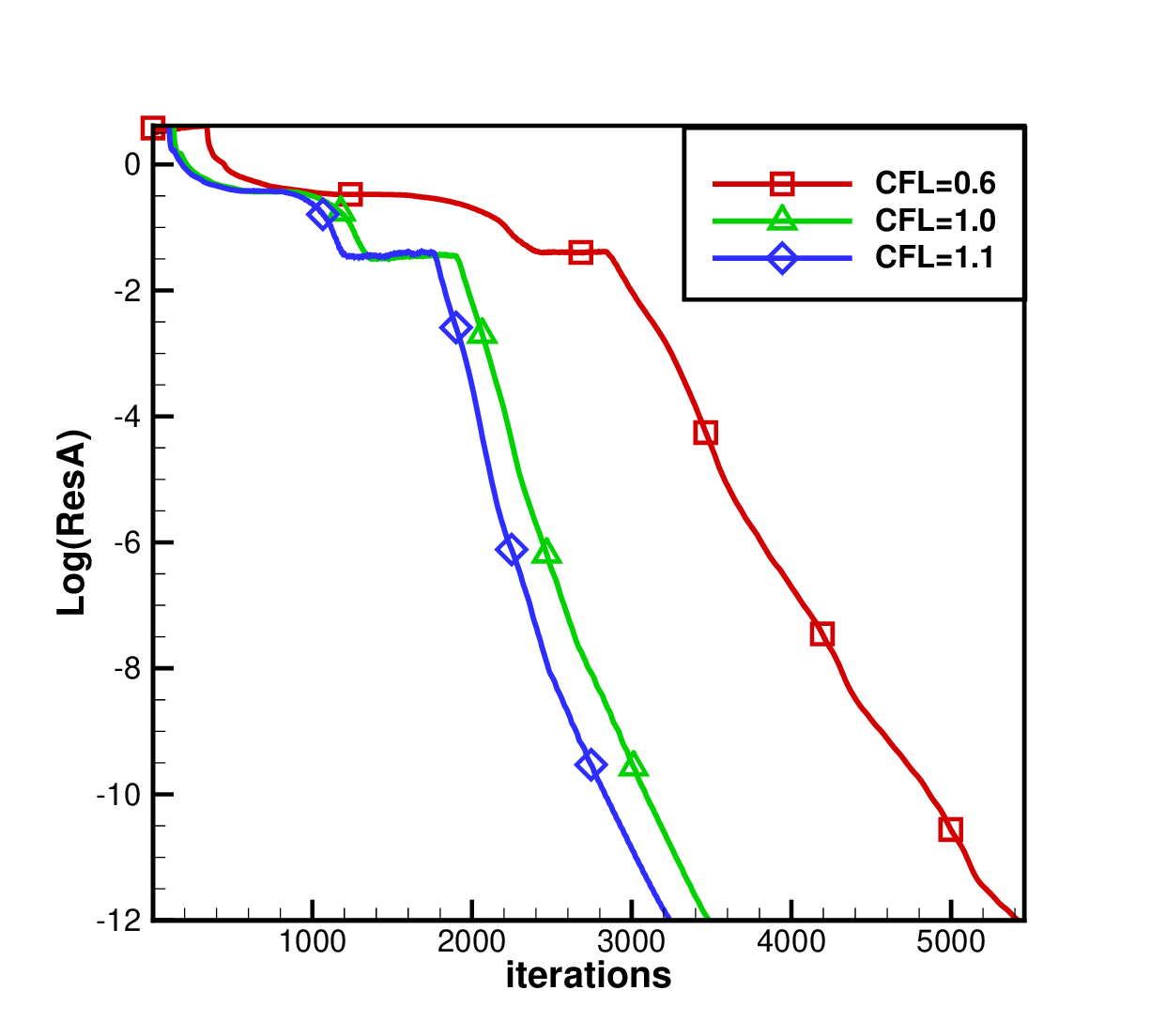}
\end{minipage}}
\caption{\label{6.3}Example 7: The convergence history of the residue as a function of number of iterations for three schemes.}
\end{figure}

\noindent This example addresses supersonic flow over a NACA0012 airfoil, as studied in \cite{NACA,NACA2}. Adopting the configuration from \cite{m2}, the flow conditions are set to Mach number $Ma=2$ and an angle of attack of $\alpha=1^{\circ}$. The simulation domain spans $[-1.5, 1.5] \times [-1.5, 1.5]$, discretized using a $200\times{200}$ grid. The boundary conditions of the airfoil are handled with the fifth-order ILW numerical boundary method \cite{ILW3,S.Tan,SIRUIT}. Table \ref{6.1} presents the number of iterations and total CPU time required to meet the convergence criterion of $10^{-12}$ for the FS-MRWENO, RK-WENOJS-AC, and FS-WENOJS-AC schemes at various CFL numbers. This test case once again verifies the high efficiency of the FS-WENOJS-AC scheme. Compared to the FS-MRWENO scheme, it reduces the fastest convergence time by $62\%$, and by $59\%$ compared to the RK-WENOJS-AC scheme. Figure \ref{6.2} shows the pressure contour plots of the three WENO schemes, which are nearly identical. Figure \ref{6.3} presents the residue evolution plots of the three schemes, all of which achieve full convergence at relatively large CFL numbers.

\begin{table}
		\centering
\begin{tabular}{|c|c|c|}\hline
			\multicolumn{3}{|c|}{FS-MRWENO}\\\hline
            $\gamma:$ CFL number & iteration number  & CPU time \\\hline
0.6	&5082		&807.28\\\hline
0.9	&3298		&517.30\\\hline
1.0	&4646		&735.80\\\hline
			\multicolumn{3}{|c|}{RK-WENOJS-AC}\\\hline
            $\gamma:$ CFL number & iteration number & CPU time \\\hline
0.6	&18108		&879.02\\\hline
1.0	&10728		&525.05\\\hline
1.2	&9747		&476.06\\\hline
			\multicolumn{3}{|c|}{FS-WENOJS-AC}\\\hline
            $\gamma:$ CFL number & iteration number  & CPU time \\\hline
0.6	&5458		&351.09\\\hline
1.0	&3482		&213.70\\\hline
1.1	&3238		&198.05\\\hline
		\end{tabular}
		\caption{\label{6.1}Example 7: Number of iterations and the total CPU time of three different iterative schemes when convergence is obtained. Convergence
criterion threshold value is $10^{-12}$. CPU time unit: second.}
	\end{table}

\bigskip
\noindent{\bf Example 8. Subsonic flows past an airfoil}

\begin{figure}
		\centering
\subfigure[FS-MRWENO]
{\begin{minipage}[t]{0.3\linewidth}
\includegraphics[width=2.0in]{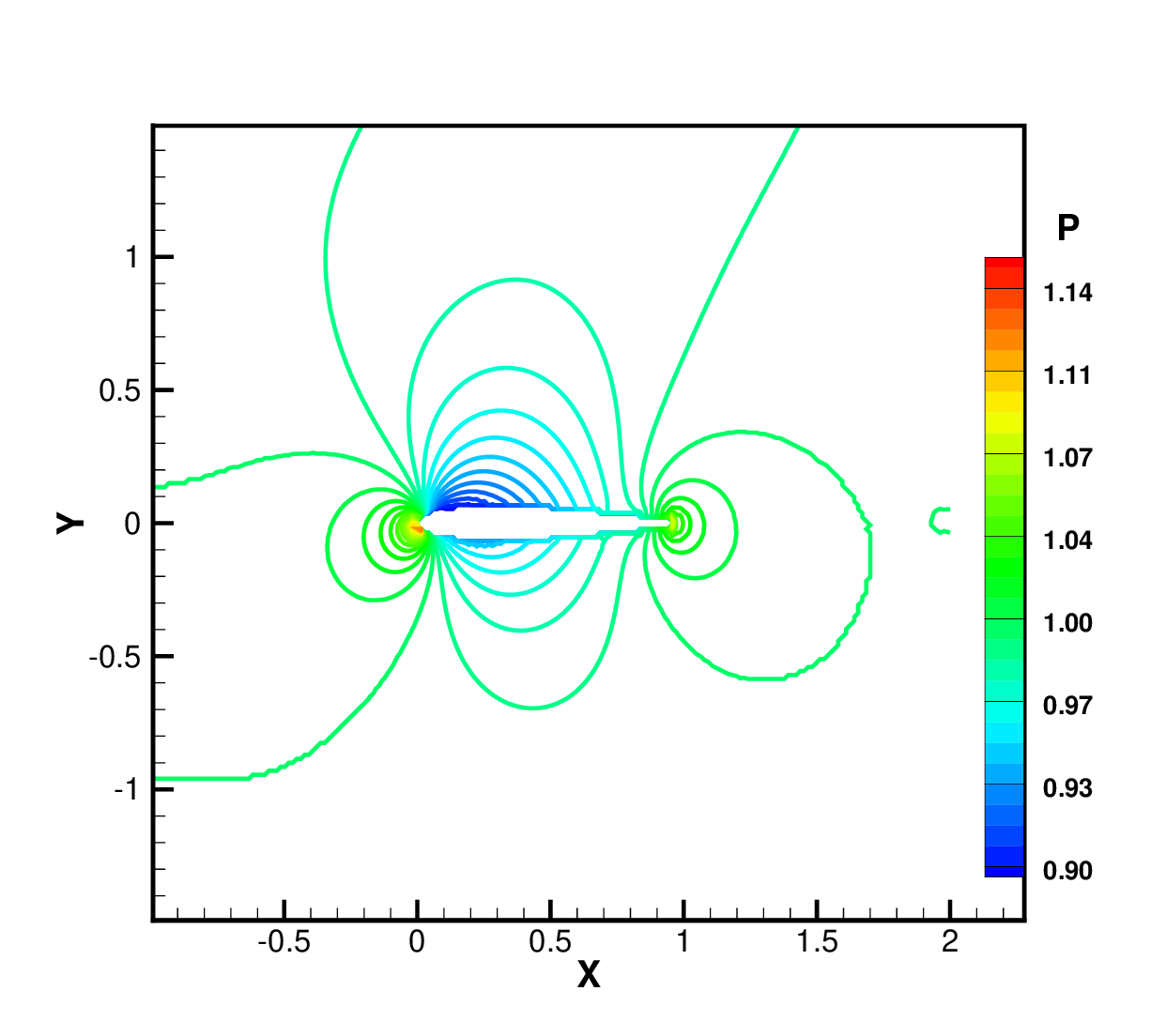}
\end{minipage}}
\subfigure[RK-WENOJS-AC]
{\begin{minipage}[t]{0.3\linewidth}
\includegraphics[width=2.0in]{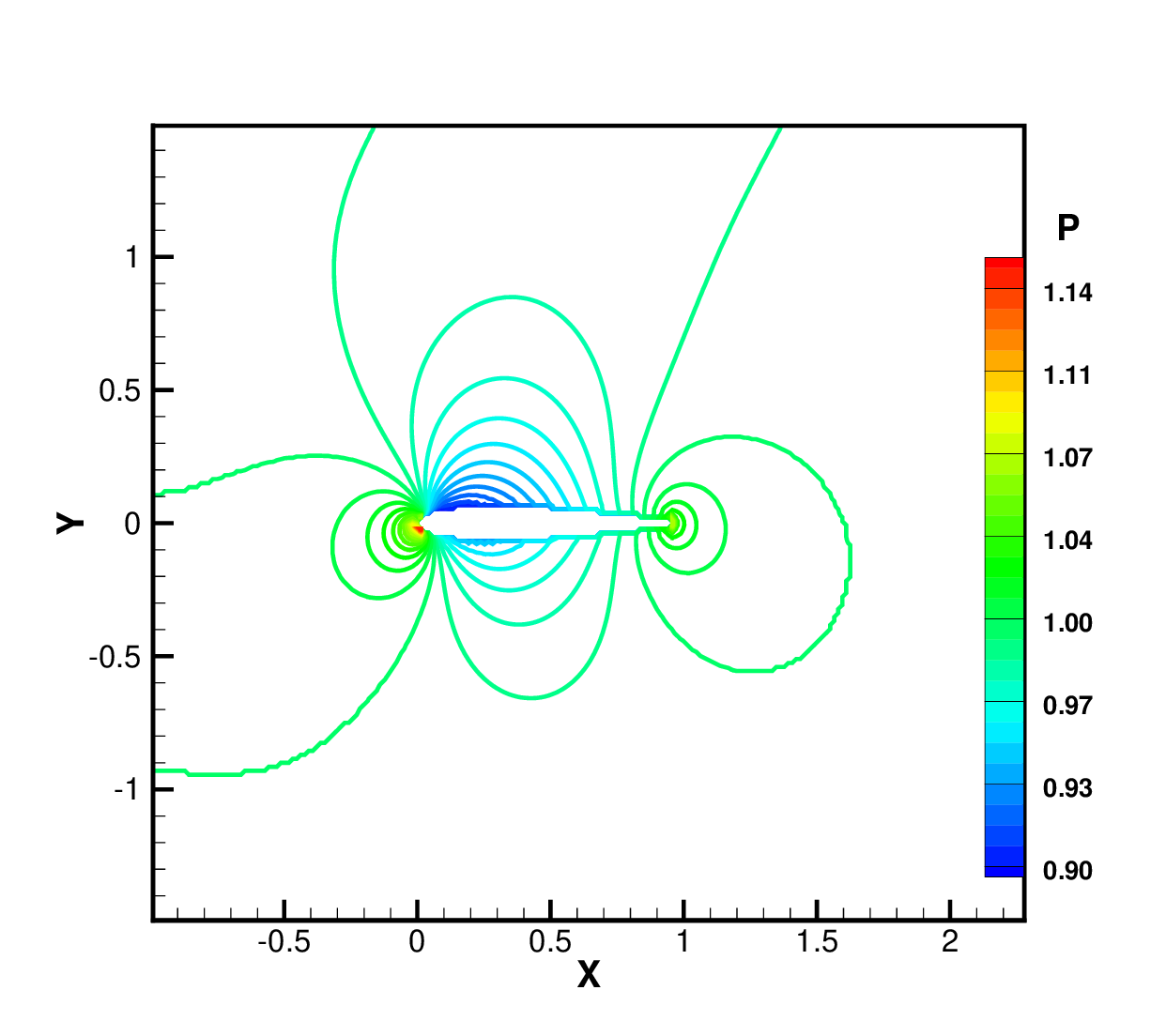}
\end{minipage}}
\subfigure[FS-WENOJS-AC]
{\begin{minipage}[t]{0.3\linewidth}
\includegraphics[width=2.0in]{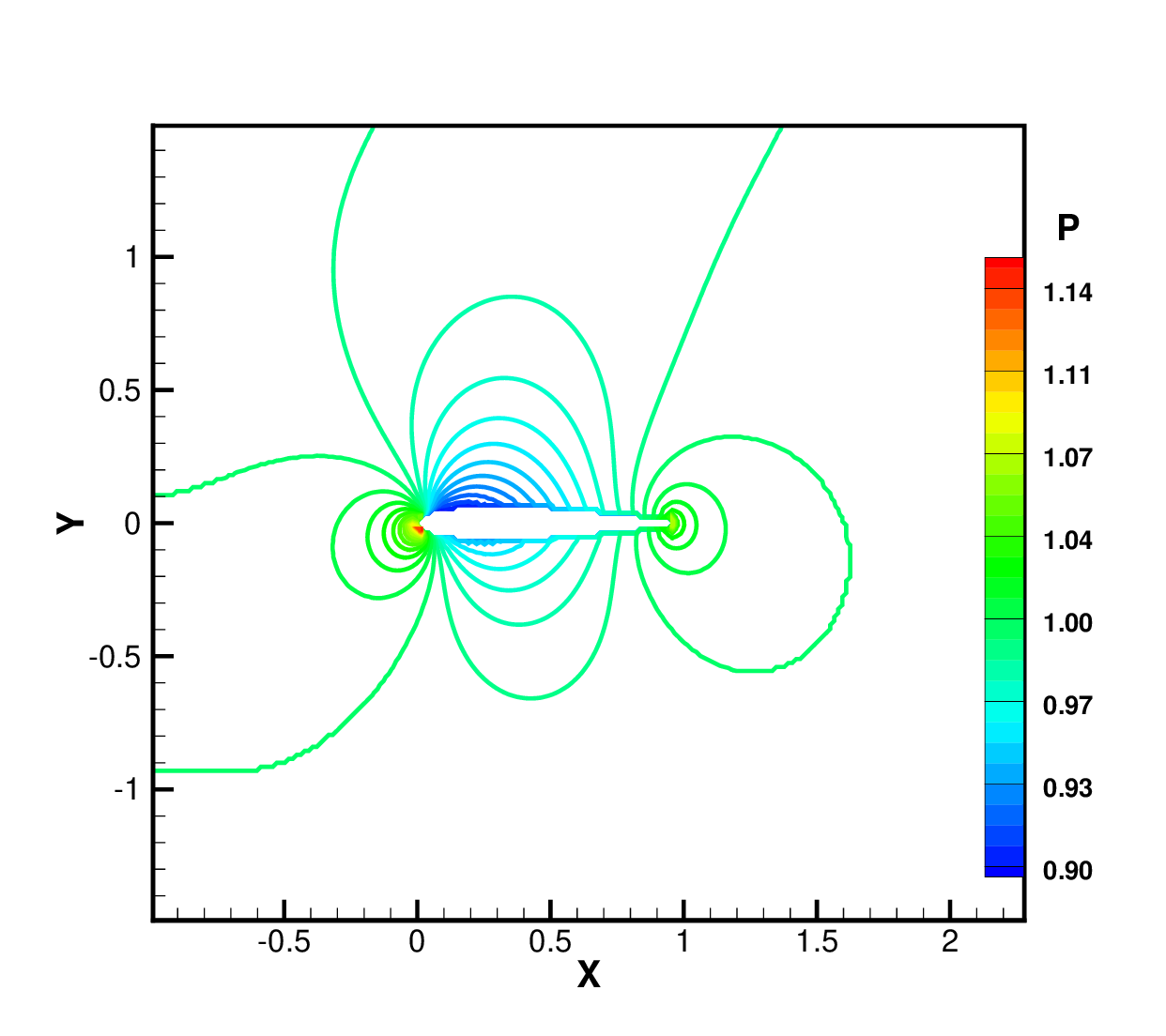}
\end{minipage}}
\caption{\label{8.2}Example 8: Thirty equally spaced pressure contours from 0.9 to 1.15 of the converged steady states of numerical solutions by three different iterative schemes.}
\end{figure}

\begin{figure}
		\centering
\subfigure[FS-MR-WENO]
{\begin{minipage}[t]{0.23\linewidth}
\includegraphics[width=1.7in]{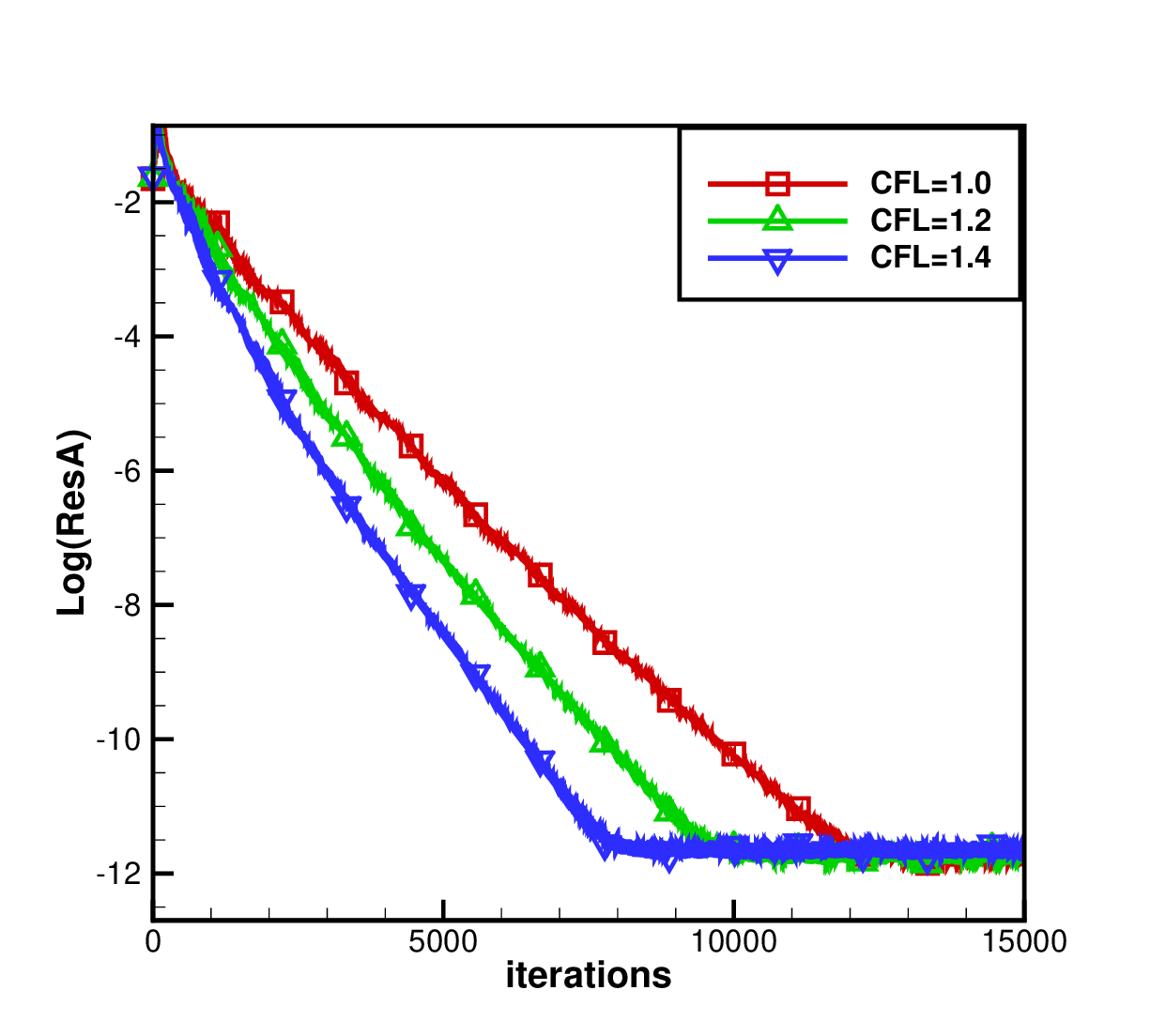}
\end{minipage}}
\subfigure[FS-WENOJS]
{\begin{minipage}[t]{0.23\linewidth}
\includegraphics[width=1.7in]{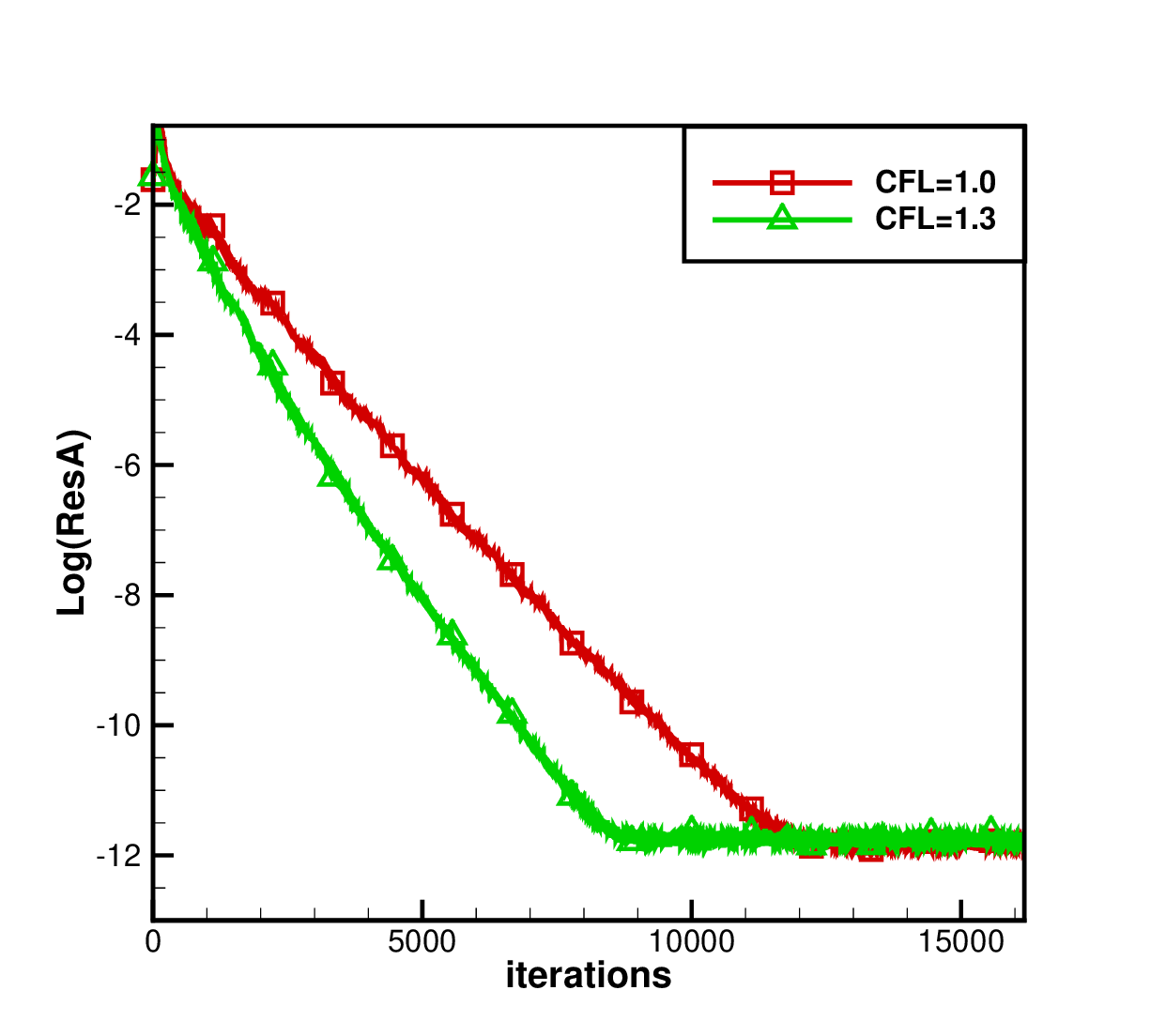}
\end{minipage}}
\subfigure[RK-WENOJS-AC]
{\begin{minipage}[t]{0.23\linewidth}
\includegraphics[width=1.7in]{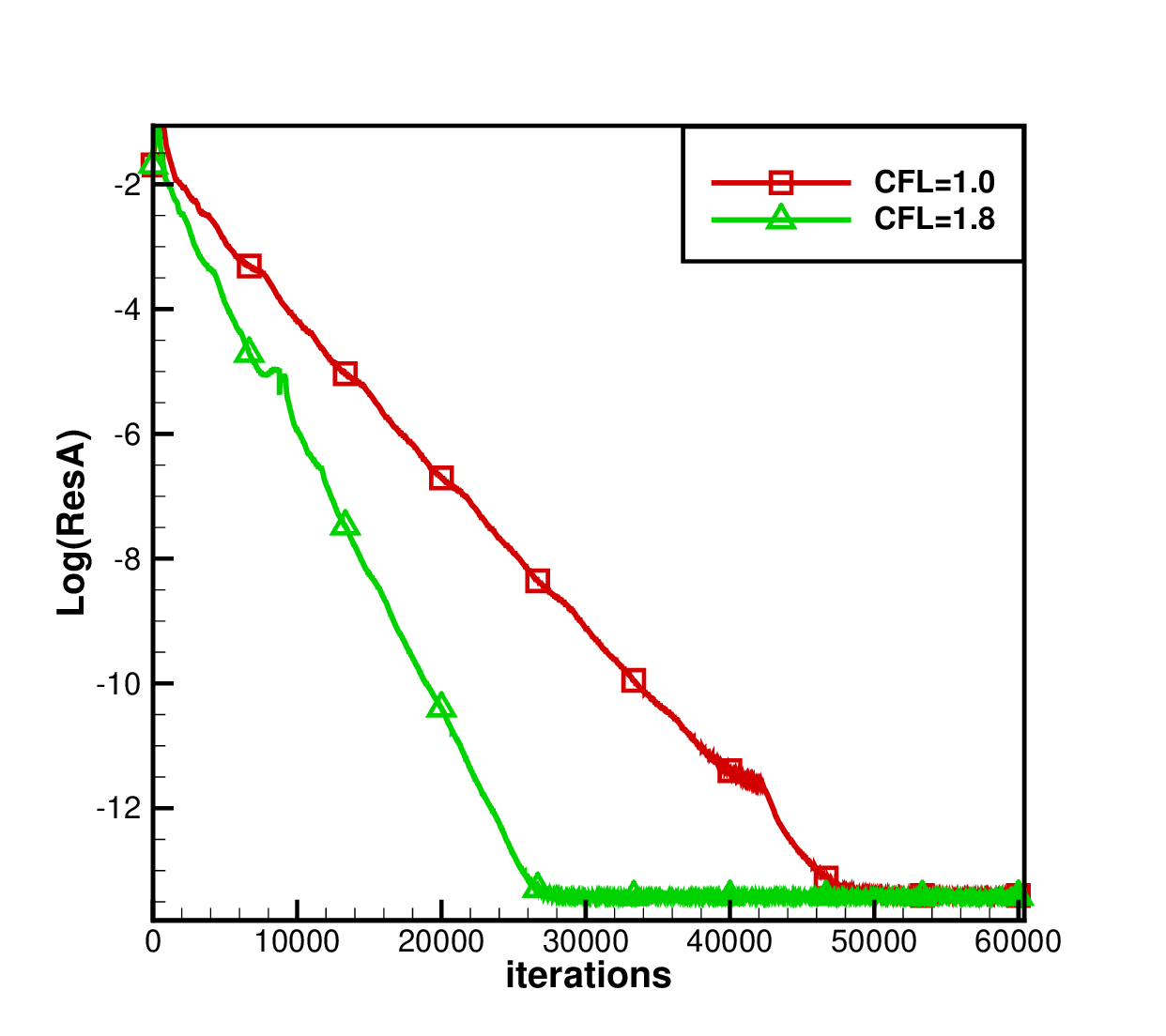}
\end{minipage}}
\subfigure[FS-WENOJS-AC]
{\begin{minipage}[t]{0.23\linewidth}
\includegraphics[width=1.7in]{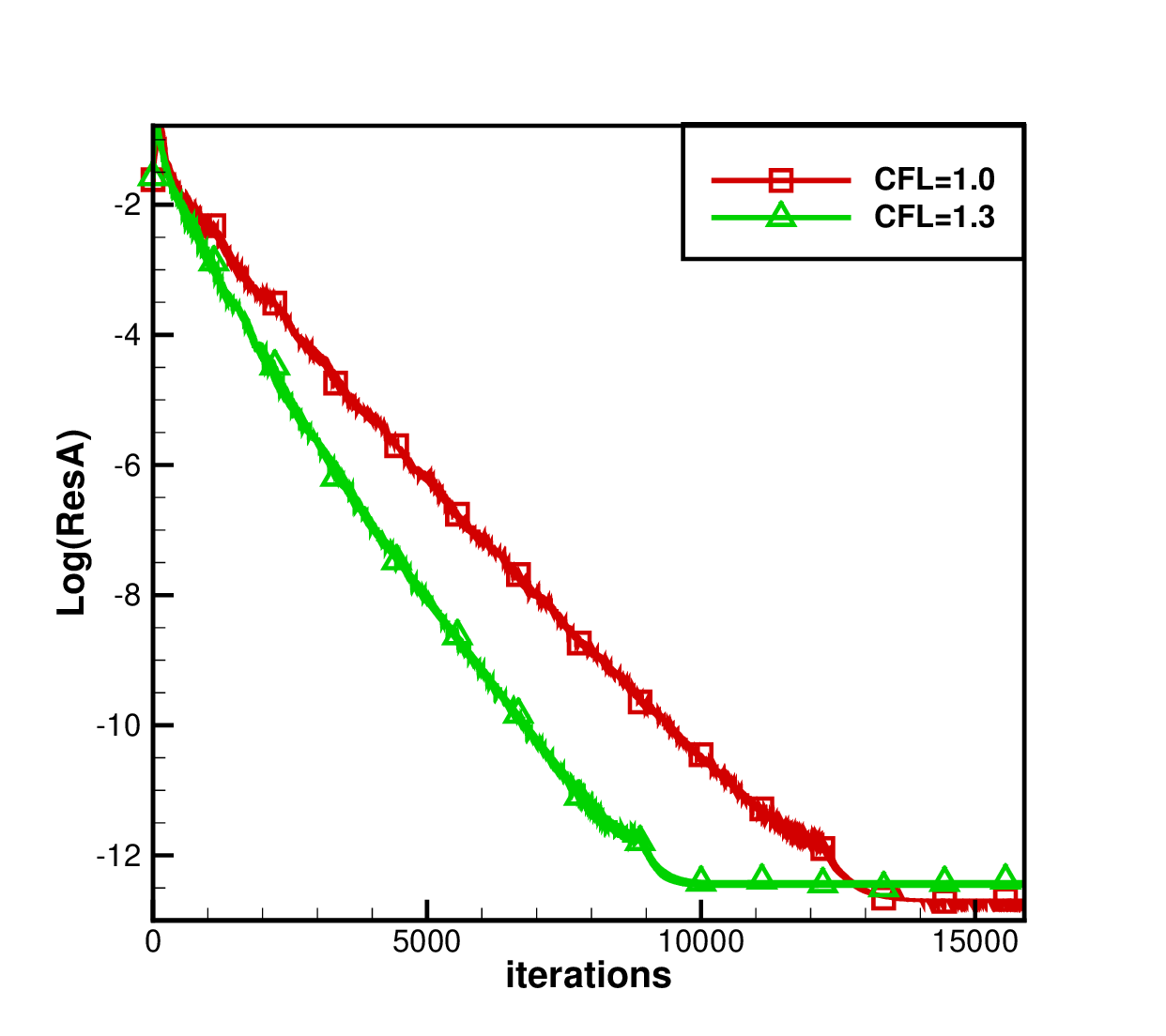}
\end{minipage}}
\caption{\label{8.3}Example 8: The convergence history of the residue as a function of number of iterations for four schemes.}
\end{figure}

\noindent Previously, we have computed many challenging numerical examples. In the end, let us test one more problem, where the FS-WENOJS scheme also achieves convergence of the iteration residues to machine zero. We will use this example to further verify the effectiveness of the FS-WENOJS-AC scheme. This test case\cite{LILW} involves a flow at Mach 0.5 impinging on a NACA0012 airfoil at an angle of attack of 1.25 degrees. The simulation domain spans $[-1, 2] \times [-1.5, 1.5]$, is resolved with a uniform $200\times{200}$ grid. There are no shock waves in the flow field, making it easier to converge. When computing this problem, the residuals of the FS-WENOJS scheme and the FS-MRWENO scheme can only decrease to about $10^{-11}$, as can be seen from Figure \ref{8.3}. When we set the convergence criterion to $10^{-11}$, the iteration counts and CPU times of the four schemes are shown in Table \ref{8.1}. When we set the convergence criterion to $10^{-12}$, only the FS-WENOJS-AC scheme and the RK-WENOJS-AC scheme are able to achieve convergence, and the iteration counts and CPU times of these two schemes are presented in Table \ref{8.11}. From these two tables, it can be seen that the FS-WENOJS-AC scheme still performs the best. Although its largest possible CFL number is only 1.3, much smaller than the 1.8 allowed in the RK-WENOJS-AC scheme, its shortest CPU time is still 40\% less than that of the RK-WENOJS-AC scheme. This test case also once again demonstrates that judging whether to freeze the linear weights does not consume excessive extra time, as can be seen from the time comparison between the FS-WENOJS-AC and FS-WENOJS schemes. Figure \ref{8.2} presents the pressure contour plots of these three fully converged schemes, and their solutions are comparable.

\begin{table}
		\centering
\begin{tabular}{|c|c|c|}\hline
			\multicolumn{3}{|c|}{FS-MRWENO}\\\hline
            $\gamma:$ CFL number & iteration number  & CPU time \\\hline
1.0	&10986		&2200.53\\\hline
1.4	&7202		&1442.58\\\hline
			\multicolumn{3}{|c|}{FS-WENOJS}\\\hline
            $\gamma:$ CFL number & iteration number & CPU time \\\hline
1.0	&10634		&1066.02\\\hline
1.3	&7689		&774.52\\\hline
			\multicolumn{3}{|c|}{RK-WENOJS-AC}\\\hline
            $\gamma:$ CFL number & iteration number  & CPU time \\\hline
1.0	&37842		&2765.47\\\hline
1.8	&21318		&1301.23\\\hline
			\multicolumn{3}{|c|}{FS-WENOJS-AC}\\\hline
            $\gamma:$ CFL number & iteration number  & CPU time \\\hline
1.0	&10634		&1089.52\\\hline
1.3	&7689		&775.89\\\hline
		\end{tabular}
		\caption{\label{8.1}Example 8: Number of iterations, the final time, and the total CPU time of four different iterative schemes when convergence is obtained. Convergence
criterion threshold value is $10^{-11}$. CPU time unit: second.}
	\end{table}

\begin{table}
		\centering
\begin{tabular}{|c|c|c|}\hline
			\multicolumn{3}{|c|}{RK-WENOJS-AC}\\\hline
            $\gamma:$ CFL number & iteration number  & CPU time \\\hline
1.0	&42981		&3048.78\\\hline
1.8	&23433		&1426.11\\\hline
			\multicolumn{3}{|c|}{FS-WENOJS-AC}\\\hline
            $\gamma:$ CFL number & iteration number  & CPU time \\\hline
1.0	&12364		&   1356.94  \\\hline
1.3	&9072		&   965.86  \\\hline
		\end{tabular}
		\caption{\label{8.11}Example 8: Number of iterations, the final time, and the total CPU time of two different iterative schemes when convergence is obtained. Convergence
criterion threshold value is $10^{-12}$. CPU time unit: second.}
	\end{table}

\section{Concluding remarks}
\label{sec3}
\setcounter{equation}{0}
\setcounter{figure}{0}
\setcounter{table}{0}
In this paper, in stead of using WENO local solvers with unequal-sized substencils to achieve fully convergence, we go back to the classical fifth-order WENO-JS local solver which has equal-sized substencils and develop a new fully convergent fixed-point fast sweeping method for solving steady-state problems of hyperbolic conservation laws. 
A simple and novel frozen-weight strategy is proposed. This simple approach eliminates the unnecessary adjustments of nonlinear weights near discontinuities, allowing the iteration residual to decrease to machine zero while preserving the order of accuracy and the essentially non-oscillatory property. 
Extensive numerical experiments on various two-dimensional steady flow problems demonstrate that, unlike the original fifth-order WENO-JS fast sweeping scheme, the new WENO-JS fast sweeping scheme consistently achieve full convergence for all examples. Furthermore, the new fast sweeping scheme exhibits excellent efficiency in the simulations. It is approximately twice faster in the computational costs to reach fully convergence, than the other fully converged fast sweeping methods such as the FS-MRWENO scheme and the RK-WENOJS-AC scheme. 


\begin{thebibliography}{99}

\bibitem{Bals1}
D. S. Balsara, S. Garain, V. Florinski and W. Boscheri. {\em An efficient class of WENO schemes with adaptive order for unstructured meshes}, J. Comput. Phys., 404: 109062, 2020.

\bibitem{BaZor}
A. Baeza, R. Bürger, P. Mulet and D. Zorío. {\em On the efficient computation of smoothness indicators for a class of WENO reconstructions}, J. Sci. Comput., 80(2): 1240--1263, 2019.

\bibitem{BorgesCarmona}
R. Borges, M. Carmona, B. Costa and W.S. Don.
{\em An improved weighted essentially non-oscillatory scheme for hyperbolic conservation laws}.
J. Comput. Phys., 227(6):3191--3211, 2008.



\bibitem{DK2}
M. Dumbser, M. K\"{a}ser, V.A. Titarev and E.F. Toro.
{\em Quadrature-free non-oscillatory finite volume schemes on unstructured
meshes for nonlinear hyperbolic systems},
J. Comput. Phys., 226:  204--243, 2007.

\bibitem{Engqui}
B. Engquist, B. D. Froese and Y.-H. R. Tsai. {\em Fast sweeping methods for hyperbolic systems of conservation laws at steady state}. J. Comput. Phys.,
255: 316--338, 2013.




\bibitem{RK1}
 S. Gottlieb, C.-W. Shu, E. Tadmor.
{\em Strong stability-preserving high-order time discretization methods}. SIAM Rev., 43(1): 89--112, 2001.

\bibitem{JYoon}
Y. Ha, C. H. Kim, Y. H. Yang and J. Yoon. {\em Improving accuracy of the fifth-order WENO scheme by using the exponential approximation space}, SIAM J. Numer. Anal. 59: 143--172, 2021.

\bibitem{HartenOsher}
A. Harten and S. Osher.
{\em Uniformly high-order accurate non-oscillatory schemes I}.
SIAM J. Numer. Anal., 24(2): 279--309, 1987.

\bibitem{HAPo}
A.K. Henrick, T.D. Aslam and J.M. Powers. {\em Mapped weighted essentially non-oscillatory
schemes: achieving optimal order near critical points}, J. Comput. Phys., 207: 542--567, 2005.


\bibitem{HS}
C. Hu and C.-W. Shu. {\em Weighted essentially non-oscillatory schemes on triangular meshes}, J. Comput. Phys., 150: 97--127, 1999.

\bibitem{CSHuang}
C.S. Huang, T. Arbogast and C. Tian. {\em Multidimensional WENO-AO reconstructions using a simplified smoothness indicator and applications to conservation laws}. J. Sci. Comput., 97: 8, 2023.

\bibitem{JS}
G.-S. Jiang and C.-W. Shu.
{\em Efficient implementation of weighted ENO schemes}.
J. Comput. Phys., 126(1): 202--228, 1996.



\bibitem{NACA}
K. Kuwaharaf and H. Takami.
{\em Computation of dynamic stall of a NACA-0012 airfoil}.
AIAA J., 25(3):408--413, 1987.

\bibitem{LILW}
L. Li, J. Zhu, C.-W. Shu and Y.-T. Zhang.
{\em A fixed-point fast sweeping WENO method with inverse Lax-Wendroff boundary treatment for steady state of hyperbolic conservation laws}.
Commun. Appl. Math. Comput., 5:403--427, 2023.

\bibitem{LZZ}
L. Li, J. Zhu and Y.-T. Zhang.
{\em Absolutely convergent fixed-point fast sweeping WENO methods for steady state of hyperbolic conservation laws}.
J. Comput. Phys., 443: 110516, 2021.



\bibitem{liangfu}
T. Liang, L. Fu. {\em A novel finite-difference converged ENO scheme for steady-state simulations of Euler equations}. J. Comput. Phys., 519:113386, 2024.



\bibitem{LiuOsherChen}
X.-D. Liu, S. Osher and T. Chan.
{\em Weighted essentially non-oscillatory schemes}.
J. Comput. Phys., 115(1): 200--212, 1994.

\bibitem{Y.Liu}
Y. Liu and Y.-T. Zhang. {\em A robust reconstruction for unstructured WENO schemes}, J.
Sci. Comput., 54: 603--621, 2013.

\bibitem{Lozano}
E. Lozano and T. D. Aslam. {\em Implicit fast sweeping method for hyperbolic systems of conservation laws}. J. Comput. Phys., 430: 110039, 2021.

\bibitem{ILW3}
J.-F. Lu, C.-W. Shu, S.-R. Tan and M.-P. Zhang.
{\em An inverse Lax-Wendroff procedure for hyperbolic conservation laws with changing wind direction on the boundary}.
J. Comput. Phys., 426:109940, 2021.

\bibitem{NACA2}
H. Luo, J.D. Baum, R. L{\"o}hner.
{\em On the computation of steady-state compressible flows using a discontinuous Galerkin method}.
Int. J. Numer. Methods Eng., 73:597--623, 2008.





\bibitem{RK2}
C.-W. Shu and S. Osher.
{\em Efficient implementation of essentially non-oscillatory shock-capturing schemes}.
J. Comput. Phys., 77(2): 439--471, 1988.

\bibitem{S.Tan}
S.-R. Tan and C.-W. Shu.
{\em Inverse Lax-Wendroff procedure for numerical boundary conditions of conservation laws}.
J. Comput. Phys., 229(21):8144--8166, 2010.

\bibitem{SIRUIT}
S.-R. Tan, C. Wang, C.-W. Shu and J.-G. Ning.
{\em Efficient implementation of high order inverse Lax-Wendroff boundary treatment for conservation laws}.
J. Comput. Phys., 231(6):2510--2527, 2012.

\bibitem{TsDum}
P. Tsoutsanis and M. Dumbser. {\em Arbitrary high order central non-oscillatory schemes on mixed-element unstructured meshes}, Computers \& Fluids, 225: 104961, 2021.

\bibitem{TsybulZhang}
E. Tsybulnik, X. Zhu and Y.-T. Zhang. {\em Efficient sparse-grid implementation of a fifth-order multi-resolution WENO scheme for hyperbolic equations}, Commun. Appl. Math.  Comput., 5(4): 1339--1364, 2023. 

\bibitem{WZ}
L. Wu and Y.-T. Zhang.
{\em A third order fast sweeping method with linear computational complexity for Eikonal equations}.
J. Sci. Comput., 62(1):198--229, 2015.

\bibitem{WuLiang}
L. Wu, Y.-T. Zhang, S. Zhang and C.-W. Shu.
{\em High order fixed-point sweeping WENO methods for steady state of hyperbolic conservation laws and its convergence study}.
Commun. Comput. Phys., 20(4):835--869, 2016.




\bibitem{XZZS}
T. Xiong, M. Zhang, Y.-T. Zhang, and C.-W. Shu.
{\em Fast sweeping fifth order WENO scheme for static Hamilton-Jacobi equations with accurate boundary treatment}.
J. Sci. Comput., 45(1):514--536, 2010.



\bibitem{SSCW}
S. Zhang, S. Jiang and C.-W. Shu.
{\em Improvement of convergence to steady state solutions of Euler equations with the WENO schemes}.
J. Sci. Comput., 47(2):216--238, 2011.

\bibitem{SCW}
S. Zhang and C.-W. Shu.
{\em A new smoothness indicator for the WENO schemes and its effect on the convergence to steady state solutions}.
J. Sci. Comput., 31(1):273--305, 2007.

\bibitem{ZCLS}
Y.-T. Zhang, S. Chen, F. Li, H. Zhao, and C.-W. Shu.
{\em Uniformly accurate discontinuous Galerkin fast sweeping methods for Eikonal equations}.
SIAM J. Sci. Comput., 33(4):1873--1896, 2011.

\bibitem{ZS}
Y.-T. Zhang and C.-W. Shu, Third order WENO scheme on three dimensional tetrahedral meshes, Commun. Comput. Phys., 5: 836--848, 2009.

\bibitem{ZZC}
Y.-T. Zhang, H.-K. Zhao and S. Chen.
{\em Fixed-point iterative sweeping methods for static Hamilton-Jacobi equations}.
Methods Appl. Anal., 13:299--320, 2006.

\bibitem{ZZQ}
Y.-T. Zhang, H.-K. Zhao, and J. Qian.
{\em High order fast sweeping methods for static Hamilton-Jacobi equations}.
J. Sci. Comput., 29(1):25--56, 2006.

\bibitem{FW}
Z. Zhang, Y. Chen and X. Deng.
{\em On steady-state solutions of nonlinear shock-capturing schemes}, J. Comput. Phys., 556:114820, 2026.

\bibitem{ZE}
H.-K. Zhao.
{\em A fast sweeping method for Eikonal equations}.
Math. Comput., 74(250):603--627, 2005.

\bibitem{ZQ3}
J. Zhu and J. Qiu. {\em New finite volume weighted essentially non-oscillatory schemes on triangular meshes}, SIAM J. Sci. Comput., 40: 903--928, 2018.


\bibitem{taijie2}
J. Zhu and C.-W. Shu.
{\em Numerical study on the convergence to steady state solutions of a new class of high order WENO schemes}.
J. Comput. Phys., 349(5):80--96, 2017.

\bibitem{JunZhu2}
J. Zhu and C.-W. Shu.
{\em A new type of multi-resolution WENO schemes with increasingly higher order of accuracy}.
J. Comput. Phys., 375(3):659--683, 2018.

\bibitem{m2}
J. Zhu and C.-W. Shu.
{\em Numerical study on the convergence to steady-state solutions of a new class of finite volume WENO schemes: triangular meshes}.
Shock Waves, 29(1):3--25, 2018.

\bibitem{JUNZ3}
J. Zhu and C.-W. Shu.
{\em Convergence to steady-state solutions of the new type of high-order multi-resolution WENO schemes: a numerical study}.
Commun. Appl. Math. Comput., 2(6):429--460, 2020.


\end{thebibliography}
\end{document}